\documentclass[a4paper,oneside,english,american]{amsart}
\usepackage[T1]{fontenc}
\usepackage{color}
\usepackage{amstext}
\usepackage{amsthm}
\usepackage{amssymb}

\makeatletter

\providecommand{\tabularnewline}{\\}

\numberwithin{equation}{section}
\theoremstyle{plain}
\newtheorem{thm}{\protect\theoremname}[section]
\theoremstyle{plain}
\newtheorem{prop}[thm]{\protect\propositionname}
\theoremstyle{remark}
\newtheorem{rem}[thm]{\protect\remarkname}
\theoremstyle{definition}
\newtheorem{defn}[thm]{\protect\definitionname}
\theoremstyle{plain}
\newtheorem{lem}[thm]{\protect\lemmaname}
\theoremstyle{plain}
\newtheorem{conjecture}[thm]{\protect\conjecturename}

\newcounter{myparagraph}[subsection]

\renewcommand{\themyparagraph}{\bf {\arabic{section}.\arabic{subsection}.\alph{myparagraph}}}

\newenvironment{myitem}{\begin{list}{}{
\setlength{\leftmargin}{0.8cm}
\setlength{\itemindent}{-0.5cm}
\setlength{\itemsep}{2pt}
}}{\end{list}}

\def\hsp#1{{\hspace{ #1 pt}}}

\def\*{\hsp{-3p}*\hsp{-3pt}}

\usepackage[all]{xy}

\usepackage{color}
\usepackage{amsthm}
\usepackage{bm}
\usepackage{graphicx, color, ulem}

\usepackage[labelformat=empty]{caption}

\makeatother

\usepackage{babel}
\addto\captionsamerican{\renewcommand{\conjecturename}{Conjecture}}
\addto\captionsamerican{\renewcommand{\definitionname}{Definition}}
\addto\captionsamerican{\renewcommand{\lemmaname}{Lemma}}
\addto\captionsamerican{\renewcommand{\propositionname}{Proposition}}
\addto\captionsamerican{\renewcommand{\remarkname}{Remark}}
\addto\captionsamerican{\renewcommand{\theoremname}{Theorem}}
\addto\captionsenglish{\renewcommand{\conjecturename}{Conjecture}}
\addto\captionsenglish{\renewcommand{\definitionname}{Definition}}
\addto\captionsenglish{\renewcommand{\lemmaname}{Lemma}}
\addto\captionsenglish{\renewcommand{\propositionname}{Proposition}}
\addto\captionsenglish{\renewcommand{\remarkname}{Remark}}
\addto\captionsenglish{\renewcommand{\theoremname}{Theorem}}
\providecommand{\conjecturename}{Conjecture}
\providecommand{\definitionname}{Definition}
\providecommand{\lemmaname}{Lemma}
\providecommand{\propositionname}{Proposition}
\providecommand{\remarkname}{Remark}
\providecommand{\theoremname}{Theorem}

\begin{document}
\selectlanguage{english}%
\global\long\def\ve{\mathcal{\varepsilon}}%
\global\long\def\st{\;\mathcal{\text{s.t.}}\;}%
\global\long\def\vp{\mathcal{\varphi}}%
\global\long\def\c{\circ}%

\global\long\def\cO{\mathcal{O}}%

\global\long\def\cA{\mathcal{A}}%

\global\long\def\cB{\mathcal{B}}%

\global\long\def\cD{\mathcal{D}}%

\global\long\def\cE{\mathcal{E}}%

\global\long\def\cF{\mathcal{F}}%

\global\long\def\cM{\mathfrak{M}}%

\global\long\def\cU{\mathcal{U}}%
\global\long\def\cV{\mathcal{V}}%
\global\long\def\cW{\mathcal{W}}%

\global\long\def\cN{\mathcal{N}}%
\global\long\def\cH{\mathcal{H}}%

\global\long\def\idR{1_{\mathbb{R}^{n}}}%

\global\long\def\id{\mathrm{id}}%

\global\long\def\bR{\mathbb{R}}%

\global\long\def\mbR{\mathbb{R}}%
\global\long\def\mbC{\mathbb{C}}%

\global\long\def\mbP{\mathbb{P}}%
\global\long\def\mbH{\mathbb{H}}%

\global\long\def\mbZ{\mathbb{Z}}%

\global\long\def\bS{\bar{S}}%
\global\long\def\bT{\bar{T}}%

\global\long\def\tx{\mathtt{x}}%
\global\long\def\ty{\mathtt{y}}%
\global\long\def\tz{\mathtt{z}}%
\global\long\def\tw{\mathtt{w}}%

\selectlanguage{american}%
\def\m{\hsp{-4}}

\def\n{\hsp{-2}}

\selectlanguage{english}%
\global\long\def\s{\mathtt{\n*\n}}%

\selectlanguage{american}%
~
\begin{flushright}
~%\today
\par\end{flushright}
\title{BCOV cusp forms of lattice polarized K3 surfaces}
\author{Shinobu Hosono and Atsushi Kanazawa}
\begin{abstract}
We introduce the BCOV formula for the lattice polarized K3 surfaces.
We find that it yields cusp forms expressed by certain eta products
for many families of rank 19 lattice polarized K3 surfaces over $\mathbb{P}^{1}$.
Moreover, for Clingher-Doran's family of $U\oplus E_{8}(-1)\oplus E_{7}(-1)$-polarized
K3 surfaces, we obtain the Igusa cusp forms $\chi_{10}$ and $\chi_{12}$
from the formula. Inspired by the arithmetic properties of mirror
maps studied by Lian-Yau, we also derive the K3 differential operators
for all the genus zero groups of type $\Gamma_{0}(n)_{+}$. 
\end{abstract}

\maketitle

\section{\textbf{\label{sec:Introduction}Introduction}}

In 1993, Bershadsky, Ceccoti, Ooguri and Vafa (BCOV) introduced the
so-called BCOV potential $F_{1}$ which gives the generating functions
of genus one Gromov-Witten invariants of Calabi-Yau threefolds generalizing
mirror symmetry at genus zero due to Candelas et al \cite{CD1}. It
involves an interesting interplay between the variation of Hodge structure
and the K\"ahler geometry on moduli spaces of Calabi-Yau manifolds.
In fact, it originates from the so-called $tt^{*}$-geometry \cite{Cecotti-Vafa}
on the moduli space of $N=2$ supersymmetric quantum field theories
in two dimensions in theoretical physics. 

One of the main ingredients in the $tt^{*}$-geometry is \textit{a
new} type of Witten index $\mathrm{F}_{1}$, introduced in \cite{Cecotti-FI-Vafa,Cecotti-Vafa},
for $N=2$ supersymmetric conformal field theories. In contrast to
the Witten index, the new index $\mathrm{F}_{1}$ is not topological
but depends on the moduli of a theory in a specific way; as a function
on the moduli space, $\mathrm{F}_{1}$ splits almost to a sum of holomorphic
and anti-holomorphic functions, but the splitting is not complete.
When the $N=2$ supersymmetric theory is given as a sigma model over
a Calabi-Yau manifold $X$, it satisfies the \textit{holomorphic anomaly
equation} 
\begin{equation}
\frac{\partial\;}{\partial x_{i}}\frac{\partial\;}{\partial\overline{x}_{j}}\mathrm{F}_{1}=\frac{1}{2}\mathrm{Tr}\,\mathcal{C}_{i}\overline{\mathcal{C}}_{\bar{j}}+\frac{\chi}{24}\,g_{i\bar{j}},\label{eq:BCOV-tt*-Tr}
\end{equation}
where the trace is taken over the space of the quantum operators whose
ring structure is determined by $\mathcal{C}_{i}=\left(\mathcal{C}_{ijk}\right)$
and $\overline{\mathcal{C}}_{\bar{j}}=\left(\overline{\mathcal{C}}_{\bar{j}\bar{k}\bar{l}}\right)$.
Also, $\chi$ is the Euler number of $X$ and $g_{i\bar{j}}=\frac{\partial\;}{\partial x_{i}}\frac{\partial\;}{\partial\overline{x}_{j}}\mathcal{K}$
is the Weil-Petersson K\"ahler metric on the moduli space (of mirror
Calabi-Yau manifold $\check{X}$). Identifying the linear space of
the quantum operators with the tangent space of the moduli space,
and using the so-called special K\"ahler geometry relation \cite{Str}
for a Calabi-Yau threefold $X$, BCOV integrate the equation (\ref{eq:BCOV-tt*-Tr})
as 
\begin{equation}
\mathrm{F}_{1}=\frac{1}{2}\log\left\{ e^{(3+h_{X}^{1,1}-\frac{\chi}{12})\mathcal{K}}\det(g_{i\bar{j}})^{-1}\,\vert f\vert^{2}\right\} ,\label{eq:BCOV-tt*}
\end{equation}
where $h_{X}^{1,1}=\dim H^{1,1}(X)$ and $f$ is a holomorphic function
which should be determined by suitable boundary conditions. As we
shall describe in detail in the text, there are special boundary points
called the large complex structure limits (LCSLs), where we introduce
special holomorphic coordinates $t_{i}=t_{i}(x)$ $(1\leq i\leq d:=h_{X}^{1,1})$
called mirror maps. In terms of the holomorphic coordinates $t_{i}$,
BCOV take the so-called topological limit $\bar{t}_{i}\rightarrow\infty$
where $\mathrm{F}_{1}$ decomposes to a sum of holomorphic and anti-holomorphic
functions. Precisely, under the limit, we have $\mathcal{K}\rightarrow-\log\big(w_{0}(x)\overline{w_{0}(x)}\big)$,
$\det(g_{i\bar{j}})^{-1}\rightarrow\vert\frac{\partial(x_{1},\cdots,x_{d})}{\partial(t_{1},\cdots,t_{d})}\vert^{2}$,
and using these relations, we obtain 
\begin{equation}
F_{1}^{top}=\frac{1}{2}\log\left\{ \Big(\frac{1}{w_{0}(x)}\Big)^{3+h_{X}^{1,1}-\frac{\chi}{12}}\frac{\partial(x_{1},\cdots,x_{d})}{\partial(t_{1},\cdots,t_{d})}f(x)\right\} \label{eq:BCOV-CY3}
\end{equation}
for the holomorphic part of the limit. The holomorphic function $F_{1}^{top}$
thus obtained is nothing but the BCOV potential $F_{1}$ when we determine
$f(x)$ by imposing suitable boundary conditions. Writing $F_{1}^{top}=\frac{1}{2}\log\tau_{_{\mathrm{BCOV}}}$,
we call $\tau_{_{\mathrm{BCOV}}}$ the BCOV formula. 

The holomorphic function $F_{1}^{top}$ is particularly interesting
for Calabi-Yau threefolds as above. Considering $F_{1}^{top}$ for
K3 surfaces, however,  has not attracted much attention, since there
is no quantum corrections coming from holomorphic curves for K3 surfaces.
In this paper, we will introduce the BCOV formula for K3 surfaces
and show that it is of considerable interest nevertheless, since it
yields interesting automorphic forms on the period domains.

~

Let $M$ be a lattice of signature $(1,r-1)$ and $M\subset L$ be
a primitive embedding to the K3 lattice $L=U^{\oplus3}\oplus E_{8}(-1)^{\oplus2}$.
Assume the orthogonal lattice $M^{\perp}=U\oplus\check{M}$ splits
a hyperbolic lattice $U$ of rank two. Under this assumption, an $M$-polarized
K3 surface $X$ is mirror symmetric to a family of $\check{M}$-polarized
K3 surfaces which are parametrized by the period domain $\Omega_{\check{M}}$.
Studying the isomorphism classes in more detail, we introduce $\check{M}$-polarizable
K3 surfaces whose isomorphism classes are parametrized by a certain
quotient $\Omega_{\check{M}}/O(\check{M}^{\perp})_{+}$ of the period
domain. Introducing the mirror map $t_{i}=t_{i}(x)$$(1\leq i\leq r)$,
we define the BCOV formula (in Definition \ref{def:tau-BCOV}) by
\[
\tau_{_{\mathrm{BCOV}}}(t)=\left(\frac{1}{w_{0}(x)}\right)^{r+1}\frac{\partial(x_{1},\cdots,x_{r})}{\partial(t_{1},\cdots,t_{r})}\prod_{i}dis_{i}(x)^{r_{i}}\,\prod_{j=1}^{r}x_{j}^{-1+a_{j}},
\]
where we impose suitable boundary conditions (\textit{conifold} and
\textit{orbifold regularity} conditions) to determine the parameters
$r_{i},a_{j}$. Note that this definition naturally follows from (\ref{eq:BCOV-CY3})
by setting $\chi=24$ and $h_{X}^{1,1}=r$ and making an ansatz for
$f$ as above. The K3 surfaces of degree $2n$ are described by a
lattice $M_{2n}:=\langle2n\rangle$, which gives rise to $\check{M}_{2n}=\langle-2n\rangle\oplus U\oplus E_{8}(-1)^{\oplus2}$.
We calculate the BCOV formula for the families of $\check{M}_{2n}$-polarizable
K3 surfaces over $\mathbb{P}^{1}$  studied in \cite{LY}. It turns
out in \textbf{Propositions \ref{prop:tauGCOV-10}, \ref{prop:tauBCOV-12}}
that $(\tau_{_{\mathrm{BCOV}}})^{-1}$ defines a cusp form on the
upper half plane $\mathbb{H}_{+}$ with the modular group $\Gamma_{0}(n)_{+}$,
which has the following special form of eta product (Definition \ref{def:tauBCOV-def-2}):
\[
\eta_{_{\mathrm{BCOV}}}(t)=\left(\prod_{r|n}\eta_{r}(t)^{\pm1}\right)^{w}.
\]
We conjecture that this is true in general for families of $\check{M}_{2n}$-polarizable
K3 surfaces whose isomorphism classes are parametrized by $\Omega_{\check{M}_{2n}}/\Gamma_{0}(n)_{+}$
with the genus zero modular group $\Gamma_{0}(n)_{+}$ (\textbf{Conjecture
\ref{conj:tauBCOV}}). 

Interestingly, if we postulate Conjecture \ref{conj:tauBCOV} and
the arithmetic properties of mirror maps observed by Lian-Yau, a certain
third order differential operator, which we call an\textit{ K3 differential
operator}, follows for each genus zero group $\Gamma_{0}(n)_{+}$.
We summarize properties of K3 differential operators in \textbf{Proposition
\ref{prop:K3-diff-ops}}. A list of K3 differential operators is provided
in Addendum of this paper \cite{addendum}. 

As another interesting example, we calculate the BCOV formula for
the family of $U\oplus E_{8}(-1)\oplus E_{7}(-1)$-polarized K3 surfaces
studied extensively by Clingher and Doran \cite{CD1}, where a coarse
moduli space is described in terms of the Siegel modular forms of
genus two. After a suitable change of coordinates, we find that the
Picard-Fuchs differential equations are given by the extended GKZ
system introduced in the 90s \cite{HKTY1,HLY}. In this case, it turns
out in \textbf{Proposition \ref{prop:BCOV-cusp-Chi10}} that we have
\[
(\tau_{_{\mathrm{BCOV}}}(\tau))^{-1}=\left(\chi_{10}(\tau)\right)^{\frac{1}{10}},\,\,\left(3\chi_{12}(\tau)+\chi_{10}(\tau)\mathcal{E}_{4}(\tau)^{\frac{1}{2}}\right)^{\frac{1}{12}}
\]
where $\chi_{10},\chi_{12}$ are the Igusa cusp forms of weight ten
and twelve, respectively, and $\mathcal{E}_{4}$ is the genus two
Eisenstein series. Here two possibilities arise from two different
boundary conditions (orbifold regularity conditions); one of which
we will interpret in terms of the extension of the lattice polarization
$\check{M}$ to $U\oplus E_{8}(-1)^{\oplus2}$. 

We describe the Weil-Petersson K\"ahler geometry on the moduli space
of Calabi-Yau manifolds. Introducing the so-called BCOV torsion $\mathcal{T}_{\mathrm{BCOV}}$
in terms of $\tau_{_{\mathrm{BCOV}}}$, we present some useful formulas
which connect the K\"ahler geometry to (quasi-)automorphic forms
on $\Omega_{\check{M}}$.

~

We briefly sketch the construction of this paper. In Sect. 2, we summarize
the theory of lattice polarized K3 surfaces and mirror symmetry among
them. We describe the period domain and automorphic forms on the domain
in detail to set up notation. In Sect. 3, we introduce the BCOV formula
$\tau_{_{\mathrm{BCOV}}}$ and calculate it for the known $\check{M}_{2n}$-polarized
K3 surfaces over $\mathbb{P}^{1}$ in \cite{LY}, and observe that
they are written by eta products of the specific form $\eta_{_{\mathrm{BCOV}}}$.
Postulating the generality of the results, we obtain K3 differential
operators. In Sect. 4, we calculate $\tau_{_{\mathrm{BCOV}}}$ for
Clingher-Doran's family of K3 surfaces. In Sect.5, we introduce the
BCOV torsion and express its relation to Weil-Petersson K\"ahler
geometry on the moduli spaces of lattice polarized K3 surfaces. We
describe detailed calculations and data which should be useful for
the reader in Appendices A--D. 

~

\noindent\textbf{ Acknowledgements. }This work came out from a collaboration
with Martin Guest on the $tt^{*}$-geometry on moduli spaces of K3
surfaces. We are grateful to him for his collaboration at early stage
and valuable discussions. This work is supported in part by Grants-in-Aid
for Scientific Research (A 18H03668 M.G., C 20K03593 S.H., and C 22K03296
A.K.).

\vskip2cm

\section{\textbf{Automorphic forms on period domains~}}

We first describe families of lattice polarized K3 surfaces and their
period domains. Then we introduce automorphic forms on the domain
by the tube domain realization.

\subsection{\label{subsec:K3-period-domain}Families of K3 surfaces and period
domains}

Let $L$ be the so-called K3 lattice which is (a unique) even unimodular
lattice of signature $(3,19)$, i.e., $L=U^{\oplus3}\oplus E_{8}(-1)^{\oplus2}$,
where $U=\left\langle \mbZ^{2},\left(\begin{smallmatrix}0 & 1\\
1 & 0
\end{smallmatrix}\right)\right\rangle $ represents the hyperbolic lattice and $E_{8}(-1)$ is the $E_{8}$
root lattice with its bilinear form being multiplied by $(-1)$. We
fix a primitive sub-lattice $M\subset L$ of signature $(1,r-1)$,
and denote by $C_{M}$ one of the (two) connected components of the
positive cone defined for $M$ in $M\otimes\mbR$. 

\subsubsection{Lattice polarized K3 surfaces}

A K3 surface $X$ is a compact 2-dimensional complex manifold satisfying
$K_{X}=0$ and $H^{1}(X,\cO_{X})=0$. By the Poincar\'e duality,
the bilinear form $(\,,\,):H^{2}(X,\mbZ)\times H^{2}(X,\mbZ)\rightarrow\mbZ$
via the cup product is non-degenerate, and $(H^{2}(X,\mbZ),(\,,\,))$
is a lattice which is isomorphic to the K3 lattice $L$. A lattice
$NS_{X}:=H^{1,1}(X)\cap H^{2}(X,\mbZ)$ is called the N\'eron-Severi
lattice, and its orthogonal lattice $T_{X}:=(NS_{X})^{\perp}$ in
$H^{2}(X,\mbZ)$ is called the transcendental lattice. We call an
isomorphism $\phi:H^{2}(X,\mbZ)\rightarrow L$ a marking of $X$,
and a pair $(X,\phi)$ a marked K3 surface. We summarize below theory
of lattice polarized K3 surfaces in order.

\vskip0.3cm\noindent(2.1.1.a) Let us consider a lattice $M$ of signature
$(1,r-1)$ and a primitive embedding $i:M\hookrightarrow L$. We write
the orthogonal decomposition $i(M)\oplus(i(M))^{\perp}$ in $L$ by
$M\oplus M^{\perp}\subset L$ with suppressing the inclusion. A K3
surface $X$ is called of type $M$ if there exists a marking $\phi$
such that $\phi^{-1}(M)\subset NS_{X}$. A marked K3 surface $(X,\phi)$
is called $M$-polarized if $X$ is of type $M$ and all divisors
in $\phi^{-1}(C_{M}^{\text{pol}})$ are ample, where $C_{M}^{\text{pol}}$
is the positive chamber in $C_{M}$ determined by the roots $\Delta_{M}=\left\{ d\in M\mid(d,d)=-2\right\} $
(see \cite{Dol,BHPvDeV}).

\vskip0.3cm\noindent(2.1.1.b)\textbf{ }Associated to the above orthogonal
decomposition $M\oplus M^{\perp}\subset L$, we define the period
domain $\Omega_{M}=\Omega(M^{\perp})$ as one of the (two) connected
components of 
\[
D_{M}=\left\{ [w]\in\mbP(M^{\perp}\otimes\mbC)\mid(w,w)=0,(w,\bar{w})>0\right\} ,
\]
i.e., $\Omega_{M}=D_{M}^{+}$ with $D_{M}=D_{M}^{+}\sqcup D_{M}^{-}$.
It is clear that we have $[\phi(\omega_{X})]\in\Omega_{M}$ for holomorphic
two forms $\omega_{X}$ of $M$-polarized K3 surfaces $(X,\phi)$.
Conversely, by the surjectivity of period map, for a point $[\omega]\in\Omega_{M}$
which we regard as a point in $\mathbb{P}(M^{\perp}\otimes\mathbb{C})\subset\mathbb{P}(L\otimes\mathbb{C})$,
there is a marked K3 surface $(X,\phi)$ such that $\phi(\mathbb{C}\omega_{X})=[\omega]$.
Note that we have $NS_{X}=(\mathbb{C}\omega_{X})^{\perp}\cap H^{2}(X,\mathbb{Z})$,
hence $\phi(NS_{X})=(\mathbb{C}\omega)^{\perp}\cap L.$ Also, using
the decomposition $M\oplus M^{\perp}\subset L$ and the fact that
$\phi(\omega_{X})\in M^{\perp}\otimes\mathbb{C}$, we see that $M\subset\phi(NS_{X})$,
i.e., $\phi^{-1}(M)\subset NS_{X}$. Composing suitable reflections
by roots in $M$ (and also in $M^{\perp}\cap(\mathbb{C}\omega_{X})^{\perp}$)
with $\phi$, we can make all divisors in $\phi^{-1}(C_{M}^{\text{pol}})$
be (pseudo-)ample. Hence, we have an (pseudo-ample) $M$-polarized
K3 surface for any point $[\omega]\in\Omega_{M}$. 

\vskip0.3cm\noindent(2.1.1.c)\textbf{ }We denote by $O(M,L)$ the
subgroup of the orthogonal group $O(L)$ of $L$ which fixes every
element in $M\subset L$. Using the decomposition $M\oplus M^{\perp}\subset L$,
the group $O(M,L)$ naturally acts on $D_{M}$. We denote by $O(M,L)_{+}$
the subgroup which fixes the connected component $D_{M}^{+}$. This
group $O(M,L)_{+}$ acts on $\Omega_{M}=D_{M}^{+}$ properly and discontinuously.
Two $M$-polarized K3 surfaces $(X_{1},\phi_{1})$ and $(X_{2},\phi_{2})$
are defined to be isomorphic if there is an isomorphism $f:X_{1}\rightarrow X_{2}$
such that $\phi_{1}\circ f^{*}\circ\phi_{2}^{-1}\in O(M,L)_{+}$.
If $(X_{1},\phi_{1})$ and $(X_{2},\phi_{2})$ are isomorphic, then
$f^{*}:H^{2}(X_{2},\mbZ)\stackrel{\sim}{\rightarrow}H^{2}(X_{1},\mbZ)$
satisfies $f^{*}(\mbC\omega_{X_{2}})=\mbC\omega_{X_{1}}$ for holomorphic
two forms $\omega_{X_{i}}$ of $X_{i}$ $(i=1,2)$. Correspondingly,
the element $\phi_{1}\circ f^{*}\circ\phi_{2}^{-1}\in O(M,L)_{+}$
maps the period point $[\phi_{2}(\omega_{X_{2}})]\in\Omega_{M}$ to
$[\phi_{1}(\omega_{X_{1}})]\in\Omega_{M}$. Hence, isomorphism classes
of (pseudo-ample) $M$-polarized K3 surfaces are parametrized by the
quotient $\Omega_{M}/O(L,M)_{+}$. The moduli space of isomorphism
classes of $M$-polarized K3 surfaces is obtained by removing hyperplanes
$H_{d}=\left\{ (x,d)=0\right\} $ for $d\in\Delta_{M^{\perp}}$ from
$\Omega_{M}$, where we can only make divisors in $\phi^{-1}(C_{M}^{\text{pol}})$
be pseudo-ample by suitable reflections by roots (see \cite{Dol,BHPvDeV}
and (\ref{eq:tubeD}) below).

\vskip0.2cm\noindent(2.1.1.d)\textbf{ }Let $O(M^{\perp})_{+}$ be
the subgroup of $O(M^{\perp})$ which fixes the connected component
$\Omega_{M}=\Omega(M^{\perp})$. Since $M^{\perp}\oplus M$ is a primitive
sub-lattice of $L$, which is even and unimodular, we have $O(M,L)_{+}\simeq O(M^{\perp})_{+}\cap O(M^{\perp})^{*}$
where
\[
O(M^{\perp})^{*}=\mathrm{Ker}\,\left\{ O(M^{\perp})\rightarrow O(A_{M^{\perp}})\right\} 
\]
and $O(A_{M^{\perp}})$ is the orthogonal group on the discriminant
$A_{M^{\perp}}=(M^{\perp})^{\vee}/M^{\perp}$ with respect to the
discriminant form $q_{A_{M^{\perp}}}$ (see e.g. \cite[Prop.3.3]{Dol})
. 

\vskip0.3cm\noindent(2.1.1.e) Mirror symmetry of K3 surfaces is described
nicely in terms of $M$-polarized K3 surfaces when the orthogonal
lattice $M^{\perp}$ of $M\subset L$ splits a hyperbolic lattice
$U$, i.e., $M^{\perp}=U\oplus\check{M}$ with a lattice $\check{M}$
of signature $(1,19-r)=:(1,\check{r}$). In fact, in such cases, we
can observe obvious symmetry in the relation $M\oplus U\oplus\check{M}\subset L$
with $r+\check{r}=20$. Describing suitable deformation spaces associated
to $M$ and $\check{M}$ precisely, it was argued in \cite{Dol} that
a family of $\check{M}$-polarized K3 surfaces is mirror symmetric
to the family of $M$-polarized K3 surfaces. 

\subsubsection{Lattice polarizable K3 surface}

Derived categories of K3 surfaces come into mirror symmetry when we
extend it to the so-called homological mirror symmetry (see \cite{Ko}).
It was shown in \cite{HLOY1} that homological mirror symmetry of
K3 surfaces is better described by introducing a family of $M$-\textit{polarizable}
K3 surfaces (see Subsect.~\ref{subsec:K3-2n-monod} below). A K3
surface $X$ is called $M$-polarizable if there is a marking $\phi:H^{2}(X,\mbZ)\rightarrow L$
such that $(X,\phi)$ is a $M$-polarized K3 surface.

\vskip0.3cm

If $M$-polarizable K3 surfaces $X_{1}$ and $X_{2}$ are generic,
i.e., satisfy $\phi_{i}^{-1}(M)=NS_{X_{i}}\,(i=1,2)$, then there
are corresponding $M$-polarized K3 surfaces $(X_{1},\phi_{1})$ and
$(X_{2},\phi_{2})$. When there is an isomorphism $f:X_{1}\stackrel{\sim}{\rightarrow}X_{2}$,
we have the corresponding Hodge isometry $f^{*}:(H^{2}(X_{2},\mbC),\mbC\omega_{X_{2}})\stackrel{\sim}{\rightarrow}(H^{2}(X_{1},\mbC),\mbC\omega_{X_{1}})$
which entails $\phi_{1}\circ f^{*}\circ\phi_{2}^{-1}:(L,\mbC\phi_{2}(\omega_{X_{2}}))\stackrel{\sim}{\rightarrow}(L,\mbC\phi_{1}(\omega_{X_{1}}))$.
Under the assumption of genericity, we have $\phi_{i}(NS_{X_{i}})\oplus\phi_{i}(T_{X_{i}})=M\oplus M^{\perp}\subset L$
and the restriction $\phi_{1}\circ f^{*}\circ\phi_{2}^{-1}\vert_{M^{\perp}}$
defines an element of $O(M^{\perp})_{+}$. The induced action of $\phi_{1}\circ f^{*}\circ\phi_{2}^{-1}\vert_{M^{\perp}}$
on the period domain $\Omega_{M}$ maps the period point $[\phi_{2}(\omega_{X_{2}})]$
to $[\phi_{1}(\omega_{X_{1}})]$. However, the restriction $\phi_{1}\circ f^{*}\circ\phi_{2}^{-1}\vert_{M}$
need not to be the identity $id_{M}$ in general. 
\begin{prop}
\label{prop:Primitive-Hodge} Let $g\cdot\omega$ be the natural action
of $g\in O(M^{\perp})_{+}$ on a point $[\omega]=\mathbb{C}\omega$
in $\Omega_{M}=\Omega(M^{\perp})\subset\Omega(L)$. Assume that a
primitive embedding $M^{\perp}\hookrightarrow L$ is unique up to
isometry of $L$, then every element $g\in O(M^{\perp})_{+}$ extends
to $\phi_{g}\in O(L)$ which defines a ``Hodge'' isometry $\phi_{g}:(L,\mbC\omega)\stackrel{\sim}{\rightarrow}(L,\mbC g\cdot\omega)$. 
\end{prop}

\begin{proof}
Recall that we have fixed our orthogonal decomposition $M\oplus M^{\perp}\subset L$
for a primitive embedding $i:M\hookrightarrow L$. This gives us a
primitive embedding $j:M^{\perp}\hookrightarrow L$. For $g\in O(M^{\perp})_{+}$,
we have a primitive embedding $j\circ g:M^{\perp}\hookrightarrow L$.
By the assumption of uniqueness, we then have an isometry $\phi_{g}\in O(L)$
such that $j\circ g=\phi_{g}\circ j$. Then, omitting $j$ for the
embedding $j:M^{\perp}\hookrightarrow L$, we obtain $\phi_{g}:(L,\mbC\omega)\stackrel{\sim}{\rightarrow}(L,\mbC g\cdot\omega)$
for any $\mathbb{C}\omega\in\Omega(M^{\perp})\subset\Omega(L)$.
\end{proof}
As described above, by the surjectivity of period map, we have $M$-polarized
K3 surfaces $(X_{1},\phi_{1})$ and $(X_{2},\phi_{2})$ satisfying
$[\phi_{1}(\omega_{X_{1}})]=[\omega]$ and $[\phi_{2}(\omega_{X_{2}})]=[g\cdot\omega]$
for $[\omega],[g\cdot\omega]\in\Omega_{M}$. Then, by the Torelli
Theorem, the ``Hodge'' isometry $\phi_{g}$ in the above proposition
gives a Hodge isometry $f_{g}^{*}=\phi_{1}^{-1}\circ\phi_{g}\circ\phi_{2}:(H^{2}(X_{2},\mbZ),\mbC\omega_{X_{2}})\stackrel{\sim}{\rightarrow}(H^{2}(X_{1},\mbZ),\mbC\omega_{X_{1}})$
which comes from an isomorphism $f_{g}:X_{1}\stackrel{\sim}{\rightarrow}X_{2}$
of K3 surfaces. Hence, under the assumption in Proposition \ref{prop:Primitive-Hodge},
the isomorphism classes of $M$-polarizable K3 surfaces are parametrized
by the quotient 
\[
\Omega_{M}/O(M^{\perp})_{+}.
\]
Here, it should be noted that a larger group $O(M^{\perp})_{+}$ than
$O(M^{\perp})_{+}\cap O(M^{\perp})^{*}\simeq O(M,L)_{+}^{*}$ naturally
appears to describe the moduli space if the above assumption holds. 

\subsubsection{\label{subsec:K3-2n-monod}Mirror to K3 surfaces of degree $2n$}

The K3 surfaces of degree $2n$ are the most general algebraic K3
surfaces described by a lattice $M_{2n}=\langle2n\rangle$. Associated
to a primitive embedding $M_{2n}\subset L$, we have the corresponding
orthogonal decomposition: $M_{2n}\oplus(M_{2n})^{\perp}=M_{2n}\oplus U\oplus\check{M}_{2n}\subset L$
with $\check{M}_{2n}=U^{\oplus2}\oplus E_{8}(-1)^{\oplus2}\oplus\langle-2n\rangle$.
We arrange this into two different orthogonal decompositions,
\[
M_{2n}\oplus(M_{2n})^{\perp}=M_{2n}\oplus U\oplus\check{M}_{2n}=(\check{M}_{2n})^{\perp}\oplus\check{M}_{2n},
\]
which indicates that a family of $M_{2n}$-polarized K3 surfaces is
mirror symmetric to a family of $\check{M}_{2n}$-polarized K3 surfaces
in the sense of \cite{Dol}. Since the primitive embedding $(\check{M}_{2n})^{\perp}=M_{2n}\oplus U\subset L$
is unique up to isometry (see \cite[Thm.1.14.4]{Ni},\cite[Thm.2.8]{Mo1}),
the quotient $\Omega_{\check{M}_{2n}}/O((\check{M}_{2n})^{\perp})_{+}$
naturally appears as the moduli space of $\check{M}_{2n}$-polarizable
K3 surfaces when describing mirror symmetry. 

%\vskip0.3cm
\begin{rem}
In \cite{HLOY1}, it was found that the number of Fourier-Mukai partners
of an $M_{2n}$-polarized K3 surface can be read from the ``monodromy''
groups $O((\check{M}_{2n})^{\perp})_{+}$ and $O((\check{M}_{2n})^{\perp})_{+}^{*}$
giving the moduli spaces of $\check{M}_{2n}$-polarizable and $\check{M}_{2n}$-polarized
K3 surfaces, respectively. For $n>1$, we can describe these groups
by 
\begin{equation}
O((\check{M}_{2n})^{\perp})_{+}\simeq\Gamma_{0}(n)_{+},\;\;O((\check{M}_{2n})^{\perp})_{+}^{*}\simeq\Gamma_{0}(n)_{+n}\label{eq:Gamma0N-iso}
\end{equation}
where $\Gamma_{0}(n):=\left\{ \left(\begin{smallmatrix}a & b\\
c & d
\end{smallmatrix}\right)\in SL(2,\mbZ)\mid\left(\begin{smallmatrix}a & b\\
c & d
\end{smallmatrix}\right)\equiv\left(\begin{smallmatrix}1 & *\\
0 & 1
\end{smallmatrix}\right)\text{ mod }n\right\} $, and $\Gamma_{0}(n)_{+}$, $\Gamma_{0}(n)_{+n}$ are normalizers
of $\Gamma_{0}(n)$ in $SL(2,\mbR)$ obtained by adjoining the Fricke
and Atkins-Lehner involutions; precisely, we define $\Gamma_{0}(n)_{+n}=\langle\Gamma_{0}(n),W_{n}\rangle$
with the Fricke involution $W_{n}=\left(\begin{smallmatrix}0 & \frac{1}{\sqrt{n}}\\
-\frac{1}{\sqrt{n}} & 0
\end{smallmatrix}\right)$, and similarly for $\Gamma_{0}(n)_{+}$ with adjoining all Atkins-Lehner
involutions $W_{e}$ for $e||n$ (i.e. $e|n$ and $(e,n/e)=1$) (see
\cite{CN}). In \cite[Thm.1.18]{HLOY1}, it was found that the index
$[\Gamma_{0}(n)_{+}:\Gamma_{0}(n)_{+n}]$ coincides with the number
of Fourier-Mukai partners of a general $M_{2n}$-polarized K3 surface
\cite{Og}, and interpreted from the viewpoint of symplectic diffeomorphisms
of mirror $\check{M}_{2n}$-polarizable K3 surfaces. 
\end{rem}

In the examples in Section \ref{sec: tauBCOV-eta-product}, we will
see that the monodromy groups (\ref{eq:Gamma0N-iso}) arise from the
Picard-Fuchs differential equations for suitable families of K3 surfaces
studied in \cite{LY,Do,VY}. 

\subsection{\label{subsec:AutForm}Automorphic forms}

For each $[\omega]\in\Omega_{\check{M}}$, the linear span $\mbR\,\mathrm{Re}\omega+\mbR\,\mathrm{Im}\omega$
defines a positive plane in $M^{\perp}\otimes\mathbb{R}$. As the
space parametrizing these positive two planes, the period domain $\Omega_{\check{M}}=\Omega(\check{M}^{\perp})$
is isomorphic to $O(2,20-\check{r})/SO(2)\times O(20-\check{r})$,
which is a symmetric domain of type IV. The BCOV cusp forms, which
we shall introduce in the next section, are special automorphic forms
on $\Omega_{\check{M}}$ with respect to the natural action by $O(\check{M}^{\perp})_{+}$. 

\vskip0.3cm

Let us assume a decomposition $\check{M}^{\perp}=U\oplus M$ with
$r+\check{r}=20$ as in the preceding subsection. Then we can describe
the period domain $\Omega_{\check{M}}$ by a tube domain; 
\begin{equation}
\Omega_{\check{M}}\simeq M\otimes\mbR+iC_{M},\label{eq:tubeD}
\end{equation}
where $C_{M}:=\left\{ x\in M\otimes\mbR\mid(x,x)>0\right\} ^{+}$
(as before $\left\{ \;\right\} ^{+}$represents one of the two connected
components). To describe the isomorphism (\ref{eq:tubeD}) explicitly,
we take a basis $e,f;b_{1},\cdots,b_{r}$ of $\check{M}^{\perp}=U\oplus M$
satisfying $(e,f)=1,(e,e)=(f,f)=0$ and $(e,b_{i})=(f,b_{i})=0$ for
$U=\mbZ e+\mbZ f$ and $M=\oplus_{i}\mbZ b_{i}$. We write a complex
vector $w\in(U\oplus M)\otimes\mbC$ by $w=(w_{0},w^{(2)},w_{1}^{(1)},\cdots,w_{r}^{(1)})=(w_{0},w^{(2)})\oplus w^{(1)}$
according to the decomposition $\check{M}^{\perp}=U\oplus M$. Then
we have 
\begin{equation}
(w,w)=2w_{0}w^{(2)}+(w^{(1)},w^{(1)})_{M},\label{eq:def-w.w}
\end{equation}
where $(\,,\,)_{M}$ represents the restriction of the bilinear form
$(\,,\,)$ to $M\otimes\mbC$. When we introduce inhomogeneous coordinates
by $t_{k}=\frac{w_{k}^{(1)}}{w_{0}}$ and $F=-\frac{w^{(2)}}{w_{0}}$,
we can express the isomorphism (\ref{eq:tubeD}) by 
\[
[w]=[1,-F,t_{1},\cdots,t_{r}]\mapsto(t_{1},\cdots,t_{r}).
\]
Here we note that quadratic relation $(w,w)=0$ for $\Omega_{\check{M}}$
is solved by $F=\frac{1}{2}(t)_{M}^{2}$. The coordinates $t_{1},\cdots,t_{r}$
(or $q^{2\pi it_{1}},\cdots,q^{2\pi it_{r}}$) will be often referred
to as the tube domain coordinates. 

\vskip0.3cm

Recall the natural action of the arithmetic group $O(\check{M}^{\perp})_{+}$
on $\Omega_{\check{M}}$, which we write by the linear action of $g\in O(\check{M}^{\perp})_{+}$
on $w\in\check{M}^{\perp}\otimes\mbC$ by $g\cdot w$. If we write
the corresponding linear action on the ``normalized'' vector $(1,-F,t_{k}^{(1)})$
by $g\cdot(1,-F,t_{k}^{(1)})$$=(D(g,t),A(g,t),B(g,t))$, then this
defines an action on the coordinates $t_{1},\cdots,t_{r}$ by 
\begin{equation}
[g\cdot(1,-F,t_{k}^{(1)})]=[1,\frac{A(g,t)}{D(g,t)},\frac{B_{k}(g,t)}{D(g,t)}]=:[1,-\tilde{F}(t),\tilde{t}_{k}^{(1)}],\label{eq:g-action-D(g,t)}
\end{equation}
where $\tilde{t}_{k}^{(1)}:=\frac{B_{k}(g,t)}{D(g,t)}\,(k=1,...,r$)
are rational functions of $t_{1},\dots,t_{r}$. We write this action
by $\tilde{t}=g\cdot t$. 
\begin{defn}
A holomorphic function $\Psi(t)=\Psi(t_{1},\dots,t_{r})$ on $\Omega_{\check{M}}$
having the property 
\begin{equation}
\Psi(g\cdot t)=D(g,t)^{m}\Psi(t)\,\,\,(g\in O(\check{M}^{\perp})_{+})\label{eq:Def-autoForm}
\end{equation}
is called an automorphic form of weight $m$ on $\Omega_{\check{M}}$. 
\end{defn}

\subsubsection{Example 1 $\mathrm{(}\check{M}^{\perp}=U\oplus U\mathrm{)}$ \label{subsec:Example-1-UU}}

Consider a primitive embedding $\check{M}=U\oplus E_{8}(-1)^{\oplus2}\hookrightarrow L$
and the orthogonal lattice $\check{M}^{\perp}$ in $L$. It is easy
to see $\check{M}^{\perp}=U\oplus U(=U\oplus M)$, and the isomorphism
$\Omega_{\check{M}}\simeq M\otimes\mbR+iC_{M}=\mbH_{+}\times\mbH_{+}$
in terms of the upper half plane $\mbH_{+}$. Since the primitive
embedding $\check{M}^{\perp}\subset L$ is unique up to isomorphism
by \cite[Thm.1.14.4]{Ni},\cite[Thm.2.8]{Mo1}, the quotient of $\Omega_{\check{M}}$
by the group $O(\check{M}^{\perp})_{+}$ parametrizes the isomorphism
classes of $\check{M}$-polarizable K3 surfaces. For the group $O(\check{M}^{\perp})_{+}=O(U\oplus U)_{+}$,
we note the following isomorphism: 
\[
O(U\oplus U)_{+}\simeq P(SL(2,\mbZ)\times SL(2,\mbZ))\rtimes\mbZ_{2}.
\]
Explicit generators of $O(U\oplus U)$ can be found in \cite[Prop.2.5]{HLOY2};
using the notation there, the first factor $P(SL(2,\mbZ)\times SL(2,\mbZ))$
is generated by $S_{1},S_{2},T_{1},T_{2}$ and the second factor $\mbZ_{2}$
is generated by $R_{1}R_{2}$. Precisely, since these generators are
given for $O(2,2,\mbZ)$, we need to translate these to our definition
of $O(U\oplus U)$. After the necessary translation, we verify the
following actions on the tube domain coordinates $(t_{1},t_{2})\in\mbH_{+}\times\mbH_{+}$:
\[
\begin{alignedat}{2} & S_{1}:(t_{1},t_{2})\mapsto(-\frac{1}{t_{1}},t_{2})\,,\; &  & S_{2}:(t_{1},t_{2})\mapsto(t_{1},-\frac{1}{t_{2}})\\
 & T_{1}:(t_{1},t_{2})\mapsto(t_{1}+1,t_{2})\,,\; &  & T_{2}:(t_{1},t_{2})\mapsto(t_{1},\frac{t_{2}}{t_{2}+1})
\end{alignedat}
\]
and $R_{1}R_{2}:(t_{1},t_{2})\mapsto(-\frac{1}{t_{2}},-\frac{1}{t_{1}})$. 

\subsubsection{Example 2 $\mathrm{(}\check{M}^{\perp}=U^{\oplus2}\oplus\langle-2\rangle\mathrm{)}$\label{subsec:Example-2-UUA1}}

Primitive embedding $\check{M}=U\oplus E_{8}(-1)\oplus E_{7}(-1)\subset L$
determines an orthogonal lattice $\check{M}^{\perp}\simeq U^{\oplus2}\oplus\langle-2\rangle$
in $L$. Since the primitive embedding $\check{M}^{\perp}\subset L$
is unique up to isomorphism as in Example 1, the quotient $\Omega_{\check{M}}/O(\check{M}^{\perp})_{+}$
parametrizes the isomorphism classes of $\check{M}$-polarizable K3
surfaces. It is known that the period domain $\Omega_{\check{M}}$
is isomorphic to the Siegel upper half space $\mbH_{2}$ of genus
two, and the action of $O(\check{M}^{\perp})_{+}=O(U^{\oplus2}\oplus\langle-2\rangle)$
on $\Omega_{\check{M}}$ coincides with the standard modular action
of $Sp(4,\mbZ)$ on $\mbH_{2}$. More precisely, there is a group
isomorphism 
\[
\Psi:Sp(4,\mbZ)/\left\{ \pm I_{4}\right\} \stackrel{\sim}{\rightarrow}SO(U^{\oplus2}\oplus\langle-2\rangle)_{+},
\]
where $SO(U^{\oplus2}\oplus\langle-2\rangle)_{+}=SO(U^{\oplus2}\oplus\langle-2\rangle)\cap O(U^{\oplus2}\oplus\langle-2\rangle)_{+}$.
There is also an isomorphism $\phi:\mbH_{2}\simeq\Omega_{\check{M}}(\simeq M\otimes\mbR+iC_{M}$
with $M=U\oplus\langle-2\rangle$) given by 
\[
\tau=\left(\begin{matrix}\tau_{11} & \tau_{12}\\
\tau_{21} & \tau_{22}
\end{matrix}\right)\mapsto\phi(\tau)=[1,-\det\tau,-\tau_{11},-\tau_{22},\tau_{12}],
\]
where we note that the condition for the elements $(-\tau_{11},-\tau_{22},\tau_{12})\in M\otimes\mbR+iC_{M}$
is represented by $\mathrm{Im}(\tau_{11})\mathrm{Im}(\tau_{22})-\mathrm{Im}(\tau_{12})^{2}>0$
and $\mathrm{Im}(\tau_{11})>0$, i.e., $\tau\in\mathbb{H}_{2}$. The
following diagram shows that the isomorphism $\phi$ is compatible
with the group actions: 
\begin{equation}
\begin{matrix}\mbH_{2} & \stackrel{g}{\longrightarrow} & \mbH_{2}\\
\downarrow\phi &  & \downarrow\phi\\
\Omega_{\check{M}} & \stackrel{\Psi(g)}{\longrightarrow} & \Omega_{\check{M}}
\end{matrix}\label{eq:H2-comm-diag}
\end{equation}
where $\Psi(g)\in SO(U^{\oplus2}\oplus\langle-2\rangle)_{+}$ and
$g\in Sp(4,\mathbb{Z})/\left\{ \pm I_{4}\right\} $. If we write the
standard action of $g=\left(\begin{matrix}A & B\\
C & D
\end{matrix}\right)\in Sp(4,\mbZ)$ on $\tau\in\mbH_{2}$ by 
\[
g\s\tau=(A\tau+B)(C\tau+D)^{-1},
\]
then we can verify the above commutative diagram by the following
relation: 
\[
\Psi(g)\phi(\tau)=[1,-\det(g\s\tau),-(g\s\tau)_{11},-(g\s\tau)_{22},(g\s\tau)_{12}]=\phi(g\s\tau),
\]
with the explicit form of $\Psi(g)$ given in Appendix \ref{subsec:AppendixA-2}.
Noting that $SO(U^{\oplus2}\oplus\langle-2\rangle)_{+}\simeq O(U^{\oplus2}\oplus\langle-2\rangle)_{+}/\left\{ \pm I_{5}\right\} $
and $\left\{ \pm I_{5}\right\} $ acts trivially on $\Omega_{\check{M}}$,
we have $\mathbb{H}_{2}/Sp(4,\mathbb{Z})\simeq\Omega_{\check{M}}/O(U^{\oplus2}\oplus\langle-2\rangle)_{+}$.
Hence the ring of the automorphic forms coming from $\check{M}$-polarizable
K3 surfaces parametrized by $\Omega_{\check{M}}$ coincides with that
of the standard modular forms on $\mbH_{2}$. 

\vskip2cm

\section{\textbf{\label{sec:BCOV-cusp-forms}BCOV cusp forms}}

We will introduce the BCOV formula to define the BCOV cusp forms.
The BCOV cusp forms are special cusp forms on the period domains of
lattice polarized (or polarizable) K3 surfaces. When we have explicit
families of lattice polarized K3 surfaces, we can compute the BCOV
formula in terms of period integrals of the families. 

\subsection{\label{subsec:MS-Mmap}Mirror symmetry and mirror map }

Suppose we have a family $\pi:\check{\mathfrak{X}}\rightarrow\mathcal{M}$
of $\check{M}$-polarized K3 surfaces over some parameter space. Here
and hereafter, we allow the fibers of the family develop ordinary
double points over a subset $D\subset\mathcal{M}$, which we call
discriminant of the family. Associated to the local system $R^{2}\pi_{*}\mbC_{\check{\mathfrak{X}}}$,
we have the period map 
\[
\mathcal{P}:\mathcal{M}\rightarrow\Omega_{\check{M}},\,\,\,x\mapsto[\int_{\gamma_{1}}\omega_{x},\cdots,\int_{\gamma_{22-\check{r}}}\omega_{x}],
\]
where $\omega_{x}:=\omega(\check{X}_{x})$ is a holomorphic two form
of the fiber $\check{X}_{x}$ over $x\in\mathcal{M}$, and $\left\{ \gamma_{1},\cdots,\gamma_{22-\check{r}}\right\} $
is a basis of the transcendental cycles in $H_{2}(X_{x_{0}},\mathbb{Z})$
with a base point $x_{0}$. Assuming the decomposition $M\oplus(M)^{\perp}=M\oplus U\oplus\check{M}$,
we have the isomorphism $\Omega_{\check{M}}\simeq M\otimes\mbR+iC_{M}$
as in (\ref{eq:tubeD}) and also the relation $r+\check{r}=20$. As
described in Subsect. \ref{subsec:AutForm}, the coordinates $t_{k}\in M\otimes\mbR+iC_{M}$
are given by ratios of period integrals, and actually give a map from
the universal covering of the parameter space $\mathcal{M}$ to the
tube domain (which is called the complexified K\"ahler moduli space
of $M$-polarized K3 surfaces). This multi-valued map $\mathcal{M}\rightarrow M\otimes\mbR+iC_{M}$
or its inverse relation is called the \textit{mirror map} (see \cite{LY,Do,VY}). 

\subsubsection{Picard-Fuchs equation and local solutions}

Period integrals can be characterized as the solutions of a certain
set of differential equations, called the Picard-Fuchs equations.
By making solutions of the Picard-Fuchs equations near special boundary
points, called large complex structure limits (LCSLs), we can determine
the mirror map. In the case of Calabi-Yau complete intersections in
toric varieties, it is known that LCSLs exist for a suitable toric
compactification $\overline{\mathcal{M}}$ of $\mathcal{M}$ (see
\cite{HLY,HLY2} and references therein). For a family of K3 surfaces,
near a LCSL, we can find period integrals $w_{0}(x),w_{2}(x);w_{k}^{(1)}(x)$
which behave as
\begin{equation}
\begin{aligned}w_{0}(x)=1+O(x),\,\,\,w_{2}(x)=-\sum K_{kl}\big(\frac{1}{2\pi i}\big)^{2}(\log x_{k})(\log x_{l})w_{0}(x)+\cdots\\
w_{k}^{(1)}(x)=\frac{1}{2\pi i}(\log x_{k})w_{0}(x)+O(x),\qquad\qquad\qquad\qquad
\end{aligned}
\label{eq:localSol}
\end{equation}
and satisfy a quadratic relation; $2w_{0}(x)w_{2}(x)+\sum K_{ij}w_{i}^{(1)}(x)w_{j}^{(1)}(x)=0$
where $(K_{ij})$ represents the Gram matrix of the lattice $M$ with
respect to the natural basis given in \cite[Eq.(2.11)]{HosIIA}. The
mirror map is determined by inverting the relations $t_{k}=\frac{1}{2\pi i}\frac{w_{k}^{(1)}(x)}{w_{0}(x)}=\frac{1}{2\pi i}\log x_{k}+O(x)$,
which give rise to the following $q$-series expansion: 
\begin{equation}
x_{k}=x_{k}(q_{1},...,q_{r})=q_{k}+O(q)\,,\,\,\,\,\,(q_{\ell}=e^{2\pi i\,t_{\ell}}).\label{eq:mirrorMap}
\end{equation}
Due to the quadratic relation, we can arrange the period integrals
as $[w]=[w_{0}(x),$ $w_{2}(x),w_{1}^{(1)},...,w_{r}^{(1)}]$ so that
we have the period map $\mathcal{P}:\mathcal{M}\rightarrow\Omega_{\check{M}}$
in the from given in Subsect. \ref{subsec:AutForm}. In this form,
the boundary point $q_{1}=q_{2}=\cdots=q_{r}=0$ (i.e., $x_{1}=\cdots=x_{r}=0)$
is expected to correspond to a cusp point which we attach to the period
domain $\Omega_{\check{M}}$ in the Baily-Borel-Satake compactification.
As we will see in the next section, this is true when $r=1$. However,
as we can see in the example in a recent work \cite{HLTY}, the correspondence
is more involved for $r>1$. 

\subsubsection{Monodromy group $\mathcal{G}_{m}$}

We denote by \textbf{$\mathcal{G}_{m}$} the global monodromy group
of the solutions (\ref{eq:localSol}) for the family $\check{\mathfrak{X}}\rightarrow\mathcal{M}$.
By definition of the family, it is a subgroup of $O(\check{M}^{\perp})_{+}$.
As we will see in Sect.~\ref{sec: tauBCOV-eta-product}, there are
many examples for which $\mathcal{G}_{m}=O(\check{M}^{\perp})_{+}$
or $O(\check{M}^{\perp})_{+}^{*}$ holds when $r=1$. However, the
monodromy group $\mathcal{G}_{m}$ is only a subgroup of $O(\check{M}^{\perp})_{+}$
in general for $r>1$. It acts on $\Omega_{\check{M}}$ and reflects
the multi-valued property of the period map $\mathcal{P}:\mathcal{M}\rightarrow\Omega_{\check{M}}$.
Then the mirror map $x_{k}=x_{k}(q_{1},...,q_{r})$ as the inverse
of the period map is invariant under the action $t\mapsto g\cdot t$
of $g\in\mathcal{G}_{m}\subset O(\check{M}^{\perp})_{+}$ on the tube
domain coordinates. 

\subsubsection{Mirror symmetry and Griffiths-Yukawa coupling}

In the study of mirror symmetry of K3 surfaces, there are no quantum
corrections from the Gromov-Witten invariants. Hence, the so-called
Griffiths-Yukawa couplings $C_{ij}=\int\omega_{x}\frac{\partial^{2}\;\;}{\partial x_{i}\partial x_{j}}\omega_{x}$,
with a suitable normalization, reproduce the Gram matrix $(K_{ab})$
of $M$ by the formula,
\begin{equation}
\frac{1}{w_{0}(x(q))^{2}}\sum_{i,j}C_{ij}\frac{\partial x_{i}}{\partial t_{a}}\frac{\partial x_{j}}{\partial t_{b}}=K_{ab}.\label{eq:YuakawaYttt}
\end{equation}

\subsection{BCOV formula and cusp forms}

Here we define the BCOV formulas for the families of $\check{M}$-polarized
K3 surface. As in the preceding subsection, we consider a primitive
embedding $\check{M}\subset L$ with the property $\check{M}^{\perp}=U\oplus M$. 

Let us assume that we have a family $\check{\mathfrak{X}}\rightarrow\mathcal{M}$
of $\check{M}$-polarized K3 surfaces over $\mathcal{M}$ as in Subsect.
\ref{subsec:MS-Mmap}. In order to have the topological limit (\ref{eq:BCOV-CY3}),
assume also that we have a LCSL at $x_{1}=x_{2}=\cdots=x_{r}=0$ by
finding a compactification $\overline{\mathcal{M}}$ of $\mathcal{M}$.
The affine chart by the coordinates $x_{1},\cdots,x_{r}$ are supposed
to be dense in $\overline{\mathcal{M}}$; e.g., this is the case if
$\overline{\mathcal{M}}$ is given by a toric variety as in the examples
in the following sections. A prototypical example is given by the
toric variety of the secondary polytope. Near the LCSL, we have period
integrals in the form of (\ref{eq:localSol}). Following the definition
(\ref{eq:mirrorMap}), we determine the mirror map $x_{k}=x_{k}(q)$
which has $q$-series expansions of the form $x_{k}=q_{k}+O(q)$.
Over the parameter space $\mathcal{M},$ the fibers of the family
degenerate on the discriminant $\cup_{i}D_{i}$, where $D_{i}=\left\{ dis_{i}(x)=0\right\} $
represent irreducible components. 
\begin{defn}[\textbf{BCOV cusp form}]
\label{def:tau-BCOV} Under the assumptions above, near the LCSL,
we define the BCOV formula $\tau_{_{\mathrm{BCOV}}}$, up to a multiplication
constant, by 
\begin{equation}
\tau_{_{\mathrm{BCOV}}}(t)=\left(\frac{1}{w_{0}(x)}\right)^{r+1}\frac{\partial(x_{1},\cdots,x_{r})}{\partial(t_{1},\cdots,t_{r})}\prod_{i}dis_{i}(x)^{r_{i}}\,\prod_{j=1}^{r}x_{j}^{-1+a_{j}},\label{eq:tauBCOV-x}
\end{equation}
where $r_{i}$ and $a_{j}$ are rational numbers which are restricted
by suitable boundary conditions (which we will describe in Subsect.
\ref{subsec:Bound-cond}). When $(\tau_{_{\mathrm{BCOV}}}(t))^{-1}$
defines a cusp form on $\Omega_{\check{M}}$, i.e., an automorphic
form which is holomorphic and vanishes at all cusps, we call it a
BCOV cusp form. 
\end{defn}

In Sects. \ref{sec: tauBCOV-eta-product} and \ref{sec:CD-family },
we will find in many examples that $(\tau_{_{\mathrm{BCOV}}}(t))^{-1}$
defines a cusp form on the period domain $\Omega_{\check{M}}$ for
suitable values of $r_{i}$ and $a_{i}$. Here we first show its modular
property. 
\begin{prop}
\label{prop:Modular-tau} The inverse power $(\tau_{_{\mathrm{BCOV}}}(t)){}^{-1}$
of the BCOV formula has weight one with respect to the action of $\mathcal{G}_{m}\subset O(\check{M}^{\perp})_{+}$.
Namely, it satisfies 
\begin{equation}
\big(\tau_{_{\mathrm{BCOV}}}(g\cdot t)\big)^{-1}=v(g)D(g,t)\big(\tau_{_{\mathrm{BCOV}}}(t)\big)^{-1}\label{eq:tau(gt)}
\end{equation}
with the weight factor $D(g,t)$ given in (\ref{eq:g-action-D(g,t)}),
and $v(g)$ is a multiplier system which is independent of $t$ and
satisfies $|v(g)|=1$. 
\end{prop}

We prove this property by showing the following. 
\begin{lem}
Let $[w]=[1,-F,t_{1},\cdots,t_{r}]$ be a point in $\Omega_{\check{M}}\simeq M\otimes\mathbb{R}+iC_{M}$
in terms of the tube domain coordinates, and write the action of $g\in\mathcal{G}_{m}$
on it by $[g\cdot w]=[1,-\tilde{F},\tilde{t}_{1},\cdots,\tilde{t}_{r}]$
as in (\ref{eq:g-action-D(g,t)}). Then it holds that 
\[
d\tilde{t}_{1}\wedge d\tilde{t}_{2}\wedge\cdots\wedge d\tilde{t}_{r}=(\det g)\,D(g,t)^{-r}dt_{1}\wedge dt_{2}\wedge\cdots\wedge dt_{r}.
\]
\end{lem}

\begin{proof}
Recall that the tube domain coordinates $t_{1},\cdots,t_{r}$ are
affine coordinates of the quadric $\left\{ (w,w)=0\right\} $ in $\mathbb{P}(\check{M}^{\perp}\otimes\mathbb{C})$.
Let us write a point $[Z]\in\mathbb{P}(\check{M}^{\perp}\otimes\mathbb{C})$
by $[Z_{0},Z_{1},\cdots,Z_{r+1}]=[u,v,z_{1},\cdots,z_{r}]$ according
to the decomposition $\check{M}^{\perp}\otimes\mathbb{C}=(U\oplus M)\otimes\mathbb{C}$.
Then the quadric relation is expressed by $(Z,Z)=2uv+(z,z)_{M}=0$.
The orthogonal group $O(\check{M}^{\perp})$ acts naturally on the
linear coordinates $Z=(Z_{0},Z_{1},\cdots,Z_{r+1})$ preserving the
quadric. It also preserves the standard $r+1$ form $d\mu_{Z}:=\sum(-1)^{i}Z_{i}\wedge dZ_{0}\wedge\cdots\widehat{dZ_{i}}\cdots\wedge dZ_{r+1}$
up to sign, since we can write $d\mu_{Z}=i_{E}dZ_{0}\wedge dZ_{1}\wedge\cdots\wedge dZ_{r+1}$
by the interior product $i_{E}$ with the Euler vector field $E=Z_{0}\frac{\partial\;}{\partial Z_{0}}+\cdots+Z_{r+1}\frac{\partial\;}{\partial Z_{r+1}}$.
Note that in the affine chart $Z_{0}\not=0$, we can solve the quadric
for $v$ and obtain $[w]=[1,-\frac{1}{2}(t,t)_{M},t_{1},\cdots,t_{r}]$
with $t_{i}=\frac{z_{i}}{u}$ for $i=1,\cdots,r$. Hence we can express
the holomorphic $r$-form on the quadric by the following residue
around the quadric:
\begin{equation}
Res_{(Z,Z)=0}\left(\frac{d\mu_{Z}}{2uv+(z,z)_{M}}\right)=\frac{u^{r}}{2}dt_{1}\wedge dt_{2}\wedge\cdots\wedge dt_{r}.\label{eq:Res-dt}
\end{equation}
Now consider the linear action $\tilde{Z}:=g\cdot(Z_{0},Z_{1},\cdots,Z_{r+1})$
of $g\in\mathcal{G}_{m}\subset O(\check{M}^{\perp})$ on $Z$. If
we solve the quadric $(\tilde{Z},\tilde{Z})=0$ for $\tilde{v}$ in
$\tilde{Z}=(\tilde{u},\tilde{v},\tilde{z_{1}},\cdots,\tilde{z}_{r})$,
then we have $[\tilde{w}]=[1,-\frac{1}{2}(\tilde{t},\tilde{t})_{M},\tilde{t}_{1},\cdots,\tilde{t}_{r}]$
with $\tilde{t}_{i}=\frac{\tilde{z_{i}}}{\tilde{u}}$. Since the l.h.s.
of (\ref{eq:Res-dt}) is invariant up to the factor $\det g$ under
the linear action of $g\in\mathcal{G}_{m}$ on $Z$, we have
\[
\frac{\tilde{u}^{r}}{2}d\tilde{t}_{1}\wedge d\tilde{t}_{2}\wedge\cdots\wedge d\tilde{t}_{r}=(\det g)\frac{u^{r}}{2}dt_{1}\wedge dt_{2}\wedge\cdots\wedge dt_{r}.
\]
Identifying $\frac{\tilde{u}}{u}$ with the factor $D(g,t)$ in (\ref{eq:g-action-D(g,t)}),
we obtain the claimed relation. 
\end{proof}
\noindent \textit{Proof of Proposition \ref{prop:Modular-tau}.}
We write the period map $\mathcal{P}:\mathcal{M}\rightarrow\Omega_{\check{M}}$
described in Subsect. \ref{subsec:MS-Mmap} by 
\[
\mathcal{P}(x)=[w_{0}(x),w_{2}(x),w_{1}^{(1)}(x),\cdots,w_{r}^{(1)}(x)]
\]
for $x\in\mathcal{M}.$ This period map is multi-valued due to the
monodromy of period integrals. Using the mirror map $x=x(t)$ in terms
of the tube domain coordinates $t=(t_{1},\cdots,t_{r})$, we can express
the monodromy property of $\mathcal{P}(x)$ by 
\[
\begin{aligned}\mathcal{P}(x(g\cdot t)) & =[w_{0}(x(g\cdot t)),w_{2}(x(g\cdot t)),w_{k}^{(1)}(x(g\cdot t))]\\
 & =[g\cdot(\,w_{0}(x(t)),w_{2}(x(t)),w_{k}^{(1)}(x(t))\,)]\\
 & =[D(g,t)w_{0}(x),A(g,t)w_{0}(x),B_{k}(g,t)w_{0}(x)]
\end{aligned}
\]
for $g\in\mathcal{G}_{m}$, where $x(g\cdot t)$ should be read as
a point on a covering $\widetilde{\mathcal{M}}$ of $\mathcal{M}$
which projects to $x(t)\in\mathcal{M}$. Note that the last two equalities
hold without the projectivization. Hence, we have 
\begin{equation}
w_{0}(x(g\cdot t))=D(g,t)\,w_{0}(x(t))\label{eq:hodge-L}
\end{equation}
with the automorphic factor $D(g,t)$. Also, by the above lemma, we
have 
\[
\frac{\partial(\tilde{t}_{1},\tilde{t}_{2},\cdots,\tilde{t}_{r})}{\partial(t_{1},t_{2},\cdots,t_{r})}=\det g\,(D(g,t))^{-r}.
\]
Using these relations, and the fact that the mirror map $x_{k}=x_{k}(q)$
is invariant under the monodromy group $\mathcal{G}_{m}$, we obtain
the claimed relation (\ref{eq:tau(gt)}). The multiplier system $v(g)$
depends on the branches of the factor $\prod x_{j}^{-1+a_{j}}$ if
some of $a_{j}$'s are not integral. $\hfill\,$$\square$

~

We remark that the relation (\ref{eq:hodge-L}) comes from the fact
that $w_{0}(x)$ is a section of the Hodge line bundle; which is given
by the $(2,0)$ part of the relative Hodge decomposition of the sheaf
$R^{2}\pi_{*}\mathbb{C}_{\check{\mathfrak{X}}}$ over $\mathcal{M}$
for the family $\pi:\check{\mathfrak{X}}\rightarrow\mathcal{M}$. 

\subsection{\label{subsec:Bound-cond}Boundary conditions}

The above modular property holds for any choice of the parameters
$r_{i}$ and $a_{i}$ in the BCOV formula. We can restrict these parameters
by imposing natural boundary conditions for the family $\pi:\check{\mathfrak{X}}\rightarrow\mathcal{M}$
(cf. \cite{BCOV1}), which we describe as follows: 

\vskip0.1cm

(1) The first condition comes from the so-called conifold degenerations.
As we defined, the fiber K3 surfaces of the family $\pi:\check{\mathfrak{X}}\rightarrow\mathcal{M}$
degenerate over the discriminant $D$, where general fibers acquire
ordinary double points (ODPs). The N\'eron-Severi lattice enlarges
after blowing up these ODPs. We require that $\log(\tau_{_{\mathrm{BCOV}}}(t))^{-1}$
is regular, i.e., $(\tau_{_{\mathrm{BCOV}}}(t))^{-1}$ is regular
and non-vanishing over the discriminant $D\subset\mathcal{M}$. We
call this condition a\textit{ conifold regularity}. 

\vskip0.1cm

(2) The second condition comes from the boundary points contained
in $\overline{\mathcal{M}}\setminus\mathcal{M}$. Suppose that there
is a boundary point given by the intersection point of normal crossing
divisors such that all monodromy actions $T_{i}$ around them are
of finite order. We call such a boundary point an \textit{orbifold
point}. Since the fiber K3 surfaces naturally extend to orbifold points,
we require that $\log(\tau_{_{\mathrm{BCOV}}}(t))^{-1}$ is regular
at the corresponding point in $\Omega_{\check{M}}$. We call this
boundary condition an \textit{orbifold regularity} at the boundary
point. If there are more than one orbifold points, the regularity
is imposed by choosing one. Hence, depending a choice of an orbifold
point, we may have different boundary conditions (see Subsect. \ref{subsec:CD-orbifold-reg}). 

\vskip0.1cm

These two conditions are natural translation of those formulated for
the case of Calabi-Yau threefolds \cite{BCOV1}. It will be observed
that the parameters $r_{i}$ and $a_{i}$ are uniquely determined
by these conditions for most cases in our examples in the next sections.

\subsubsection{Conifold regularity\label{subsec:Conifold-regularity}}

We can implement the condition (1) by using the mirror symmetry relation
(\ref{eq:YuakawaYttt}). Let us calculate the determinant of (\ref{eq:YuakawaYttt})
as 
\[
\det(C_{ij})\left(\frac{\partial(x_{1},\cdots,x_{r})}{\partial(t_{1},\cdots,t_{r})}\right)^{2}=\det(K_{ab})\omega_{0}(x)^{2r}.
\]
Using this relation, we have 
\begin{equation}
\tau_{_{\mathrm{BCOV}}}(t)=\det(K_{ab})^{\frac{1}{2}}\frac{\prod x_{j}^{-1+a_{j}}}{w_{0}(x)}\left(\det(C_{ij})\right)^{-\frac{1}{2}}\prod_{i}dis_{i}(x)^{r_{i}}.\label{eq:tauBCOV-Cij}
\end{equation}
In this form, it is easy to implement the conifold regularity since
it is straight forward to evaluate $\det(C_{ij})$ for a given family
$\pi:\check{\mathfrak{X}}\rightarrow\mathcal{M}$. In particular,
for all the examples in the present paper, we will see that the following
relation holds:
\begin{equation}
\det(C_{ij})\prod_{i}dis_{i}(x)\,\prod_{j=1}^{r}x_{j}^{2}=\text{const.}\,\,.\label{eq:detCij}
\end{equation}
Using this, we see that the conifold regularity holds if $r_{i}=-\frac{1}{2}$
for all $i$. Namely, if we verify (\ref{eq:detCij}), the parameters
$r_{i}$ are determined completely by the conifold regularity. 

\subsubsection{Orbifold regularity\label{subsec:Orbifold-regularity}}

We implement the orbifold regularity by choosing an orbifold point
from $\overline{\mathcal{M}}\setminus\mathcal{M}$. Let us assume
that the relation (\ref{eq:detCij}) holds. Then by the conifold regularity,
we have the BCOV formula in the following form: 
\begin{equation}
\tau_{_{\mathrm{BCOV}}}(t)=\text{const.}\frac{\prod x_{j}^{a_{j}}}{w_{0}(x)}.\label{eq:tauBCOV2}
\end{equation}
Let $\tilde{x}_{i}$ be the affine coordinates centered at the orbifold
point. For simplicity, we assume that $\tilde{x}_{i}$ are related
to the affine coordinates $x_{i}$ centered at the LCSL by Laurent
monomials $x_{i}=\varphi_{i}(\tilde{x})$ (see the examples in the
following sections). Recall that the period integrals around the LCSL
are arranged so that $\mathcal{P}(x)=[w_{0}(x),w_{2}(x),w_{1}^{(1)}(x),\cdots,w_{r}^{(1)}(x)]$
is a point in $\Omega_{\check{M}}$; hence $w_{0}(x)\not=0$ holds
for $\mathcal{P}(x)$ corresponding to smooth K3 surfaces. Consider
an analytic continuation of $w_{0}(x)$ to the orbifold point, and
express it by 
\[
w_{0}(x)=c_{0}\tilde{w}_{0}(\tilde{x})+c_{1}\tilde{w}_{1}(\tilde{x})+\cdots+c_{r+2}\tilde{w}_{r+2}(\tilde{x}),
\]
where $\tilde{w}_{k}(\tilde{x})$ are (ordered) local solutions of
the Picard-Fuchs equation around the orbifold point. If we order the
solutions so that $\tilde{w}_{k}(\tilde{x})/\tilde{w}_{0}(\tilde{x})$
are regular, we have $w_{0}(x)=c_{0}\tilde{w}_{0}(\tilde{x})(1+O(\tilde{x}))$
with $c_{0}\not=0$ since $\mathcal{P}(x)$ is a point in $\Omega_{\check{M}}$
at $\tilde{x}=0$. We can now write the analytic continuation of the
BCOV formula (\ref{eq:tauBCOV2}) as 
\begin{equation}
\tau_{_{\mathrm{BCOV}}}(t)=\text{const.}\frac{\prod\tilde{x}_{j}^{\tilde{a}_{j}}}{c_{0}\tilde{w}_{0}(\tilde{x})(1+\cdots)},\label{eq:tauBCOV-orbifold}
\end{equation}
where we have set $\prod x_{j}^{a_{j}}=\prod\varphi_{j}(\tilde{x})^{a_{j}}=:\prod\tilde{x}_{j}^{\tilde{a}_{j}}$.
In this form, we require that $\log(\tau_{BCOV}(t))^{-1}$ is regular
around $\tilde{x}=0$. Note that the parameters $\tilde{a}_{j}$ are
uniquely determined by the exponent $\rho$ in $\tilde{w}_{0}(\tilde{x})=\tilde{x}^{\rho}(1+O(\tilde{x}))$. 

As above, the orbifold regularity is implemented for a choice of an
orbifold point from $\overline{\mathcal{M}}\setminus\mathcal{M}$.
Depending on the choice, the resulting $\tau_{_{\mathrm{BCOV}}}(t)$
can be different. Also, if there is no orbifold point, the regularity
condition is vacant. 

\section{\textbf{\textcolor{black}{\label{sec: tauBCOV-eta-product}$\tau_{_{\mathrm{BCOV}}}$
and eta products}}}

Here we present some explicit calculations of $\tau_{_{\mathrm{BCOV}}}$
coming from $\check{M}_{2n}$-polarizable K3 surfaces whose isomorphism
classes are parametrized by $\mbH_{+}/\Gamma_{0}(n)_{+}$. An example
of a higher dimensional period domain will be considered in the next
section. 

\subsection{Genus zero groups }

Let us consider a family of $\check{M}_{2n}$-polarizable K3 surfaces
as a mirror family of an $M_{2n}$-polarized K3 surface. As summarized
in Subsect. \ref{subsec:K3-period-domain}, the quotient $\Omega_{\check{M}_{2n}}/O(U\oplus\langle2n\rangle)_{+}\simeq\mbH_{+}/\Gamma_{0}(n)_{+}$
gives a coarse moduli space of the family in general. We will calculate
the BCOV formula when we have an explicit realization $\check{\mathfrak{X}}\rightarrow\mathcal{M}$
of the family with the parameter space being compactified by $\overline{\mathcal{M}}\simeq\mbP^{1}$.
Here, we assume the isomorphism $\overline{\mathcal{M}}\simeq\mathbb{P}^{1}$
so that we can calculate the period integrals by solving the Picard-Fuchs
differential equation over $\mathbb{P}^{1}$. If we have such a family
and the monodromy group $\mathcal{G}_{m}$ of the family coincides
with $\Gamma_{0}(n)_{+}$, then we should have $\overline{\mbH_{+}/\Gamma_{0}(n)_{+}}\simeq\mbP^{1}$
(see \cite{LY,VY,Do}). Similarly, if we have an explicit family of
$\check{M}_{2n}$-polarized K3 surfaces $\check{\mathfrak{X}}\rightarrow\mathcal{M}$
with $\overline{\mathcal{M}}\simeq\mathbb{P}^{1}$ and $\mathcal{G}_{m}\simeq\Gamma_{0}(n)_{+n}$,
we should have an isomorphism $\overline{\mbH_{+}/\Gamma_{0}(n)_{+n}}\simeq\mbP^{1}$.
The modular groups $\Gamma_{0}(n)_{+}$ and $\Gamma_{0}(n)_{+n}$
which give rise to the isomorphisms as above are called \textit{genus
zero groups}. 

The list of genus zero groups $\Gamma_{0}(n)_{+n}$ is classical due
to Fricke \cite[p.367]{Fricke}, which is given by 
\[
n=2,3,\cdots,21,23,24,25,26,27,29,31,32,35,36,39,41,47,49,50,59,71.
\]
It was extended in \cite{CN} to $\Gamma_{0}(n)_{+}$ adjoining all
involutive normalizers, which we list in Appendix \ref{subsec:AppB-table-Gamma}
with some additional data. In \cite{CN}, certain generator of the
function field of $\overline{\mbH_{+}/\Gamma_{0}(n)_{+}}\simeq\mbP^{1}$
is called the Thompson series, and its concrete formula in terms of
elliptic modular forms (in particular, the eta function) is given.
For a genus zero group, if we have an explicit realization $\check{\mathfrak{X}}\rightarrow\mathcal{M}$
as above, the Thompson series coincides exactly with the mirror map
of the family up to an additive constant (see \cite{LY,VY,Do}). 

\subsection{Family of $\check{M}_{20}$-polarizable K3 surfaces}

In \cite{LY}, many examples of the Picard-Fuchs differential equations
for families of K3 surfaces over $\mathbb{P}^{1}$ have been listed,
and the mirror maps have been identified with the Thompson series.
Many of them are $\check{M}_{2n}$-polarizable K3 surfaces with the
modular group $\Gamma_{0}(n)_{+}.$ 

As an example, let us consider the case $n=10$. To describe the family,
we start with a K3 surface $X=(1,1)^{4}\subset\mathbb{P}^{3}\times\mathbb{P}^{3}$
which is given as a complete intersection of four general $(1,1)$
divisors in $\mbP^{3}\times\mbP^{3}$, whose mirror family can be
obtained by the Batyrev-Borisov toric mirror construction \cite{BaBo}.
In fact, $X$ is a Cayley model of the Reye congruence, and its mirror
family is parametrized by $\mbP^{2}$. A family of $\check{M}_{20}$-polarizable
K3 surfaces appears when we restrict the parameter space to the diagonal
$\mbP^{1}\subset\mbP^{2}$. The Picard-Fuchs differential equation
can be read from the table of \cite{LY} as 
\begin{equation}
\left\{ \theta_{x}^{3}-2x(2\theta_{x}+1)(3\theta_{x}^{2}+3\theta_{x}+1)-x^{2}(4\theta_{x}+3)(4\theta_{x}+4)(4\theta_{x}+5)\right\} w=0,\label{eq:PFdeg20}
\end{equation}
where $\theta_{x}:=x\frac{d\;}{dx}$. This equation is singular at
$x=0,\infty$ and also at the zeros of $dis_{0}\,dis_{1}=0$ with
$dis_{0}=1+4x,$ $dis_{1}=1-16x$. We summarize the data of the singularities
in Riemann's $\mathcal{P}$ scheme as\def\RS1{\left\{ \small\begin{array}{cccc} 0 & -\frac{1}{4} & \frac{1}{16} & \infty \\ 
\hline  
0 & 0 & 0 & \frac{3}{4}\\ 
0 & \frac{1}{2} & \frac{1}{2} & 1\\ 
0 & 1 & 1 & \frac{5}{4}
\end{array}  \right\}  }
\[
\RS1.
\]
As we see in this form, the origin $x=0$ is a LCSL, which corresponds
to a cusp of $\overline{\mbH_{+}/\Gamma_{0}(10)_{+}}$. Around this
point, for the regular solution $w_{0}(x)$ in (\ref{eq:localSol}),
we have 
\[
w_{0}(x)=\sum_{n,m\in\mbZ_{\geq0}}\frac{((n+m)!)^{4}}{(n!)^{4}(m!)^{4}}x^{n+m}.
\]
Finding all other solutions, it is straightforward to have the mirror
map (\ref{eq:mirrorMap}) as 
\[
x(q)=q-4q^{2}-6q^{3}+56q^{4}-45q^{5}-360q^{6}+894q^{7}+\cdots.
\]
Here we identify $1/x(q)$ with the Thompson series $t_{10A}$ for
the group $\Gamma_{0}(10)_{+}$ up to an additive constant; we have
\[
\frac{1}{x(q)}=t_{10C}(q)+\frac{25}{t_{10C}}+6,\,\,\,\,\,\,t_{10C}(q)=\left(\frac{\eta_{1}(t)\eta_{2}(t)}{\eta_{5}(t)\eta_{10}(t)}\right)^{2}
\]
where $\eta_{k}(t):=\eta(kt)$ with the eta function $\eta(t)=q^{\frac{1}{24}}\prod_{n\geq1}(1-q^{n})$~$(q=e^{2\pi it})$.
By the method described in Appendix \ref{subsec:AppD-GY-couplings},
we determine the Griffiths-Yukawa coupling as 
\[
C_{xx}=\frac{20}{x^{2}(1+4x)(1-16x)}
\]
and verify (\ref{eq:YuakawaYttt}) with the Gram matrix $K_{11}=20$. 
\begin{prop}
\label{prop:tauGCOV-10} The conifold and orbifold regularities uniquely
determine the parameters in $\tau_{_{\mathrm{BCOV}}}$ as $r_{0}=r_{1}=-\frac{1}{2}$
and $a=-\frac{3}{4}$. Then, we have 
\[
\tau_{_{\mathrm{BCOV}}}(t)=\left(\frac{1}{w_{0}(x)}\right)^{2}\frac{dx}{dt}dis_{0}^{r_{0}}dis_{1}^{r_{1}}x^{-1+a}=\frac{1}{\eta_{1}(t)\eta_{2}(t)\eta_{5}(t)\eta_{10}(t)},
\]
and $(\tau_{_{\mathrm{BCOV}}}(t))^{^{-1}}$ defines a BCOV cusp form
of weight two on $\mathbb{H}_{+}$ with respect to $\Gamma_{0}(10)_{+}$. 
\end{prop}

\begin{proof}
We first note that the relation (\ref{eq:detCij}) holds, i.e., we
have $C_{xx}(1+4x)(1-16x)x^{2}=20$. Hence, by the general arguments
for the conifold regularity in Subsect. \ref{subsec:Conifold-regularity},
we obtain $r_{0}=r_{1}=-\frac{1}{2}$. Furthermore, we can implement
the orbifold regularity by using the form $\tau_{_{\mathrm{BCOV}}}=\text{const.}\frac{x^{a}}{w_{0}(x)}$
and (\ref{eq:tauBCOV-orbifold}). From Riemann's scheme, we recognize
that $x=\infty$ is an orbifold point, and read the ordered solutions
around $\infty$ as $\tilde{w}_{0}(\tilde{x})=\tilde{x}^{\frac{3}{4}}(1+O(\tilde{x}))$,
$\tilde{w}_{1}(\tilde{x})=\tilde{x}(1+O(\tilde{x}))$ and $w_{2}(\tilde{x})=\tilde{x}^{\frac{5}{4}}(1+O(\tilde{x}))$
with $\tilde{x}=\frac{1}{x}$. Then, from the equation (\ref{eq:tauBCOV-orbifold}),
we see that the orbifold regularity holds only for $\tilde{a}=\frac{3}{4}$.
Since $\tilde{a}=-a$, we obtain $a=-\frac{3}{4}$.

\textcolor{black}{By Proposition \ref{prop:Modular-tau}, we know
that $(\tau_{_{\mathrm{BCOV}}}(t))^{-1}$ is a modular form of weight
one with respect to $O(\check{M}_{20}^{\perp})_{+}$, which means
weight two with respect to $\Gamma_{0}(10)_{+}$(see Appendix \ref{subsec:AppendixA-1}).
From the modular property, the equality follows by comparing the $q$-expansion
of $(\tau_{_{\mathrm{BCOV}}}(t))^{-1}$ and that of the eta product
in the r.h.s. Finally, it is clear that the eta product is a cusp
from (cf. Lemma \ref{lem:eta-cusp} below). }
\end{proof}
A finite product of the eta functions, $\prod_{m\geq1}\eta(mt)^{d_{m}}$
with $d_{m}\in\mathbb{Z}$, is called an eta product. Defining $N:=\mathrm{lcm}\left\{ m\mid d_{m}\not=0\right\} $,
we can write the eta product as 
\[
f(t)=\prod_{m|N}\eta(mt)^{d_{m}}=\prod_{m|N}\eta_{m}(t)^{d_{m}}.
\]
Eta products of this form are said of level $N$, since they are modular
forms of $\Gamma_{0}(N)$ (see \cite{Koe} for example). Since the
eta function $\eta(t)$ has no zeros in $\mathbb{H}_{+}$, the eta
products are holomorphic on $\mathbb{H}_{+}$. The only possibilities
are zero or poles at the cusps in $\overline{\mathbb{H}}_{+}\setminus\mathbb{H}_{+}$.
If an eta product $f(t)$ vanishes at all cusps, then it is called
a cusp form. The following Lemma from \cite[Cor. 2.3]{Koe} is useful
for our calculations. 
\begin{lem}
\label{lem:eta-cusp}An eta product $f(t)$ is holomorphic on $\overline{\mathbb{H}}_{+}$if
and only if the inequality 
\[
\sum_{m|N}\frac{\left({\rm gcd}(D,m)\right)^{2}}{m}d_{m}\geq0
\]
holds for all positive divisors $D$ of $N$. If the inequality holds
strictly for all $D$, then $f(t)$ defines a cusp from.
\end{lem}

For other examples in \cite{LY} of the families of $\check{M}_{2n}$-polarizable
K3 surfaces with the modular group $\Gamma_{0}(n)_{+}$, it is straightforward
to calculate the BCOV formula. In particular, if a family has an orbifold
point, we can verify the following properties: 1) the parameters $r_{i}$
and $a$ in the BCOV formula are uniquely determined by the conifold
regularity and the orbifold regularity, 2) the resulting $\tau_{_{\mathrm{BCOV}}}$
is an eta product, 3) $(\tau_{_{\mathrm{BCOV}}}(t))^{-1}$ is a cusp
form of weight two with respect to $\Gamma_{0}(n)_{+}$.

\subsection{Family of $\check{M}_{24}$-polarizable K3 surfaces}

Let us consider a family of $\check{M}_{24}$-polarizable K3 surfaces
in \cite{LY} whose modular group is $\Gamma_{0}(12)_{+}$. The Picard-Fuchs
equation of this family is given by 
\[
\left\{ \theta_{x}^{3}-2x(2\theta_{x}+1)(5\theta_{x}^{2}+5\theta_{x}+2)+64x^{2}(\theta_{x}+1)^{3}\right\} w=0,
\]
and its singularities are summarized as \def\RS2{\left\{ \small\begin{array}{cccc} 0 & \frac{1}{4} & \frac{1}{16} & \infty \\ 
\hline  
0 & 0 & 0 & 1\\ 
0 & \frac{1}{2} & \frac{1}{2} & 1\\ 
0 & 1 & 1 & 1
\end{array}  \right\}  }
\[
\RS2.
\]
As we read from this, the zeros of $dis=(1-4x)(1-16x)$ give conifold
singularities, while $x=0$ and $\infty$ are both LCSL. Since we
do not have an orbifold point in this case, the form of $\tau_{_{\mathrm{BCOV}}}(t)$
is determined only up to a parameter $a$. 
\begin{prop}
\label{prop:tauBCOV-12}With a parameter $a$, we have the following
form of $\tau_{_{\mathrm{BCOV}}}$ in terms of an eta product: 
\[
\tau_{_{\mathrm{BCOV}}}(t)=\text{const.}\frac{x^{a}}{\omega_{0}(x)}=\text{const. }x^{1+a}\left(\frac{\eta_{2}(t)\eta_{6}(t)}{\eta_{1}(t)\eta_{3}(t)\eta_{4}(t)\eta_{12}(t)}\right)^{2}.
\]
\end{prop}

We leave verifying the above equality to the reader, since calculations
are similar to the preceding example. The parameter $a$ is left undetermined
as above. However, we note that the mirror map $x(q)$ (or the inverse
of the Thompson series up to additive constant) is given by the following
eta product, 
\[
x(q)=\left(\frac{\eta_{1}(t)\eta_{3}(t)\eta_{4}(t)\eta_{12}(t)}{\eta_{2}(t)^{2}\eta_{6}(t)^{2}}\right)^{6}.
\]
Namely, $\tau_{_{\mathrm{BCOV}}}$ is determined only up to multiplying
some power of the eta product of $x(q)$. If we set $a=-\frac{5}{6}$,
for example, we can verify that $(\tau_{_{\mathrm{BCOV}}}(t))^{-1}$
is a cusp form. 

\subsection{$\tau_{_{\mathrm{BCOV}}}$ in terms of eta products }

The results in Propositions \ref{prop:tauGCOV-10},\ref{prop:tauBCOV-12}
suggest the following definition for all genus zero groups $\Gamma_{0}(n)_{+}$
listed in Appendix \ref{subsec:AppB-table-Gamma}. 
\begin{defn}
\label{def:tauBCOV-def-2} For the genus zero group $\Gamma_{0}(n)_{+}$
$(n\not=8,24,27,54$ and $56$), we define the BCOV eta product by
\begin{equation}
\eta_{_{\mathrm{BCOV}}}(t)=\left(\prod_{r|n}\eta_{r}(t)^{\pm1}\right)^{w},\label{eq:tauBCOV-etaProd}
\end{equation}
where $+1$ is taken when $(r,n/r)\not=1$ and $-1$ when $(r,n/r)=1$.
The power $w$ is determined so that the weight of $\eta_{_{\mathrm{BCOV}}}$
is two (see Appendix \ref{subsec:AppendixA-1}). For $n=8,24,27,54$
and $56$, respectively, we define $\eta_{_{\mathrm{BCOV}}}(t)$ by
\[
\begin{aligned} & \frac{\eta_{1}(t)^{3}\eta_{8}(t)^{3}}{\eta_{2}(t)\eta_{4}(t)},\,\,\frac{\eta_{1}(t)^{2}\eta_{3}(t)^{2}\eta_{8}(t)^{2}\eta_{24}(t)^{2}}{\eta_{2}(t)\eta_{4}(t)\eta_{6}(t)\eta_{12}(t)},\,\,\frac{\eta_{1}(t)^{3}\eta_{27}(t)^{3}}{\eta_{3}(t)\eta_{9}(t)},\\
 & \frac{\eta_{1}(t)^{2}\eta_{2}(t)^{2}\eta_{27}(t)^{2}\eta_{54}(t)^{2}}{\eta_{3}(t)\eta_{6}(t)\eta_{9}(t)\eta_{18}(t)}\,\,\,\,\text{and}\,\,\,\,\frac{\eta_{1}(t)^{2}\eta_{7}(t)^{2}\eta_{8}(t)^{2}\eta_{56}(t)^{2}}{\eta_{2}(t)\eta_{4}(t)\eta_{14}(t)\eta_{28}(t)}.
\end{aligned}
\]

\vskip0.3cm
\end{defn}

By calculating $\tau_{_{\mathrm{BCOV}}}$ explicitly, we verify the
following:
\begin{prop}
For all families of $\check{M}_{2n}$-polarizable K3 surfaces listed
in \cite{LY}, the following properties hold: 

\begin{myitem} 

\item[$(1)$] By imposing the conifold regularity, we obtain 
\begin{equation}
\tau_{_{\mathrm{BCOV}}}(t)=\text{const.}\frac{x^{a}}{\omega_{0}(x)}=\text{const. }\frac{x(t)^{\gamma+a}}{\eta_{_{\mathrm{BCOV}}}(t)}\label{eq:tau-eta-BCOV}
\end{equation}
where $\gamma$ is determined by the condition $x(t)^{\gamma}/\eta_{_{\mathrm{BCOV}}}(t)=1+O(q)$. 

\item[$(2)$] The orbifold regularity condition, if there is, determines
the parameter $a$ uniquely to be $a=-\gamma$. The resulting $\tau_{_{\mathrm{BCOV}}}(t)$
defines a BCOV cusp from, i.e., $(\tau_{_{\mathrm{BCOV}}}(t))^{-1}=\eta_{_{\mathrm{BCOV}}}(t)$
is a cusp form of weight two with respect to $\Gamma_{0}(n)_{+}$.

\item[$(3)$] If the orbifold regularity condition is vacant, the
parameter $a$ is left undetermined. However, in such a case, mirror
map $x(q)$ is written by an eta product, and there exists a choice
of $a$ for which $(\tau_{_{\mathrm{BCOV}}}(t))^{-1}=x^{-a-\gamma}\eta_{_{\mathrm{BCOV}}}(t)$
is a cusp from. 

\end{myitem}
\end{prop}

\begin{rem}
In the list of \cite{LY}, the case (3) occurs for $n=4,12,20$ and
$28$. We observe that, for these cases, the number of cusps of $\mathbb{H}_{+}/\Gamma_{0}(n)_{+}$
is greater than one. Let us denote by $N_{c}(n)$ the number of cusps
of $\mathbb{H}_{+}/\Gamma_{0}(n)_{+}$. We can see in \cite[Tables 2,3]{CN}
that the Thompson series of $\Gamma_{0}(n)_{+}$ with some additive
constant, $T_{n}(t)+c_{n}$, can be written by some eta product if
$N_{c}(n)\geq2$ except for $n=25,27,49,50,54$, see Appendix \ref{subsec:AppB-table-Gamma}.
When $N_{c}(n)=k$ and $n$ is not in these exceptional cases, we
can read from \cite[Tables 2,3]{CN} that there are $k-1$ different
forms of eta products. 
\end{rem}

For all genus zero groups $\Gamma_{0}(n)_{+}$, except for $n=25,27,49,50,54$,
we can verify the following property of the BCOV eta product in general. 
\begin{prop}
\label{prop:eta-prod-cusp}If $N_{c}(n)=1$, then $\eta_{_{\mathrm{BCOV}}}$
is a cusp form of weight two of $\Gamma_{0}(n)_{+}$. If $N_{c}(n)\geq2$,
there exists an expression of the Thompson series $T_{n}(t)+c_{n}$
by an eta product such that $(T_{n}(t)+c_{n})^{b}\eta_{BCOV}(t)$
is a cusp form for some $b\in\mathbb{Q}$.
\end{prop}

\begin{proof}
When $N_{c}(n)=1$, it is straightforward to verify the claimed properties
by Lemma \ref{lem:eta-cusp}. When $N_{c}(n)\geq2$, we can find the
claimed form of the Thompson series from \cite[Table 3]{CN}. In Appendix
\ref{subsec:AppendixB-Tm-product}, we list the product forms and
the values of $b$ with the resulting cusp forms. 
\end{proof}
Based on this proposition, we conjecture the form of the BCOV formula
under the conifold and orbifold regularity conditions as follows:
\begin{conjecture}
\label{conj:tauBCOV}For families of $\check{M}_{2n}$-polarizable
K3 surfaces over $\mathbb{P}^{1}$, if $N_{c}(n)=1$, then $(\tau_{_{\mathrm{BCOV}}}(t))^{-1}=\text{const. }\eta_{_{\mathrm{BCOV}}}(t)$
and this defines a BCOV cusp on $\mathbb{H}_{+}/\Gamma_{0}(n)_{+}$.
If $N_{c}(n)\geq2$ and $n\not=25,27,49,50,54$, then $(\tau_{_{\mathrm{BCOV}}}(t))^{-1}=\text{const. }x(t)^{-b}\eta_{_{\mathrm{BCOV}}}(t)$
with a parameter $b$ and this defines a cusp form for some $b\in\mathbb{Q}$.
\end{conjecture}

\subsection{K3 differential operators\label{subsec:K3-diff-eqs}}

We need to know explicit families of $\check{M}_{2n}$-polarizable
K3 surfaces to verify Conjecture \ref{conj:tauBCOV}. Here, instead
of seeking such families, we shall determine differential equations
of 3rd order by postulating that the conjecture holds. 

Let us assume that for each genus zero group $\Gamma_{0}(n)_{+}$
there exists a series $w_{0}(x)=1+O(x)$ which is a solution of a
3rd order differential equation. If we assume Conjecture \ref{conj:tauBCOV}
and also that the mirror map $x(t)=q+O(q)$ is given by the Thompson
series $T_{n}(t)+c_{n}$ with some additive constant, then by using
(\ref{eq:tau-eta-BCOV}) we have 
\begin{equation}
w_{0}(x)=x(t)^{-\gamma}\eta_{_{\mathrm{BCOV}}}(t),\label{eq:w0(x)-by-eta}
\end{equation}
where we determine the parameter $\gamma$ by requiring that the $q$-series
expansion of the r.h.s has the form $1+O(q)$. The $q$-series expansion
can be transformed to $x$-series by using the inverse relation $q=q(x)$
of 
\[
x(t)=\frac{1}{T_{n}(t)+c_{n}}.
\]
Substituting the resulting $x$-series $q=q(x)$ into the r.h.s. of
(\ref{eq:w0(x)-by-eta}), we obtain an $x$-series expansion of $w_{0}(x)$.
Now we can seek a differential operator which annihilate the power
series $w_{0}(x)$, which we shall call \textit{a K3 differential
operator of $\Gamma_{0}(n)_{+}$}. 
\begin{prop}
\label{prop:K3-diff-ops}For each genus zero group $\Gamma_{0}(n)_{+}$
listed in Table 1 of Appendix \ref{subsec:AppB-table-Gamma}, there
exists a K3 differential operator which has the following properties:\begin{myitem}

\item{$(1)$} It has a LCSL point at $x=0$, and the mirror map reproduces
the Thompson series $T_{n}(t)+c_{n}$ with an additive constant $c_{n}$. 

\item{$(2)$} The number of LCSL points coincides with the number
of cusps of $\mathbb{H}_{+}/\Gamma_{0}(n)_{+}$. 

\item{$(3)$} The regular power series solution $w_{0}(x)$ at $x=0$
satisfies the identity $w_{0}(x)=x(t)^{-\gamma}\eta_{_{\mathrm{BCOV}}}(t)$
in terms of the mirror map and the eta product. 

\end{myitem}
\end{prop}

\begin{proof}
Our proof is based on explicit calculations done for all cases. As
an example, we sketch the calculations to determine the K3 differential
operator for the group $\Gamma_{0}(36)_{+}$. First, by our definition
of $\eta_{_{\mathrm{BCOV}}}(t)$, it is straightforward to have the
following $q$-series 
\[
\eta_{_{\mathrm{BCOV}}}(t)=\left(\frac{\eta_{2}(t)\eta_{3}(t)\eta_{6}(t)\eta_{12}(t)\eta_{18}(t)}{\eta_{1}(1)\eta_{4}(t)\eta_{9}(t)\eta_{36}(t)}\right)^{4}=q^{-\frac{3}{2}}(1+4q+10q^{2}+20q^{3}+\cdots),
\]
from which we read $\gamma=-\frac{3}{2}$. We change this $q$-series
to $x$-series by substituting the inverse series of $x=1/(T_{36A}(t)+c_{36A})=q(1-2q+q^{2}+\cdots)$
where $T_{36A}(t)$ is the Thompson series given in \cite[Table 4]{CN}.
We choose $c_{36A}=2$ for the additive constant so that we have the
following product form:
\[
T_{36A}(t)+c_{36A}=\frac{\eta_{2}(t)^{2}\eta_{3}(t)^{4}\eta_{12}(t)^{4}\eta_{18}(t)^{2}}{\eta_{1}(t)^{2}\eta_{4}(t)^{2}\eta_{6}(t)^{4}\eta_{9}(t)^{2}\eta_{36}(t)^{2}}=\frac{1}{q}+2+3q+2q^{2}+\cdots,
\]
which can be found in \cite[Table 3]{CN}. Inverting the $q$-series
$x=1/(T_{36A}(t)+c_{36A})$, we obtain the following $x$-series 
\[
q=x+2x^{2}+7x^{3}+28x^{4}+125x^{5}+598x^{6}+\cdots.
\]
Substituting the above $x$-series into the r.h.s. of the relation
$w_{0}(x)=x^{-\gamma}\eta_{_{\mathrm{BCOV}}}(t)$, we obtain 
\[
w_{0}(x)=1+x+3x^{2}+15x^{3}+75x^{4}+387x^{5}+2037x^{6}+\cdots.
\]
Making the above series expansion up to order 60 and seeking a differential
operator which annihilate the series, we can find the following differential
operator: \def\ttx{\theta_x}
\[\begin{aligned}
\mathcal{D}_{36A}=&  \ttx^3
 -x (3 \ttx+1)\left(3 \ttx^2+2 \ttx+1\right) 
 -6 x^2 \ttx \left(12 \ttx^2-3 \ttx-1\right)  \\
& +2 x^3 \ttx \left(284 \ttx^2+405 \ttx+199\right) 
 +6 x^4 \ttx \left(1156 \ttx^2+75 \ttx+89\right)  \\
& -6 x^5 \ttx \left(11927 \ttx^2+10401 \ttx+4939\right) \\
& +18 x^6 \left(8968 \ttx^3+11586 \ttx^2+5960 \ttx+2553\right) \\
& +18 x^7 \left(11788 \ttx^3+14184\ttx^2-5086 \ttx-19947\right) \\ 
& -27 x^8 \left(30109 \ttx^3+44628 \ttx^2+7040 \ttx-6990\right) \\
& -27 x^9 \left(19871 \ttx^3+39147 \ttx^2+9715 \ttx+29949\right) \\
& + 486 x^{10} \left(2664 \ttx^3+4503\ttx^2+2623 \ttx+561\right) \\
& + 486 x^{11} \left(2892 \ttx^3+6453 \ttx^2+5465 \ttx+1657\right) 
 + 360126 x^{12} (\ttx+1)^3.
\end{aligned}
\] Clearly the differential equation $\mathcal{D}_{36A}w(x)=0$ has
the claimed properties (1) and (3). As for (2), we describe the singularities
of the equation $\mathcal{D}_{36A}w(x)=0$ in Riemann's $\mathcal{P}$-scheme,
\[\left\{ 
\small\begin{array}{cccccc} 
-1 & 0 & \frac{1}{3} & -1-\frac{2\sqrt{3}}{3}& -1+\frac{2\sqrt{3}}{3} & \infty \\ 
\hline  
1 & 0 & 1& 0           & 0           & 1\\ 
1 & 0 & 1& \frac{1}{2} & \frac{1}{2} & 1\\ 
1 & 0 & 1& 1           & 1           & 1
\end{array}  \right\}  
\]and we see that $\mathcal{D}_{36A}$ has four LCSL points at $x=-1,0,\frac{1}{3},\infty$,
which correspond to the four cusp points shown in the table of Appendix
\ref{subsec:AppB-table-Gamma}. By similar calculations for all other
genus zero groups $\Gamma_{0}(n)_{+}$, we can find K3 differential
operators with the claimed properties. We list the explicit forms
of K3 differential operators in Addendum to the arXiv version \cite{addendum}
of this paper. 
\end{proof}
\begin{rem}
Here we remark on the additive constant in the relation $\frac{1}{x}=T_{n}(t)+c_{n}$.
Changing the additive constant $c_{n}$ simply results in a linear
fractional transformation on $x$, i.e., an automorphism of the corresponding
K3 differential equation. In our listing \cite{addendum}, we have
chosen $c_{n}$ so that we have integer coefficients for the power
series expansion of $w_{0}(x)$ up to degree 20, if possible. However,
when $T_{n}(t)+c_{n}$ is expressed by an eta product, we have kept
$c_{n}$ to preserve this property. $\hfill\square$
\end{rem}

It should be noted that the K3 differential operators listed in \cite{addendum}
are two dimensional analogues of the so-called Calabi-Yau differential
operators considered for Calabi-Yau threefolds in \cite{AESZ}.

\vskip3cm

\section{\textbf{\textcolor{black}{\label{sec:CD-family }Clingher-Doran's
family and Igusa cusp forms $\chi_{10},\chi_{12}$}}}

As an example of higher dimensional period domains, we consider a
lattice $\check{M}=U\oplus E_{8}(-1)\oplus E_{7}(-1)$ for which we
have $\Omega_{\check{M}}\simeq\mbH_{2}$. We fix a (unique) primitive
lattice embedding $\check{M}\subset L$ and its orthogonal lattice
$\check{M}^{\perp}\simeq U^{\oplus2}\oplus\langle-2\rangle$. 

\subsection{Normal form of $\check{M}$-polarized K3 surfaces}

In \cite{CD1}, Clingher and Doran studied special quartic surfaces
$\left\{ f=0\right\} \subset\mathbb{P}^{3}$ with 
\begin{equation}
\begin{alignedat}{1}f= & \ty^{2}\tz\tw-4\tx^{3}\tz+3\alpha\tx\tz\tw^{2}+\beta\tz\tw^{3}+\gamma\tx\tz^{2}\tw-\frac{1}{2}(\delta\tz^{2}\tw^{2}+\tw^{4})\end{alignedat}
\label{eq:def-eq-CD}
\end{equation}
where $[\tx,\ty,\tz,\tw]\in\mbP^{3}$ and $\alpha,\beta,\gamma,\delta\in\mbC$
are parameters. For fixed $\alpha,\beta,\gamma,\delta,$ the minimal
resolution of $\left\{ f=0\right\} \subset\mbP^{3}$ is denoted by
$X(\alpha,\beta,\gamma,\delta)$. The quartic $f$ is called a normal
form of an $\check{M}$-polarized K3 surface since the following properties
hold \cite[Thm. 1.2]{CD1}: \begin{myitem}

\item{(1)} When $\gamma\not=0$ or $\delta\not=0$, $X(\alpha,\beta,\gamma,\delta)$
is an $\check{M}$-polarized K3 surface.

\item{(2)} For a given $\check{M}$-polarized K3 surface $X$, there
exists $(\alpha,\beta,\gamma,\delta)\in\mbC^{4}$, satisfying $\gamma\not=0$
or $\delta\not=0$, such that $X(\alpha,\beta,\gamma,\delta)$ is
isomorphic to $X$ as $\check{M}$-polarized K3 surface. 

\end{myitem}

Taking the natural isomorphism $X(\alpha,\beta,\gamma,\delta)\simeq X(t^{2}\alpha,t^{3}\beta,t^{5}\gamma,t^{6}\delta)$
$(t\in\mbC^{*})$ into account, the open variety
\[
\mathcal{M}_{\mathrm{CD}}:=\left\{ [\alpha,\beta,\gamma,\delta]\in\mathbb{W}\mbP^{3}(2,3,5,6)\mid\gamma\not=0\,\text{or}\,\delta\not=0\right\} 
\]
has been identified with a coarse moduli space of the $\check{M}$-polarized
K3 surfaces. Since we have an isomorphism $\Omega_{\check{M}}\simeq\mbH_{2}$
as shown in Subsect.$\,$\ref{subsec:Example-2-UUA1}, we have a period
map $\mathcal{P}:\mathcal{M}_{\mathrm{CD}}\rightarrow\Omega_{\check{M}}\simeq\mbH_{2}$
which is a multi-valued map. Note that, due to the isomorphism $Sp(4,\mbZ)/\left\{ \pm I_{4}\right\} \simeq SO(U^{\oplus2}\oplus\langle-2\rangle)_{+}\simeq O(U^{\oplus2}\oplus\langle-2\rangle)_{+}^{*}$,
we have 
\[
\Omega_{\check{M}}/O(U^{\oplus2}\oplus\langle-2\rangle)_{+}\simeq\Omega_{\check{M}}/O(U^{\oplus2}\oplus\langle-2\rangle)_{+}^{*}\simeq\mbH_{2}/PSp(4,\mbZ),
\]
where we use the fact that $SO(U^{\oplus2}\oplus\langle-2\rangle)_{+}\simeq O(U^{\oplus2}\oplus\langle-2\rangle)_{+}/\left\{ \pm I_{5}\right\} $
and $I_{5}$ acts trivially on $\Omega_{\check{M}}$. Namely, there
is no distinction between $\check{M}$-polarized and $\check{M}$-polarizable
K3 surfaces in this case. 
\begin{thm}[{Clingher and Doran \cite[Thm.1.5]{CD1}}]
 \label{thm:CD-Thm1.5}The inverse period map $\mathcal{P}^{-1}:\mbH_{2}\rightarrow\mathcal{M}_{\mathrm{CD}}$
is given by 
\begin{equation}
\mathcal{P}^{-1}(\tau)=\left[\mathcal{E}_{4}(\tau),\mathcal{E}_{6}(\tau),2^{12}3^{5}\chi_{10}(\tau),2^{12}3^{6}\chi_{12}(\tau)\right]\label{eq:period-inv-CD}
\end{equation}
where $\mathcal{E}_{4}$ and $\mathcal{E}_{6}$ are genus two Eisenstein
series of weight four and six, and $\chi_{10}$ and $\chi_{12}$ are
Igusa's cusp forms of weight ten and twelve, respectively. 
\end{thm}

Clingher and Doran obtained Thm. \ref{thm:CD-Thm1.5} based on an
argument which is purely geometric and does not involve any computations
of period integrals. For a general $\check{M}$-polarized K3 surface
$X$ corresponding to $[w]\in\Omega_{\check{M}}$, they construct
a new K3 surface $Y$ which is isogeny to $X$ and also carries a
canonical Kummer structure. Then the beautiful result above is derived
by identifying an elliptic fibration $Y\rightarrow\mathbb{P}^{1}$
obtained in \cite{CD2} and that for the Kummer surface representing
$Y$. We shall present explicit computations of period integrals of
the family. Together with the result (\ref{eq:period-inv-CD}), we
will calculate the BCOV formula. 

\subsection{\label{subsec:GKZ-for-CD-family}Period integrals and Picard-Fuchs
equations}

We change the normal form (\ref{eq:def-eq-CD}) to a different form
which allows us to write the period integrals explicitly, in particular,
by hypergeometric series. To this aim, let us introduce a slightly
general quartic polynomial 
\[
\begin{alignedat}{1}F(a)= & a_{0}\,\tx\ty\tz\tw+a_{1}\,\ty^{2}\tz\tw+a_{2}\,\tx^{3}\tz+a_{3}\,\tx\tz^{2}\tw+a_{4}\,\tz^{2}\tw^{2}+a_{5}\,\tw^{4}+a_{6}\,\tz\tw^{3}\\
 & \hsp{135}+a_{7}\,\ty\tz\tw^{2}+a_{8}\,\tx\tz\tw^{2}+a_{9}\,\tx^{2}\tz\tw
\end{alignedat}
\]
with $a_{0,}a_{1},\cdots,a_{9}\in\mathbb{C}$. We denote by $F_{toric}(a)$
a special form of $F(a)$ with vanishing $a_{7},a_{8},a_{9}$, i.e.,
$F_{toric}(a)=a_{0}\,\tx\ty\tz\tw+\cdots+a_{6}\,\tz\tw^{3}$. In order
to describe symmetry of $F(a)$, we introduce a (non-reductive) group
\begin{equation}
G=\left\{ \left(\begin{smallmatrix}\lambda_{1} & 0 & 0 & c_{1}\\
c_{3} & \lambda_{2} & 0 & c_{2}\\
0 & 0 & \lambda_{3} & 0\\
0 & 0 & 0 & \lambda_{4}
\end{smallmatrix}\right)\Big|\lambda_{1},...,\lambda_{4}\in\mbC^{*},c_{1},c_{2},c_{3}\in\mbC\right\} .\label{eq:non-red-G}
\end{equation}

~
\begin{prop}
\label{prop:CD-to-Toric} The following properties of $F(a)$ hold:

\begin{myitem}

\item[$(1)$] The quartic polynomial $F(a)$ preserves its form under
the linear change of coordinates $\tx,\ty,\tz,\tw$ by $G\subset GL(4,\mbC)$. 

\item[$(2)$]The normal form $f=f(\alpha,\beta,\gamma,\delta)$ is
transformed to a form $F_{toric}(a)$ by a coordinate transformation
of $G$ 
\[
\left(\begin{smallmatrix}\tx\\
\ty\\
\tz\\
\tw
\end{smallmatrix}\right)\mapsto\left(\begin{smallmatrix}\tx-\frac{\sqrt{\alpha}}{2}\tw\\
\ty-i\sqrt{6}\alpha^{\frac{1}{4}}\tx\\
\tz\\
\tw
\end{smallmatrix}\right),
\]
where the parameters $a_{i}$ in $F_{toric}(a)$ are related to $\alpha,\beta,\gamma,\delta$
by 
\begin{equation}
(a_{0},a_{1},a_{2},a_{3},\cdots,a_{6})=(-2i\sqrt{6}\alpha^{\frac{1}{4}},1,-4,\gamma,-\frac{1}{2}(\delta+\sqrt{\alpha}\gamma),-\frac{1}{2},(\beta-\alpha^{\frac{3}{2}})).\label{eq:ak-by-albeta}
\end{equation}

\end{myitem}

\vskip0.5cm
\end{prop}

Here, instead of giving a proof, we briefly describe how one can find
the above properties. For the polynomial $F(a)=\sum a_{i}\tx^{m_{i}}$,
let us denote (a translation of) its Newton polytope by $N_{F}=$
$\text{Conv\ensuremath{\left\{  m_{1}-m_{0},\cdots,m_{9}-m_{0}\right\} } in \ensuremath{\mbR^{4}}}$.
Since $F$ is homogeneous, all the vertices of $N_{F}$ lie on $H_{\mathbb{Z}}:=\left\{ m\in\mathbb{Z}^{4}\mid(1,1,1,1).m=0\right\} $.
Projecting $H_{\mathbb{Z}}$ to $\mathbb{Z}^{3}$ of the first three
components, we have an isomorphism $H_{\mathbb{Z}}\simeq\mathbb{Z}^{3}.$
Then it is easy to verify that $m_{i}-m_{0}$ ($i=1,...,6$) are the
corners of the polytope $N_{F}$, and, under the isomorphism $H_{\mathbb{Z}}\simeq\mathbb{Z}^{3}$,
the polytope is identified with 
\begin{equation}
\Delta=\text{Conv.}\left\{ \left(\begin{smallmatrix}-1\\
1\\
0
\end{smallmatrix}\right),\left(\begin{smallmatrix}2\\
-1\\
0
\end{smallmatrix}\right),\left(\begin{smallmatrix}0\\
-1\\
1
\end{smallmatrix}\right),\left(\begin{smallmatrix}-1\\
-1\\
1
\end{smallmatrix}\right),\left(\begin{smallmatrix}-1\\
-1\\
-1
\end{smallmatrix}\right),\left(\begin{smallmatrix}-1\\
-1\\
0
\end{smallmatrix}\right)\right\} \subset\mbR^{3},\label{eq:Delta-CD}
\end{equation}
which turns out to be a reflexive polytope \cite{Ba}. We also find
that the Laurent monomials in $\hat{F}(a):=F(a)/\tx\ty\tz\tw$ are
in one to one correspondence to the integral points $\Delta\cap\mbZ^{3}$;
in particular, the monomial $a_{0}1$ corresponds to the origin $\,^{t}(0,0,0)$
contained in $\Delta$, and monomials with coefficients $a_{7},a_{8},a_{9}$
correspond to integral points in the relative interior of codimension-one
faces of $\Delta$. Let $\mathbb{P}_{\Delta}$ be the projective toric
variety associated to $\Delta$, and consider a hypersurface $Z_{\hat{F}}:=\{\hat{F}(a)=0\}\subset(\mathbb{C}^{*})^{3}$
parametrized by $a_{0},a_{1},\cdots,a_{9}$. Then, due to \cite{Ba},
its closure $\overline{Z}_{\hat{F}}$ in $\mathbb{P}_{\Delta}$ defines
a smooth K3 surface after making a MPCP (maximally, projective, crepant,
and partial) resolution of $\mbP_{\Delta}$. Then the group $G$ in
(\ref{eq:non-red-G}) represents the group of automorphisms of the
toric variety $\mbP_{\Delta}$, and hence induces naturally isomorphisms
of the K3 surfaces $\{\hat{F}(a)=0\}\subset\mathbb{P}_{\Delta}$,
which give rise to those of the quartics $\left\{ F(a)=0\right\} \subset\mathbb{P}^{3}$. 

Actually, the form $F_{toric}(a)$ is a preferred form to analyze
period integral in terms of the extended GKZ system introduced in
\cite[Sect.2.3]{HLY},\cite{HKTY1}, i.e., a Gel'fand-Kapranov-Zelevinski
(GKZ) hypergeometric system \cite{GKZ1} with additional linear differential
operators representing infinitesimal group actions of $G$ in (\ref{eq:non-red-G}).
In fact, for the K3 surfaces $\left\{ F(a)=0\right\} \subset\mathbb{P}^{3}$,
the reduced form $F_{toric}(a)$ gives a canonical form of period
integral as 
\begin{equation}
w_{0}(a)=\int_{C_{0}}\frac{a_{0}}{\hat{F}_{toric}(a)}\frac{d\tx}{\tx}\frac{d\ty}{\ty}\frac{d\tz}{\tz}\frac{d\tw}{\tw}\label{eq:w0-CD}
\end{equation}
where $\hat{F}_{toric}(a)=F_{toric}(a)/\tx\ty\tz\tw$ and $C_{0}=\left\{ |\tx|=|\ty|=|\tz|=|\tw|=\ve\right\} (\ve>0)$
is a standard torus cycle in $(\mbC^{*})^{3}\subset\mbP_{\Delta}$. 

\newpage
\begin{prop}
\label{prop:Perid-map-xyz}For the period integral $w_{0}(a)$, the
following properties hold: 

\begin{myitem} 

\item{$(1)$} The period integral (\ref{eq:w0-CD}) can be evaluated
as 
\begin{equation}
\begin{aligned} & w_{0}(a)=\sum_{n,m,l\in\mbZ_{\geq0}}c(l,m,n)x^{l}y^{m}z^{n},\end{aligned}
\label{eq:CD-hypergeom}
\end{equation}
where $x:=\frac{a_{4}a_{5}}{a_{6}^{2}},y:=\frac{a_{3}^{2}a_{6}}{a_{2}a_{4}^{3}},z:=\frac{a_{1}a_{2}a_{4}}{a_{0}^{2}a_{3}}$
and the coefficients are given by 
\[
c(l,m,n)=\frac{(2n)!}{n!l!(n-m)!(3m-n)!(l-3m+n)!(m-2l)!}
\]
with understanding $k!=\Gamma(k+1)$ for negative $k$. 

\item{$(2)$} The series $w_{0}(a)$ is a hypergeometric series of
GKZ type, which is a unique regular solution of the Picard-Fuchs differential
equations (\ref{eq:exGKZ}) in Appendix \ref{sec:Appendix-CD-family}.
In particular, the origin $(x,y,z)=(0,0,0)$ is a LCSL.

\end{myitem}
\end{prop}

\begin{proof}
The residue calculation is straightforward, see \cite{BaC} for example.
As described briefly above, the group $G$ represents the group of
automorphisms of $\mbP_{\Delta}$, and acts naturally on the period
integral $w_{0}(a)$. This entails exactly the extended GKZ system
in \cite{HKTY1,HLY}. We determine the Picard-Fuchs system (\ref{eq:exGKZ})
as a minimal set of differential operators whose solutions give rise
to period integrals. Practically, we can determine them by searching
differential operators which annihilate $w_{0}(x)$. The claimed property
at the origin follows directly by constructing local solutions around
it.
\end{proof}
The extended GKZ system is a system of differential operators which
is naturally defined over a toric variety $\mathbb{P}_{SecP}$ associated
to the so-called secondary polytope \cite[2]{GKZ1}; and the affine
coordinates $x,y,z$ are for one of the affine charts of $\mathbb{P}_{SecP}$.
It should be noted that the above nice properties of the period integral
indicate that, if we change the normal form (\ref{eq:def-eq-CD})
to $F_{toric}(a)$, then we have a good family 
\begin{equation}
\check{\mathfrak{X}}\rightarrow\mathcal{M}=(\mathbb{C}^{*})^{3}\label{eq:CD-toric}
\end{equation}
to which our framework in Sect.~\ref{sec:BCOV-cusp-forms} applies
with $\overline{\mathcal{M}}=\mathbb{P}_{SecP}$. Here we define the
family $\check{\mathfrak{X}}$ from the universal family $\left\{ F_{toric}(a)=0\right\} \subset\mathbb{P}_{\Delta}\times\mathbb{C}^{7}$
by taking the quotient by the natural action of the tori $(\mathbb{C}^{*})^{3}\times\mathbb{C}^{*}$
where $(\mathbb{C}^{*})^{3}\subset\mathbb{P}_{\Delta}$ and $\mathbb{C}^{*}$
is the scaling action on $\mathbb{C}^{7}$ (see e.g. \cite{GKZ2}
for more details). 

Since the origin $(x,y,z)=(0,0,0)$ is a LCSL, we can determine the
mirror map from the local solutions around the origin following Subsect.
\ref{subsec:MS-Mmap}. In fact, we have the following results:
\begin{prop}
\label{prop:Mirror-Map-xyz}The relations (\ref{eq:ak-by-albeta})
and (\ref{eq:period-inv-CD}) determine the mirror map as 
\[
x=\frac{2^{10}3^{5}(3\chi_{12}+\chi_{10}\cE_{4}^{\frac{1}{2}})}{(\cE_{4}^{\frac{3}{2}}-\cE_{6})^{2}},\,\,y=\frac{-2\chi_{10}^{3}(\cE_{4}^{\frac{3}{2}}-\cE_{6})}{(3\chi_{12}+\chi_{10}\cE_{4}^{\frac{1}{2}})^{3}},\,\,z=\frac{-(3\chi_{12}+\chi_{10}\cE_{4}^{\frac{1}{2}})}{12\chi_{10}\cE_{4}^{\frac{1}{2}}}.
\]
Substituting these into $w_{0}(a)$, we have $w_{0}(x,y,z)=\cE_{4}(\tau)^{\frac{1}{4}}$. 
\end{prop}

\begin{proof}
Using the relations (\ref{eq:ak-by-albeta}),(\ref{eq:period-inv-CD})
and also the definitions $x=\frac{a_{4}a_{5}}{a_{6}^{2}},y=\frac{a_{3}^{2}a_{6}}{a_{2}a_{4}^{3}},z=\frac{a_{1}a_{2}a_{4}}{a_{0}^{2}a_{3}}$,
it is straightforward to verify the first half of the claims. We derive
the second half by substituting these mirror maps into the hypergeometric
series $w_{0}(x,y,z)$ in (\ref{eq:CD-hypergeom}), and verify the
equality up to sufficiently high degree in the $q$-expansions. The
claim follows since we know the modular property of $w_{0}(x,y,z)$
as described in Sect. \ref{sec:BCOV-cusp-forms}. 
\end{proof}
From the Picard-Fuchs differential equations, we can determine the
Griffiths-Yukawa couplings $C_{ij}(x,y,z)$. Also, we find that, except
the divisor $\left\{ xyz=0\right\} $, the discriminant of the Picard-Fuchs
equations consists of only one irreducible component, i.e., $dis=dis_{0}$.
Since both $C_{ij}$ and $dis_{0}$ are lengthy polynomials, we only
present some of them in Appendix \ref{subsec:AppD-GY-couplings}.
Using these, we can verify the relation (\ref{eq:detCij}) 
\[
x^{2}y^{2}z^{2}dis_{0}\det(C_{ij})=2.
\]
Combining this relation with the form of $\tau_{_{\mathrm{BCOV}}}$
in (\ref{eq:tauBCOV-Cij}), we can impose the conifold regularity,
and then obtain 
\begin{equation}
\tau_{_{\mathrm{BCOV}}}=\text{const. }\frac{x^{a_{1}}y^{a_{2}}z^{a_{3}}}{w_{0}(x)}.\label{eq:tau-conifold-abc}
\end{equation}

\subsection{\label{subsec:CD-orbifold-reg}Regularity conditions at orbifold
points}

The coordinates $x,y,z$ are those of an affine chart of a toric variety
$\mbP_{SecP}$ of the secondary polytope $SecP$ defined for the regular
triangulations of the integral polytope $\Delta$ in (\ref{eq:Delta-CD})
\cite{GKZ2}. As we summarize briefly in Appendix \ref{subsec:AppendixD-local-sol-GKZ},
there are 16 regular triangulations $T_{i}\,(i=1,...,16)$ of $\Delta$.
The affine coordinates $x,y,z$ correspond to the triangulation $T_{3}$,
and the origin is a LCSL where the hypergeometric series $w_{0}(x)$
in (\ref{eq:CD-hypergeom}) is the only regular solution of the Picard-Fuchs
differential equations and all other solutions develop logarithmic
singularities near the origin. 

As we described in Subsect.~\ref{subsec:Bound-cond} in general,
there are orbifold points among the toric boundary points in $\mbP_{SecP}$
which correspond to some triangulations $T_{k}$. Such boundary points
are characterized by the property that all local solutions of the
Picard-Fuchs equations have finite monodromy around the boundary divisors.
Precisely, since toric boundary points of $\mbP_{SecP}$ are singular
in general, we have to make a suitable toric resolutions to express
local solutions. Since the condition that all local solutions have
finite monodromy does not depend on the resolution, we simply say
that a triangulation $T_{k}$ represents an orbifold point if there
is a resolution for which all local solutions have finite monodromy. 
\begin{prop}
The orbifold points in $\mathbb{P}_{SecP}$ are given by the regular
triangulations $T_{9},T_{12}$ and $T_{16}$. The local solutions
from $T_{12}$ are isomorphic to those from $T_{16}$. 
\end{prop}

\begin{proof}
We have verified the claimed property by constructing all local solutions
explicitly for a toric resolution of the affine chart for each $T_{k}(k=1,...,16)$.
For the triangulation $T_{k}(k=9,12,16)$, the forms of local solutions
and also the affine coordinates of the resolutions are presented in
Appendix \ref{subsec:AppendixD-local-sol-GKZ}. 
\end{proof}
~~
\begin{prop}
\label{prop:BCOV-cusp-Chi10} The orbifold regularity condition from
each orbifold point in $\mbP_{SecP}$ determines the BCOV formula
(\ref{eq:tau-conifold-abc}) as follows:

\begin{myitem}

\item[$(1)$] From the regularity at the orbifold point represented
by $T_{12}$ (or $T_{16}$), we obtain 
\[
\tau_{_{\mathrm{BCOV}}}=\text{const. }\left(\frac{1}{\chi_{10}(\tau)}\right)^{\frac{1}{10}}.
\]

\item[$(2)$]From the regularity at the orbifold point represented
by $T_{9}$, we obtain 
\[
\tau_{_{\mathrm{BCOV}}}=\text{const. }\left(\frac{1}{3\chi_{12}(\tau)+\chi_{10}(\tau)\cE_{4}(\tau)^{\frac{1}{2}}}\right)^{\frac{1}{12}}.
\]

\end{myitem}
\end{prop}

\begin{proof}
We start with $\tau_{_{\mathrm{BCOV}}}$ in (\ref{eq:tau-conifold-abc})
where the conifold regularity is already satisfied. We implement the
regularity at an orbifold point by doing analytic continuation of
(\ref{eq:tau-conifold-abc}) to the corresponding boundary point.
Note that the affine coordinates $x,y,z$ in (\ref{eq:tau-conifold-abc})
are simply transformed to the affine coordinates $\tilde{x},\tilde{y},\tilde{z}$
around the orbifold point, while the local solution $w_{0}(x)$ should
be replaced by 
\[
w_{0}(x)=c_{0}\tilde{w}_{0}(\tilde{x})+c_{1}\tilde{w}_{1}(\tilde{x})+c_{2}\tilde{w}_{2}(\tilde{x})+c_{3}\tilde{w}_{3}(\tilde{x})+c_{4}\tilde{w}_{4}(\tilde{x}),
\]
where $\tilde{w}_{k}(\tilde{x})$ are the local solutions around the
orbifold point. As shown in Appendix \ref{subsec:AppendixD-local-sol-GKZ},
we fix the solution $\tilde{w}_{0}(\tilde{x})$ so that $\tilde{w}_{k}(\tilde{x})/\tilde{w}_{0}(\tilde{x})$
are regular. Recall that the period integrals around the LCSL are
arranged so that $\mathcal{P}(x)=[w_{0}(x),w_{2}(x),w_{1}^{(1)}(x),$
$\cdots,w_{r}^{(1)}(x)]$ is a point on $\Omega(M^{\perp})$, hence
$w_{0}(x)\not=0$ holds for the period point $\mathcal{P}(x)$ corresponding
to a smooth K3 surface; which is the case for the orbifold points.
Therefore we have $c_{0}\not=0$ and  
\[
\tau_{_{\mathrm{BCOV}}}=\text{const.}\frac{\tilde{x}^{a_{1}}\tilde{y}^{a_{2}}\tilde{z}^{a_{3}}}{c_{0}\tilde{w}_{0}(\tilde{x})(1+O(\tilde{x}))},
\]
where $x^{a_{1}}y^{a_{2}}x^{a_{3}}=\tilde{x}^{\tilde{a}_{1}}\tilde{y}^{\tilde{a}_{2}}\tilde{z}^{\tilde{a}_{3}}$.
In this form, the orbifold regularity condition is equivalent to the
local behavior $\tau_{_{\mathrm{BCOV}}}=\text{const.}(1+O(\tilde{x}))$. 

For the orbifold point $T_{12}$, we have 
\[
\tilde{x}=\frac{1}{yz^{2}},\,\,\tilde{y}=\frac{1}{xz},\,\,\tilde{z}=z
\]
and also $\tilde{w}_{0}(\tilde{x})=\tilde{x}^{\frac{1}{5}}\tilde{y}^{\frac{1}{10}}(1+\cdots)$
with $\tilde{x}^{\frac{1}{5}}\tilde{y}^{\frac{1}{10}}=x^{-\frac{1}{10}}y^{-\frac{1}{5}}z^{-\frac{1}{2}}$.
Using these data, the parameters $a_{1},a_{2},a_{3}$ in (\ref{eq:tau-conifold-abc})
are determined uniquely to $(a_{1},a_{2},a_{3})=(-\frac{1}{10},-\frac{1}{5},-\frac{1}{2})$.
Then the result follows from $\tau_{_{\mathrm{BCOV}}}$ in (\ref{eq:tau-conifold-abc})
and formulas in Proposition \ref{prop:Mirror-Map-xyz}. The arguments
are exactly the same for $T_{16}$ and $T_{9}$ with the data given
in Appendix \ref{subsec:AppendixD-local-sol-GKZ}. The parameters
in (\ref{eq:tau-conifold-abc}) are determined to be $(a_{1},a_{2},a_{3})=(-\frac{1}{10},-\frac{1}{5},-\frac{1}{2})$
and $(-\frac{1}{12},-\frac{1}{6},-\frac{1}{2})$ for $T_{16}$ and
$T_{9}$, respectively. 
\end{proof}
\begin{rem}
The modular form $\chi_{10}(\tau)$ is a well-studied cusp from on
$\mathbb{H}_{2}$ of the modular group $Sp(2,\mathbb{Z})$ \cite{Ig,Bo}.
This also appeared in the study of heterotic-type II string duality
\cite{Ka1} (see also \cite{GN1}), which is close to our subject
here. In any case, the BCOV formula defines a BCOV cusp from on $\Omega_{\check{M}}\simeq\mathbb{H}_{2}$
when we implement the orbifold regularity at $T_{12}$ or $T_{16}$.
The modular form $3\chi_{12}(\tau)+\chi_{10}(\tau)\cE_{4}(\tau)^{\frac{1}{2}}$
from (2) is a cusp form on $\mathbb{H}_{2}$. However, as we see in
the fractional power of $\mathcal{E}_{4}^{\frac{1}{2}}$, it is a
cusp form for a smaller subgroup $\mathcal{G}_{m}\subset O(U^{\oplus2}\oplus\langle-2\rangle)_{+}$
where $\mathcal{G}_{m}$ is the monodromy group of our family $\check{\mathfrak{X}}\rightarrow\mathbb{P}_{SecP}$
in (\ref{eq:CD-toric}). We leave the problem of determining $\mathcal{G}_{m}$
in future study. 
\end{rem}

\subsection{Family of $U\oplus E_{8}(-1)^{\oplus2}$-polarized K3 surfaces }

Over the divisor $\left\{ \gamma=0\right\} \subset\mathbb{P}(2,3,5,6)$,
the lattice polarization of Clingher-Doran's family extends from $\check{M}=U\oplus E_{8}(-1)\oplus E_{7}(-1)$
to $\check{M}^{e}:=U\oplus E_{8}(-1)^{\oplus2}$ (see \cite[Sect.3.3]{CD1}).
Under the period map $\mathcal{P}:\mathcal{M}_{\mathrm{CD}}\rightarrow\mathbb{H}_{2}$,
this divisor is mapped to the Humbert surface given as the vanishing
locus of the cusp form $\chi_{10}$. We can observe this fact in the
isomorphism $\Omega_{\check{M}^{e}}=\Omega(U\oplus U)\simeq\mathbb{H}_{+}\times\mathbb{H}_{+}$
where we use $(\check{M}^{e})^{\perp}=U\oplus U$ (see Example 1 in
Subsect. \ref{subsec:Example-1-UU}). Since the dimension of the deformation
space reduces while the lattice polarization extends, we call the
family over the divisor a reduced family. Geometric aspects as well
as modular property of the reduced family were studied in detail in
\cite{CD2}. Also, since all calculations are parallel to those for
Clingher-Doran's family presented in the preceding subsections, we
only sketch them below for this reduced family. 

\subsubsection{Period integrals and mirror map}

We continue to work on the family $\check{\mathfrak{X}}\rightarrow\mathbb{P}_{SecP}$
obtained from (\ref{eq:CD-toric}). Then since $a_{3}=\gamma$ by
Proposition \ref{prop:CD-to-Toric}, the divisor $\{\gamma=0\}$ corresponds
to a toric divisor $\{a_{3}=0\}$ in $\mathbb{P}_{SecP}$. To describe
the reduced family, let $\Delta^{red}$ be the integral polytope $\Delta$
in (\ref{eq:Delta-CD}) with the third vertex $\,^{t}(0,-1,1)$ being
removed. It turns out that this reduced polytope is still reflexive,
and all toric constructions in the preceding subsections work with
\[
F_{toric}^{red}(a):=F_{toric}(a)\vert_{a_{3}=0}.
\]
In particular, the period integral (\ref{eq:w0-CD}) defines the corresponding
period integral by $w_{0}^{red}(a)=w_{0}(a)\vert_{a_{3}=0}$. Also,
we can describe the reduced family $\check{\mathfrak{X}}^{red}\rightarrow\mathbb{P}_{SecP_{red}}$
as the quotient of $\left\{ F_{toric}^{red}(a)=0\right\} \subset\mathbb{P}_{\Delta^{red}}\times\mathbb{C}^{6}$
by the action $(\mathbb{C}^{*})^{3}\times\mathbb{C}^{*}$ with the
secondary polytope $SecP_{red}$ for $\Delta^{red}$. Actually, this
reduced family is known as a mirror family of a hypersurface of degree
twelve in the weighted projective space $\mathbb{P}(6,4,1,1)$. 
\begin{prop}
For the reduced family, we have the following results:

\begin{myitem}

\item[$(1)$] The period integral $w_{0}^{red}(a)$ is given by 
\[
w_{0}^{red}(a)=\sum_{m,n\geq0}\frac{(6n)!}{(3n)!(2n)!(n-2m)!(n!)^{2}}u^{n}v^{m}
\]
with $u:=\frac{a_{1}^{3}a_{2}^{2}a_{6}}{a_{0}^{6}}$, $v:=\frac{a_{4}a_{5}}{a_{6}^{6}}$. 

\item[$(2)$] The origin $(u,v)=(0,0)$ is a LCSL. 

\item[$(3)$] Around the LCSL in (2), we have the mirror map
\begin{equation}
u=\frac{1}{864}\frac{(E_{4}\widetilde{E}_{4})^{\frac{3}{2}}-E_{6}\widetilde{E}_{6}}{(E_{4}\widetilde{E}_{4})^{\frac{3}{2}}},\,\,v=(864)^{2}\frac{\eta^{24}\widetilde{\eta}^{24}}{(E_{4}\widetilde{E}_{4})^{\frac{3}{2}}-E_{6}\widetilde{E}_{6}},\label{eq:reduced-mirror-map}
\end{equation}
where we define $E_{k}:=E_{k}(\tau_{11})$ , $\widetilde{E}_{k}:=E_{k}(\tau_{22})$
with the elliptic Eisenstein series $E_{k}$, and set $\eta:=\eta(\tau_{11})$
and $\widetilde{\eta}:=\eta(\tau_{22})$. Substituting the mirror
map into $w_{0}^{red}(a)$, we also have $w_{0}^{red}(a)=E_{4}^{\frac{1}{4}}\,\widetilde{E}_{4}^{\frac{1}{4}}.$

\end{myitem}
\end{prop}

\begin{proof}
It is straightforward to have (1) from (\ref{eq:CD-hypergeom}). The
property (2) follows by solving the Picard-Fuchs equations of the
reduced family, which we leave to the reader. The mirror map and the
relation $w_{0}^{red}(a)=E_{4}^{\frac{1}{4}}\,\widetilde{E}_{4}^{\frac{1}{4}}$
follow from Proposition \ref{prop:Mirror-Map-xyz}. In the derivation,
we use the relations 
\[
\chi_{10}(\tau)\rightarrow0,\,\,\,\,\chi_{12}(\tau)\rightarrow\eta^{24}\widetilde{\eta}^{24},\,\,\,\,\mathcal{E}_{4}(\tau)\rightarrow E_{4}\widetilde{E}_{4},\,\,\,\,\mathcal{E}_{6}(\tau)\rightarrow E_{6}\widetilde{E}_{6}
\]
under the limit $\tau_{12}\rightarrow0$ where $\tau=\left(\begin{smallmatrix}\tau_{11} & \tau_{12}\\
\tau_{21} & \tau_{22}
\end{smallmatrix}\right)\in\mathbb{H}_{2}$. 
\end{proof}

\subsubsection{BCOV cusp form }

We sketch the computation of the BCOV cusp from of the reduced family.
Since we can verify the relation (\ref{eq:detCij}) for the Griffiths-Yukawa
couplings, we start from 
\[
\tau_{_{\mathrm{BCOV}}}^{red}(t)=\text{const.}\frac{u^{a}v^{b}}{w_{0}^{red}(u,v)}.
\]
The secondary polytope $SecP_{red}$ has four vertices corresponding
to regular triangulations of $\Delta^{red}$. Making local solutions,
we can identify one of the four triangulations with an orbifold point,
which we denote by $T_{4}^{red}$. Imposing the orbifold regularity
on $\tau_{_{\mathrm{BCOV}}}^{red}(t)$ at $T_{4}^{red}$ in $\mathbb{P}_{SecP_{red}}$,
we find that $(a,b)=(-\frac{1}{6},-\frac{1}{12})$. Then, using (\ref{eq:reduced-mirror-map})
and the relation $w_{0}^{red}(u,v)=E_{4}^{\frac{1}{4}}\,\widetilde{E}_{4}^{\frac{1}{4}}$,
we obtain 
\[
\tau_{_{\mathrm{BCOV}}}^{red}=\text{const.}\left(\frac{1}{\eta^{24}(\tau_{11})\eta^{24}(\tau_{22})}\right)^{\frac{1}{12}}=\mathrm{const.}\frac{1}{\eta^{2}\tilde{\eta}^{2}}.
\]
The inverse power $(\tau_{_{\mathrm{BCOV}}}^{red})^{-1}$ clearly
defines a cusp form on $\mathbb{H}_{+}\times\mathbb{H}_{+}$, hence
we obtain a BCOV cusp form from the BCOV formula. 
\begin{rem}
As we notice, the BCOV cusp form $(\tau_{_{\mathrm{BCOV}}}^{red})^{-1}$
directly follows as the limit of $(\tau_{_{\mathrm{BCOV}}})^{-1}$
for the case (2) in Proposition \ref{prop:BCOV-cusp-Chi10} where
we impose the orbifold regularity at the orbifold point corresponding
to the triangulation $T_{9}$ in $\mathbb{P}_{SecP}$. Note that one
of the faces of $SecP$ can be identified with $SecP_{red}$ since
$\mathbb{P}_{SecP_{red}}$ is a toric divisor on $\mathbb{P}_{SecP}$.
The triangulation $T_{4}^{red}$ corresponds to $T_{9}$ under this
relation $SecP_{red}\prec SecP$. This explains that $(\tau_{_{\mathrm{BCOV}}}^{red})^{-1}$
follows from $(\tau_{_{\mathrm{BCOV}}})^{-1}$ for the case (2) in
Proposition \ref{prop:BCOV-cusp-Chi10}. On the other hand, the cusp
from $(\tau_{_{\mathrm{BCOV}}})^{-1}=\text{const.}(\chi_{10})^{\frac{1}{10}}$
for the case (1) does not reduce to a non-trivial BCOV cusp from.
This can be explained by the fact that the orbifold points corresponding
to $T_{10}$ and $T_{16}$ do not reduce to any point on the divisor
$\mathbb{P}_{SecP_{red}}$ of $\mathbb{P}_{SecP}$. Reversing the
arguments, one can say that $(\tau_{_{\mathrm{BCOV}}})^{-1}=\text{const.}(\chi_{10})^{\frac{1}{10}}$
is a genuine BCOV cusp form which arises for the $\check{M}$-polarized
K3 surfaces. 
\end{rem}

~

\section{\textbf{\label{sec:Weil-Peterson}Weil-Petersson geometry and BCOV
torsion}}

~

Here we collect some formulas which relate the BCOV formula with the
Weil-Petersson K\"ahler geometry on the moduli space of K3 surfaces.
Related geometry will be studied more in detail in \cite{GHK}.

\subsection{The Weil-Petersson metric}

Let $\check{\mathfrak{X}}\rightarrow\mathcal{M}$ be a family of $\check{M}$-polarized
(or polarizable) K3 surfaces. As in Sect. \ref{sec:BCOV-cusp-forms},
we describe the period map $\mathcal{P}:\mathcal{M}\rightarrow\Omega_{\check{M}}$
by $\mathcal{P}(x)=[w(x)]\in\Omega_{\check{M}}$ in terms of the period
integrals 
\[
w(x)=\big(w_{0}(x),w^{(2)}(x),w_{1}^{(1)}(x),\cdots,w_{r}^{(1)}(x)\big)
\]
which satisfy the quadratic relation $(w,w)=2w_{0}w^{(2)}+(w^{(1)},w^{(1)})_{M}=0$.
Here by definition of $\Omega_{\check{M}}=\Omega(\check{M}^{\perp})$
and the assumption $\check{M}^{\perp}=U\oplus M$, we use the bilinear
form $(\,,\,)$ of the lattice $U\oplus M$, and denote by $(\,,\,)_{M}$
its restriction to $M$. 

Due to \cite{Tian}, we introduce the Weil-Petersson metric on $\mathcal{M}$
in terms of the following K\"ahler potential 
\[
\mathcal{K}(x,\bar{x})=-\log\big(w(x),\overline{w(x)}\big),
\]
and denote the K\"ahler metric by $g_{i\bar{j}}:=\frac{\partial\;}{\partial x_{i}}\frac{\partial\;}{\partial\overline{x_{j}}}\mathcal{K}(x,\bar{x})$.
In order to write the corresponding metric on the tube domain $M\otimes\mathbb{R}+iC_{M}$,
we set $(t)_{M}^{2}:=(t,t)_{M}$ for the tube domain coordinates $t=(t_{a}):=(t_{1},\cdots,t_{r})$.
In what follows, we use the indices $a,b,c,...$ for the coordinates
$t=(t_{a})$, while $i,j,k,...$ for the coordinates $x=(x_{1},...,x_{r})$. 
\begin{lem}
Using the definition $t_{a}=\frac{w_{a}^{(1)}(x)}{w_{0}(x)}$, we
have 
\end{lem}

\begin{equation}
\big(w(x),\overline{w(x)}\big)=-\frac{1}{2}w_{0}(x)\overline{w_{0}(x)}(t-\bar{t})_{M}^{2}.\label{eq:norm-(t-bt)^2}
\end{equation}

\begin{proof}
From the quadratic relation $(w,w)=0$, we have $w(x)=w_{0}(x)(1,-\frac{1}{2}(t)_{M}^{2},t_{1},$
$\cdots,t_{r})$. The equality follows immediately from this.
\end{proof}
With the above equality, we may write the K\"ahler metric on the
tube domain as $g_{a\bar{b}}=-\frac{\partial\;}{\partial t_{a}}\frac{\partial\;}{\partial\bar{t}_{b}}\log\big(-(t-\bar{t})_{M}^{2}\big)$.
Here it should be noted that $(t-\bar{t})_{M}^{2}<0$. The K\"ahler
metric $g_{i\bar{j}}$ (or $g_{a\bar{b}}$) defines the so-called
special K\"ahler geometry on the moduli spaces, and this provides
a profound connection to the $tt^{*}$-geometry which will be the
subject in \cite{GHK}. 

\subsection{BCOV torsion}

Let us note that the equation (\ref{eq:norm-(t-bt)^2}) is monodromy
invariant, hence it is globally defined over $\mathcal{M}$. Also,
the period integral $w_{0}(x)$ is regarded as a section of the Hodge
bundle over $\mathcal{M}$, and transforms as automorphic form of
weight one as we saw in Sect. \ref{sec:BCOV-cusp-forms}. By Proposition
\ref{prop:Modular-tau}, we know that the BCOV formula $(\tau_{_{\mathrm{BCOV}}}(t))^{-1}$
has weight one with some multiplier system $v(g)$. From these properties,
we naturally come to the following definition:
\begin{defn}[\textbf{BCOV torsion}]
 As a $C^{\infty}$ function over $\Omega_{\check{M}}\simeq M\otimes\mathbb{R}+iC_{M}$,
or its quotient $\Omega_{\check{M}}/\mathcal{G}_{m}$ by the monodromy
group $\mathcal{G}_{m}$, we define the BCOV torsion by $\log\mathcal{T}_{\mathrm{BCOV}}$
with 
\begin{equation}
\mathcal{T}_{\mathrm{BCOV}}(t,\bar{t}):=-(t-\bar{t})_{M}^{2}(\tau_{_{\mathrm{BCOV}}}(t))^{-1}\overline{(\tau_{_{\mathrm{BCOV}}}(t))^{-1}}.\label{eq:T-BCOV}
\end{equation}

In the original work by BCOV \cite{BCOV1}, it was conjectured that
the BCOV torsion coincides with the so-called Ray-Singer analytic
torsion \cite{RSinger} for Calabi-Yau manifolds in general. We refer
\cite{JT,Yo} for some results on the analytic torsions of K3 and
Enriques surfaces. Here it should be useful to see how these torsions
are related for the case of elliptic curve: For a family of elliptic
curves, the Ray-Singer analytic torsion is defined as the (modified)
determinant of Laplace operator \cite{RSinger} \cite[eq.(5.17)]{Atiyah}
by 
\begin{equation}
\mathcal{T}_{\mathrm{RS}}^{E}=\text{det}'\partial_{\tau}\overline{\partial_{\tau}}=\log\left\{ (-i)(\tau-\overline{\tau})\eta(\tau)^{2}\overline{\eta(\tau)}^{2}\right\} ,\label{eq:T-ellp}
\end{equation}
where $\det'$ means omitting zero eigenvalues. In this form, the
Ray-Singer analytic torsion represents the Quillen's norm of determinant
line bundle \cite{Quillen}. As a family of elliptic curves, we take
the Weierstrass normal form of the family. For this family, similar
calculations as in Sect. \ref{sec:CD-family } apply. We have the
BCOV formula $\tau_{_{\mathrm{BCOV}}}$ similar to (\ref{eq:tauBCOV-x})
and obtain $(\tau_{_{\mathrm{BCOV}}}(\tau))^{-1}=\eta(\tau)^{2}$
by calculating the mirror map. Also, we have the K\"ahler potential
of the Poincar\'e metric on $\mathbb{H}_{+}$ as $\mathcal{K}(\tau,\bar{\tau})=-\log\big(w_{0}(x)\overline{w_{0}(x)}(-i)(\tau-\bar{\tau})\big)$
with $\tau=\frac{I_{1}(x)}{I_{0}(x)}$ and $I_{0}(x)=w_{0}(x)$ (see
Appendix \ref{sec:AppendixD-Elliptic} for details). Then, replacing
the ``norm'' factor $-(t-\bar{t})_{M}^{2}$ in (\ref{eq:T-BCOV})
with $(-i)(\tau-\bar{\tau})$, we verify the relation $\mathcal{T}_{\mathrm{RS}}^{E}=\log\mathcal{T}_{\mathrm{BCOV}}^{E}$. 
\end{defn}

\subsection{Quasi-automorphic forms}

The function $\mathcal{T}_{\mathrm{BCOV}}^{E}$ is a $C^{\infty}$
function on the fundamental domain $\mathbb{H}_{+}/\mathrm{PSL(2,\mathbb{Z})}$.
Taking the derivative of $\log\mathcal{T}_{\mathrm{BCOV}}^{E}$, we
have
\begin{equation}
\frac{1}{2\pi i}\frac{\partial\;}{\partial\tau}\log\mathcal{T}_{\mathrm{BCOV}}^{E}=\frac{1}{\tau-\bar{\tau}}+\frac{1}{12}E_{2}(t)=:\frac{1}{12}E_{2}^{*}(\tau,\bar{\tau}),\label{eq:E2*}
\end{equation}
which obviously has modular weight two, since $\frac{\partial\;}{\partial\tau'}=(c\tau+d)^{2}\frac{\partial\;}{\partial\tau}$
for the modular transform $\tau'=\frac{a\tau+b}{c\tau+d}$, but is
not holomorphic. Formally, the $\bar{\tau}\rightarrow\infty$ limit
of $E_{2}^{*}(\tau,\bar{\tau})$ gives the standard elliptic Eisenstein
series $E_{2}(\tau$). Since $E_{2}(\tau)$ itself does not have a
modular property, it is called a quasi-modular form (see \cite{KZ}).
For K3 surfaces, corresponding to (\ref{eq:E2*}), we introduce \textit{propagator
function} by
\[
S^{a}=\sum_{b}K^{ab}\frac{\partial\;}{\partial t_{b}}\log\mathcal{T}_{\mathrm{BCOV}}(t,\bar{t})
\]
where $K^{ab}$ is the matrix components of $(K_{ab})^{-1}$ for the
Gram matrix $(K_{ab})$ of the lattice $M$ (see (\ref{eq:localSol})). 
\begin{prop}
The propagator function transforms under the action $t\mapsto\tilde{t}:=g\cdot t$
of the monodromy group $g\in\mathcal{G}_{m}\subset O(U\oplus M)_{+}$
as 
\[
S^{a}(\tilde{t},\overline{\tilde{t}})=D(g,t)^{2}\sum_{b}\frac{\partial\tilde{t}_{a}}{\partial t_{b}}S^{b}(t,\bar{t})
\]
where $D(g,t)$ is the automorphic factor introduced in (\ref{eq:g-action-D(g,t)}). 
\end{prop}

\begin{proof}
Since $\mathcal{T}_{\mathrm{BCOV}}(t,\bar{t})$ is monodromy invariant,
it is sufficient to prove 
\begin{equation}
\sum_{b}K^{ab}\frac{\partial\;}{\partial\tilde{t}_{b}}=D(g,t)^{2}\sum_{b,c}\frac{\partial\tilde{t}_{a}}{\partial t_{b}}K^{bc}\frac{\partial\;}{\partial t_{c}},\label{eq:propa-pd}
\end{equation}
which generalizes the relation $\frac{\partial\;}{\partial\tau'}=(c\tau+d)^{2}\frac{\partial\;}{\partial\tau}$
for (\ref{eq:E2*}). Recall the mirror symmetry relation (\ref{eq:YuakawaYttt})
\begin{equation}
\frac{1}{w_{0}(x(t))^{2}}\sum_{i,j}C_{ij}(x(t))\frac{\partial x_{i}(t)}{\partial t_{a}}\frac{\partial x_{j}(t)}{\partial t_{b}}=K_{ab}\label{eq:MS-relation}
\end{equation}
with the mirror map $x(t)=x(q)$. Under the action $t\mapsto g\cdot t$
of the monodromy group $g\in\mathcal{G}_{m}$, the period integral
$w_{0}(x(t))$ transforms as a section of the Hodge bundle, i.e.,
$w_{0}(x(g\cdot t))=D(g,t)w_{0}(x(t))$, while we have $C_{ij}(x(g\cdot t))=C_{ij}(x(t))$
for the Griffiths-Yukawa couplings since these are rational functions
of $x=x(t)$ which is invariant under the monodromy, i.e., $x(g\cdot t)=x(t)$.
Since the relation (\ref{eq:MS-relation}) is an identity on $\Omega_{\check{M}}$,
which is invariant under the monodromy group $\mathcal{G}_{m}$, we
have 
\[
\begin{aligned} & K_{ab} & = & \frac{1}{w_{0}(x(g\cdot t))^{2}}\sum_{i,j}C_{ij}(x(g\cdot t))\frac{\partial x_{i}(g\cdot t)}{\partial(g\cdot t)_{a}}\frac{\partial x_{j}(g\cdot t)}{\partial(g\cdot t)_{b}}\\
 &  & = & \frac{1}{D(g,t)^{2}}\frac{1}{w_{0}(x(t))^{2}}\sum_{i,j,c,d}C_{ij}(x(t))\frac{\partial x_{i}(t)}{\partial t_{c}}\frac{\partial x_{j}(t)}{\partial t_{d}}\frac{\partial t_{c}}{\partial(g\cdot t)_{a}}\frac{\partial t_{d}}{\partial(g\cdot t)_{b}}\\
 &  & = & \frac{1}{D(g,t)^{2}}\sum_{c,d}K_{cd}\frac{\partial t_{c}}{\partial\tilde{t}_{a}}\frac{\partial t_{d}}{\partial\tilde{t}_{b}},
\end{aligned}
\]
where we use the relation $\tilde{t}_{a}=(g\cdot t)_{a}$ . Calculating
the inverse $(K_{ab})^{-1}$, we obtain 
\[
K^{ab}=D(g,t)^{2}\sum_{c,d}K^{cd}\frac{\partial\tilde{t}_{a}}{\partial t_{c}}\frac{\partial\tilde{t}_{b}}{\partial t_{d}}.
\]
The relation (\ref{eq:propa-pd}) follows immediately from this.
\end{proof}
We can write the propagator function explicitly as 
\[
S^{a}(t,\bar{t})=\frac{t_{a}-\bar{t}_{a}}{(t-\bar{t})_{M}^{2}}-\sum_{b}K^{ab}\frac{\partial\;}{\partial t_{b}}\log\tau_{_{\mathrm{BCOV}}}(t).
\]
In this form, it is clear that we have a (vector-valued) quasi-automorphic
form of weight two from the $\bar{t}\rightarrow\infty$ limit. More
geometric aspects of the propagator function will be discussed in
\cite{GHK}.

~

\vskip2cm~

\selectlanguage{english}%

\appendix
\renewcommand{\themyparagraph}{{\Alph{section}.\arabic{subsection}.\alph{myparagraph}}}
\selectlanguage{american}%

\section{\textbf{Isomorphisms of modular groups}}

\subsection{\label{subsec:AppendixA-1}Homomorphism $\Gamma_{0}(n)_{+}\rightarrow O(U\oplus\langle2n\rangle)$}

For reader's convenience, we present an explicit form of group homomorphism
$\Gamma_{0}(n)_{+}\rightarrow O(U\oplus\langle2n\rangle)$ which gives
the isomorphisms (\ref{eq:Gamma0N-iso}) for $n>1$. Its derivation
can be found in \cite[Thm.7]{Do} for example. Translating the convention
in \cite{Do} to ours, we obtain the group homomorphism $\psi:\Gamma_{0}(n)_{+}\rightarrow O(U\oplus\langle2n\rangle)$
by 

\[
\psi\left(\left(\begin{matrix}a & b\\
c & d
\end{matrix}\right)\right)=\left(\begin{matrix}d^{2} & -\frac{c^{2}}{n} & 2cd\\
-nb^{2} & a^{2} & -2nab\\
bd & -\frac{ac}{n} & ad+bc
\end{matrix}\right).
\]
It is an easy exercise to verify $\,^{t}\psi(A)\Sigma_{2n}\psi(A)=\Sigma_{2n}$
with $\Sigma_{2n}=\left(\begin{smallmatrix}0 & 1 & 0\\
1 & 0 & 0\\
0 & 0 & 2n
\end{smallmatrix}\right)$ and also $\psi(AB)=\psi(A)\psi(B)$ holds for $A,B\in\Gamma_{0}(n)_{+}$.
Following our convention in Subsect. \ref{subsec:AutForm}, we define
$\phi:\mathbb{H}_{+}\rightarrow\Omega_{\check{M}_{2n}}=\Omega(U\oplus\langle2n\rangle)$
by $t\mapsto\phi(t):=[1,-nt^{2},t]$. Then writing the standard action
of $A=\left(\begin{smallmatrix}a & b\\
c & d
\end{smallmatrix}\right)$ on $t\in\mathbb{H}_{+}$ by $A*t=\frac{at+b}{ct+d}$, we can verify
\[
\psi(A)\cdot\phi(t)=[(ct+d)^{2},-n(at+b)^{2},(ct+d)(at+b)]=\phi(A*t),
\]
where $\psi(A)\cdot\phi(t)$ represents the natural linear action
of $\psi(A)$ on $\phi(t)$. This indicates that we should read $D(A,t)=(ct+d)^{2}$
for the automorphic factor $D(A,t)$ in (\ref{eq:Def-autoForm}).
This means that the weight of an automorphic form of the group $O(U\oplus\langle2n\rangle)$
should be doubled when we read it as a modular form of $\Gamma_{0}(n)_{+}$. 

\subsection{\label{subsec:AppendixA-2}Isomorphism $\Psi:Sp(4,\protect\mbZ)/\left\{ \pm I_{4}\right\} \simeq SO(U^{\oplus2}\oplus\langle-2\rangle)_{+}$}

Let $V_{\mbZ}$ be a $\mbZ$-module generated by $e_{1},e_{2},e_{3},e_{4}$.
We consider the natural action of $SL(4,\mbZ)$ on $V_{\mbZ}$. Let
$\wedge^{2}V_{\mbZ}$ be the exterior power and take a basis $e_{i}\wedge e_{j}\,(1\leq i<j\leq4)$
of it. We define a non-degenerate symmetric bilinear form $(\,,\,)$
on $\wedge^{2}V_{\mbZ}$ by 
\[
X\wedge Y=:(X,Y)\,e_{1}\wedge e_{2}\wedge e_{3}\wedge e_{4}
\]
for $X=\sum x_{ij}e_{i}\wedge e_{j}$ and $Y=\sum y_{ij}e_{i}\wedge e_{j}$.
We denote the natural action of $SL(4,\mbZ)$ on $\wedge^{2}V$ by
$\rho_{\wedge}(g)(X):=gX\,^{t}g$. Then it is easy to verify that
$(\rho_{\wedge}(g)(X),\rho_{\wedge}(g)(Y))=(\det g)\,(X,Y)=(X,Y)$
for $g\in SL(4,\mbZ)$. 

Consider an element $W:=e_{1}\wedge e_{3}-e_{2}\wedge e_{4}$ and
define a group 
\[
\tilde{\Gamma}=\left\{ g\in SL(4,\mbZ)\mid\rho_{\wedge}(g)(W)=W\right\} .
\]
Due to the invariance $(\rho_{\wedge}(g)(X),\rho_{\wedge}(g)(Y))=(X,Y)$,
$\tilde{\Gamma}$ acts on the sub-module $W^{\perp}\subset\wedge^{2}V_{\mbZ}$
which is generated by 
\[
E_{1},F_{1},E_{2},F_{2},A:=e_{1}\wedge e_{2},\,e_{3}\wedge e_{4},\,e_{2}\wedge e_{3},\,e_{1}\wedge e_{4},\,e_{1}\wedge e_{3}+e_{2}\wedge e_{4}.
\]
It is easy to verify that non-vanishing pairings of these generators
are $(E_{1},F_{1})=(E_{2},F_{2})=1,(A,A)=-2$. Hence $W^{\perp}\simeq U^{\oplus2}\oplus\langle-2\rangle$,
and $\tilde{\Gamma}$ acts on $W^{\perp}$ preserving the bilinear
form, i.e., $\rho_{\wedge}(\tilde{\Gamma})\subset O(W^{\perp})$.
Now, the isomorphism $\rho_{\wedge}:\tilde{\Gamma}/\left\{ \pm I_{4}\right\} \stackrel{\sim}{\rightarrow}SO(U^{\oplus2}\oplus\langle-2\rangle)_{+}$
follows by studying the homomorphism $\rho_{\wedge}:\tilde{\Gamma}\rightarrow O(W^{\perp})$
(see e.g. \cite{GH} for more details). 

We can write the condition for $g\in\tilde{\Gamma}$ explicitly as
\[
g\left(\begin{smallmatrix}0 & 0 & 1 & 0\\
0 & 0 & 0 & -1\\
-1 & 0 & 0 & 0\\
0 & 1 & 0 & 0
\end{smallmatrix}\right)\,^{t}g=\left(\begin{smallmatrix}0 & 0 & 1 & 0\\
0 & 0 & 0 & -1\\
-1 & 0 & 0 & 0\\
0 & 1 & 0 & 0
\end{smallmatrix}\right).
\]
Comparing this with the definition $Sp(4,\mbZ):=\left\{ x\mid\,^{t}x\left(\begin{smallmatrix}0 & E_{2}\\
-E_{2} & 0
\end{smallmatrix}\right)x=\left(\begin{smallmatrix}0 & E_{2}\\
-E_{2} & 0
\end{smallmatrix}\right)\right\} $, we describe the (anti-)isomorphism $Sp(4,\mbZ)\simeq\tilde{\Gamma}$
by $x\mapsto I\,^{t}xI$ with $I=\mathrm{diag}(1,1,1,-1)$. Then we
can write the isomorphism $\Psi:Sp(4,\mbZ)/\left\{ \pm I_{4}\right\} \simeq SO(U^{\oplus2}\oplus\langle-2\rangle)_{+}$
explicitly by 
\[
\Psi(g)=\rho_{\wedge}(IgI),\,\,\,\,(g\in Sp(4,\mbZ)/\left\{ \pm I_{4}\right\} ).
\]
Using this, we verify the commutative diagram (\ref{eq:H2-comm-diag}). 

\section{\textbf{Table of $\Gamma_{0}(n)_{+}$}}

\subsection{Table of $\Gamma_{0}(n)_{+}$\label{subsec:AppB-table-Gamma}}

In the following Table 1, we list all genus zero groups $\Gamma_{0}(n)_{+}$
which we read from a larger list in \cite{CN}. The ``type'' in
the second column means the name used in the listing there. In the
third column, the numbers $N_{c}(n)$ of cusps of $\mathbb{H}_{+}/\Gamma_{0}(n)_{+}$
are listed. The case $n=1$ is included by setting $\Gamma_{0}(1)_{+}\equiv PSL(2,\mbZ)$
as a convention. We see in \cite[Table 3]{CN} that, if the number
cusps $N_{c}(n)$ is greater than 1, the Thompson series $T_{n}(t)$
of $\Gamma_{0}(n)_{+}$ has an expression by an eta product except
for $n=25,27,49,50$ and $54$ (which we indicate by {*} in the table).

\noindent%
\begin{tabular}{|c|c|c|}
\hline 
$n$ & type & \textbf{$c$}\tabularnewline
\hline 
\hline 
1 & 1A & 1\tabularnewline
\hline 
2 & 2A & 1\tabularnewline
\hline 
3 & 3A & 1\tabularnewline
\hline 
4 & 4A & 2\tabularnewline
\hline 
5 & 5A & 1\tabularnewline
\hline 
6 & 6A & 1\tabularnewline
\hline 
7 & 7A & 1\tabularnewline
\hline 
8 & 8A & 2\tabularnewline
\hline 
9 & 9A & 2\tabularnewline
\hline 
10 & 10A & 1\tabularnewline
\hline 
11 & 11A & 1\tabularnewline
\hline 
12 & 12A & 2\tabularnewline
\hline 
13 & 13A & 1\tabularnewline
\hline 
\end{tabular}~%
\begin{tabular}{|c|c|c|}
\hline 
$n$ & type & \textbf{$c$}\tabularnewline
\hline 
\hline 
14 & 14A & 1\tabularnewline
\hline 
15 & 15A & 1\tabularnewline
\hline 
16 & 16C & 3\tabularnewline
\hline 
17 & 17A & 1\tabularnewline
\hline 
18 & 18B & 2\tabularnewline
\hline 
19 & 19A & 1\tabularnewline
\hline 
20 & 20A & 2\tabularnewline
\hline 
21 & 21A & 1\tabularnewline
\hline 
22 & 22A & 1\tabularnewline
\hline 
23 & 23AB & 1\tabularnewline
\hline 
24 & 24B & 2\tabularnewline
\hline 
25 & 25A & 3{*}\tabularnewline
\hline 
26 & 26A & 1\tabularnewline
\hline 
\end{tabular}~%
\begin{tabular}{|c|c|c|}
\hline 
$n$ & type & $c$\tabularnewline
\hline 
\hline 
27 & 27A & 3{*}\tabularnewline
\hline 
28 & 28B & 2\tabularnewline
\hline 
29 & 29A & 1\tabularnewline
\hline 
30 & 30B & 1\tabularnewline
\hline 
31 & 31AB & 1\tabularnewline
\hline 
32 & 32A & 4\tabularnewline
\hline 
33 & 33B & 1\tabularnewline
\hline 
34 & 34A & 1\tabularnewline
\hline 
35 & 35A & 1\tabularnewline
\hline 
36 & 36A & 4\tabularnewline
\hline 
38 & 38A & 1\tabularnewline
\hline 
39 & 39A & 1\tabularnewline
\hline 
41 & 41A & 1\tabularnewline
\hline 
\end{tabular}~%
\begin{tabular}{|c|c|c|}
\hline 
$n$ & type & \textbf{$c$}\tabularnewline
\hline 
\hline 
42 & 42A & 1\tabularnewline
\hline 
44 & 44AB & 2\tabularnewline
\hline 
45 & 45A & 2\tabularnewline
\hline 
46 & 46CD & 1\tabularnewline
\hline 
47 & 47AB & 1\tabularnewline
\hline 
49 & 49Z & 4{*}\tabularnewline
\hline 
50 & 50A & 3{*}\tabularnewline
\hline 
51 & 51A & 1\tabularnewline
\hline 
54 & 54A & 3{*}\tabularnewline
\hline 
55 & 55A & 1\tabularnewline
\hline 
56 & 56A & 2\tabularnewline
\hline 
59 & 59AB & 1\tabularnewline
\hline 
60 & 60B & 2\tabularnewline
\hline 
\end{tabular}~%
\begin{tabular}{|c|c|c|}
\hline 
$n$ & type & \textbf{$c$}\tabularnewline
\hline 
\hline 
62 & 62AB & 1\tabularnewline
\hline 
66 & 66A & 1\tabularnewline
\hline 
69 & 69AB & 1\tabularnewline
\hline 
70 & 70A & 1\tabularnewline
\hline 
71 & 71AB & 1\tabularnewline
\hline 
78 & 78A & 1\tabularnewline
\hline 
87 & 87AB & 1\tabularnewline
\hline 
92 & 92AB & 2\tabularnewline
\hline 
94 & 94AB & 1\tabularnewline
\hline 
95 & 95AB & 1\tabularnewline
\hline 
105 & 105A & 1\tabularnewline
\hline 
110 & 110A & 1\tabularnewline
\hline 
119 & 119AB & 1\tabularnewline
\hline 
\end{tabular}

\begin{center}
\textbf{Table 1. }
\par\end{center}

\vskip0.2cm

\noindent 

\subsection{\label{subsec:AppendixB-Tm-product}Thompson series $T_{n}(t)+c_{n}$
by eta products }

When $N_{c}(n)\geq2$ and $n\not=25,27,49,50,54$, the Thompson series
with some additive constant $c_{n}$ is expressed by an eta product.
In Table 2, we list the product forms following the notation in \cite{CN},
i.e., $a_{1}^{e_{1}}a_{2}^{e_{2}}\cdots:=\eta_{a_{1}}^{e_{1}}\eta_{a_{2}}^{e_{2}}\cdots$.
For each $n$, we also list a value $b\in\mathbb{Q}$ for which $(T_{n}(t)+c_{n})^{b}\eta_{_{\mathrm{BCOV}}}(t)$
defines a cusp form (cf. Lemma \ref{lem:eta-cusp}) with the expression
of the cusp form. In the table, the blank sign indicates that we couldn't
find a value $b\in\mathbb{Q}$ with the required property.
\begin{center}
\begin{tabular}{|c|c|c|c|}
\hline 
$n$ & $T_{n}+c_{n}$ & cusp form & $b$\tabularnewline
\hline 
\hline 
4 & $2^{48}$/$1^{24}4^{24}$ & $2^{4}$ & $1/6$\tabularnewline
\hline 
8 & $2^{8}4^{8}/1^{8}8^{8}$ & $1\,2\,4\,8$ & $1/4$\tabularnewline
\hline 
9 & $3^{12}/1^{6}9^{6}$ & $1\,3^{2}9$ & $1/2$\tabularnewline
\hline 
12 & $2^{12}6^{12}/1^{6}3^{6}4^{6}12^{6}$ & $1\,3\,4\,12$ & $1/6$\tabularnewline
\hline 
16 & $\begin{alignedat}{1}2^{6}8^{6}/1^{4}4^{4}16^{4}\\
4^{10}/1^{2}2^{3}8^{3}16^{2}
\end{alignedat}
$ & $\begin{matrix}4^{8}/2^{2}8^{2}\\
2^{6}8^{6}/1^{\frac{8}{3}}4^{\frac{8}{3}}16^{\frac{8}{3}}
\end{matrix}$ & $\begin{matrix}-1\\
-2/3
\end{matrix}$\tabularnewline
\hline 
18 & $3^{4}6^{4}/1^{2}2^{2}9^{2}18^{2}$ & $3^{2}6^{2}$ & $1$\tabularnewline
\hline 
20 & $2^{8}10^{8}/1^{4}4^{4}5^{4}20^{4}$ & $2^{2}10^{2}$ & $1/2$\tabularnewline
\hline 
24 & $2^{2}4^{4}6^{2}12^{2}/1^{2}3^{3}8^{2}24^{2}$ & $2\,4\,6\,12$ & $1$\tabularnewline
\hline 
28 & $2^{6}14^{6}/1^{3}4^{3}7^{3}28^{3}$ & $1\,4\,7\,28$ & $1/3$\tabularnewline
\hline 
32 & $2^{3}16^{3}/1^{2}4\,8\,32^{2}$  & $4^{3}8^{3}/2\,16$ & $-1$\tabularnewline
\hline 
36 & $\begin{aligned} & 1\,4\,6^{16}9\,36/2^{4}3^{6}12^{6}18^{4}\\
 & 2^{6}3^{2}12^{2}18^{6}/1^{3}4^{3}6^{4}9^{3}36^{3}\\
 & 2^{2}3^{4}12^{4}18^{2}/1^{2}4^{2}6^{4}9^{2}36^{2}
\end{aligned}
$ & $\begin{matrix}-\\
-\\
6^{12}/3^{4}12^{4}
\end{matrix}$ & $\begin{aligned}-\\
-\\
-2
\end{aligned}
$\tabularnewline
\hline 
44 & $2^{4}22^{4}/1^{2}4^{2}11^{2}44^{2}$ & $2^{2}22^{2}$ & 1\tabularnewline
\hline 
45 & $3^{2}15^{2}/1\,5\,9\,45$ & $1\,5\,9\,45$ & 1\tabularnewline
\hline 
56 & $2\,4\,14\,28/1\,7\,8\,56$ & $1\,7\:8\,56$ & $1$\tabularnewline
\hline 
60 & $2^{2}6^{2}10^{2}30^{2}/1\,3\,4\,5\,12\,15\,20\,60$ & $2\,6\,10\,30$ & 1\tabularnewline
\hline 
92 & $2^{2}46^{2}/1\,4\,23\,92$ & $1\,4\,23\,92$ & 1\tabularnewline
\hline 
\end{tabular}
\par\end{center}

\begin{center}
\textbf{Table 2. }
\par\end{center}

\section{\textbf{\label{sec:Appendix-CD-family}Clingher-Doran's family}}

\subsection{Picard-Fuchs equations}

The Period integral (\ref{eq:w0-CD}) is a solution of the Picard-Fuchs
equations $\cD_{i}w_{0}(x)=0$ $(i=1,2,...,5)$ with
\begin{equation}
\begin{alignedat}{3} &  & \cD_{1} &  &  & =\theta_{x}(\theta_{x}-3\theta_{y}+\theta_{z})-x(2\theta_{x}-\theta_{y}+1)(2\theta_{x}-\theta_{y}),\\
 &  & \cD_{2} &  &  & =(\theta_{y}-\theta_{z})(\theta_{x}-3\theta_{y}+\theta_{z})+2z(2\theta_{z}+1)(3\theta_{y}-\theta_{z}),\\
 &  & \cD_{3} &  &  & =(2\theta_{x}-\theta_{y})(3\theta_{y}-\theta_{z})+2yz(\theta_{x}-2\theta_{y}-1)(\theta_{x}-3\theta_{y}+\theta_{z})\\
 &  &  &  &  & \hsp{103}+2yz^{2}(2\theta_{z}+1)(3\theta_{x}-6\theta_{y}+2\theta_{z}),\\
 &  & \cD_{4} &  &  & =(2\theta_{x}-\theta_{y})(\theta_{y}-\theta_{z})-4yz^{2}(2\theta_{y}+3)(\theta_{x}-3\theta_{y}+\theta_{z})\\
 &  &  &  &  & \hsp{98}-8yz^{3}(2\theta_{z}+1)(6\theta_{y}+4\theta_{z}+15),\\
 &  & \cD_{5} &  &  & =\theta_{x}(3\theta_{y}-\theta_{z})+2xyz(2\theta_{x}-\theta_{y})(\theta_{x}-3\theta_{y}+\theta_{z})\\
 &  &  &  &  & \hsp{68}+12xyz^{2}(2\theta_{x}-\theta_{y})(2\theta_{z}+1).
\end{alignedat}
\label{eq:exGKZ}
\end{equation}
These differential operators should follow from the extended GKZ system
as done in \cite{HKTY1,HLY}. Here we have determined these as a system
of differential operators which annihilate the hypergeometric series
$w_{0}(x)$ in (\ref{eq:CD-hypergeom}). The completeness of these
to determine all other period integrals can be verified by calculating
the indicial ideal at $x=y=z=0$ (see \cite{HLY}). 

\subsection{Griffiths-Yukawa couplings\label{subsec:AppD-GY-couplings}}

Let $\check{\mathfrak{X}}\rightarrow\mathcal{M}$ be Clingher-Doran's
family, and $\Omega:=\Omega(x)$ be a holomorphic two form of the
fiber K3 surface $X_{x}$ over $x\in\mathcal{M}$. Using the notation
$(A,B):=\int_{X_{x}}A\wedge B$ for $A,B\in H^{2}(X_{x},\mathbb{C})$,
we define the Griffiths-Yukawa couplings by $C_{ij}:=(\Omega,\partial_{i}\partial_{j}\Omega)$
where we set $(\partial_{1},\partial_{2},\partial_{3}):=(\frac{\partial\;}{\partial x},\frac{\partial\;}{\partial y},\frac{\partial\;}{\partial z})$.
By the Griffiths transversality, we have $(\Omega,\Omega)=(\Omega,\partial_{i}\Omega)=0$.
Using these relations, we have
\[
(\Omega,\partial_{i}\partial_{j}\partial_{k}\Omega)=\frac{1}{2}\partial_{i}(\Omega,\partial_{j}\partial_{k}\Omega)+\frac{1}{2}\partial_{j}(\Omega,\partial_{i}\partial_{k}\Omega)+\frac{1}{2}\partial_{k}(\Omega,\partial_{i}\partial_{j}\Omega).
\]
Then, by evaluating $(\Omega,\partial_{i}\cD_{a}\Omega)$ for $i=1,2,3,a=1,...,5$,
we obtain a system of first order differential equations for $C_{ij}$.
By integrating these equations, we can determine $C_{ij}$ up to a
multiplicative constant (see \cite{Candelas,HKTY1}). For example,
we obtain $C_{11}=\tilde{C}_{11}/dis_{0}$ with \textbf{\footnotesize{}
\[
\begin{alignedat}{2} & \tilde{C}_{11} &  & =2y^{2}z^{2}\big(115200xy^{3}z^{6}-82944xy^{2}z^{6}+9600xy^{3}z^{5}-27648xy^{2}z^{5}-976xy^{2}z^{4}\\
 &  &  & +88xy^{2}z^{3}-9xy^{2}z^{2}-103680y^{3}z^{6}+103680y^{2}z^{6}+25920y^{2}z^{5}-3456y^{2}z^{4}\\
 &  &  & -48y^{2}z^{3}+3456yz^{4}+912yz^{3}+12yz^{2}+4yz-4z-1\big)
\end{alignedat}
\]
} and {\footnotesize{}
\[
\begin{alignedat}{2} & dis_{0} & = & 1-4x+4z-16xz-4yz+18xyz+8xyz^{2}-27x^{2}y^{2}z^{2}-864yz^{3}+3456xyz^{3}\\
 &  &  & -48xy^{2}z^{3}+180x^{2}y^{2}z^{3}-3456yz^{4}+13824xyz^{4}+3456y^{2}z^{4}-15264xy^{2}z^{4}\\
 &  &  & -1280x^{2}y^{2}z^{4}-127872xy^{2}z^{5}+483840x^{2}y^{2}z^{5}+24000x^{2}y^{3}z^{5}\\
 &  &  & +186624y^{2}z^{6}-1990656xy^{2}z^{6}+4976640x^{2}y^{2}z^{6}+518400xy^{3}z^{6}\\
 &  &  & -2304000x^{2}y^{3}z^{6}+746496y^{2}z^{7}-5971968xy^{2}z^{7}+11943936x^{2}y^{2}z^{7}\\
 &  &  & -746496y^{3}z^{7}+6220800xy^{3}z^{7}-13824000x^{2}y^{3}z^{7}+3200000x^{3}y^{4}z^{7}
\end{alignedat}
\]
}By adjusting the multiplicative constant in $C_{ij}$, we can verify
the mirror symmetry relation (\ref{eq:YuakawaYttt}) with the Gram
matrix $(K_{ab})$ given by
\[
(K_{ab})=\left(\begin{smallmatrix}0 & 1 & 3\\
1 & 2 & 6\\
3 & 6 & 16
\end{smallmatrix}\right)=\,^{t}P\left(\begin{smallmatrix}0 & 1 & 0\\
1 & 0 & 0\\
0 & 0 & -2
\end{smallmatrix}\right)P,\,\,\,P:=\left(\begin{smallmatrix}1 & -1 & 0\\
0 & 1 & -3\\
0 & 0 & 1
\end{smallmatrix}\right).
\]

\subsection{\label{subsec:AppendixD-local-sol-GKZ}Local solutions around orbifold
points }

The extended GKZ system (\ref{eq:exGKZ}) is an irreducible system
of the GKZ system determined by the point configuration $\mathcal{A}_{3}=\left\{ v_{1},v_{2},...,v_{6}\right\} $
which we define by writing (\ref{eq:Delta-CD}) $\Delta=\mathrm{Conv.}\left\{ v_{1},\cdots,v_{6}\right\} $.
Precisely, adding the origin $v_{0}$ and putting the points on an
affine plane by $\bar{v}_{i}:=\,^{t}(1,v_{i})$, the point configuration
is given by 
\[
\mathcal{A}=\left\{ \left(\begin{smallmatrix}1\\
0\\
0\\
0
\end{smallmatrix}\right),\left(\begin{smallmatrix}1\\
-1\\
1\\
0
\end{smallmatrix}\right),\left(\begin{smallmatrix}1\\
2\\
-1\\
0
\end{smallmatrix}\right),\left(\begin{smallmatrix}1\\
0\\
-1\\
1
\end{smallmatrix}\right),\left(\begin{smallmatrix}1\\
-1\\
-1\\
1
\end{smallmatrix}\right),\left(\begin{smallmatrix}1\\
-1\\
-1\\
-1
\end{smallmatrix}\right),\left(\begin{smallmatrix}1\\
-1\\
-1\\
0
\end{smallmatrix}\right)\right\} ,
\]
and the GKZ system is the $\mathcal{A}$-hypergeometric system with
exponent $\beta=\,^{t}\left(-1,0,0,0\right)$. As an irreducible system
of GKZ system, the extended GKZ system (\ref{eq:exGKZ}) shares the
compactification of the parameter space with the GKZ system, which
is given by $\mathbb{P}_{SecP}$ with the secondary polytope $SecP$.
The secondary polytope is defined as the convex hull of vertices $v_{T}$
parametrized by regular triangulations $T$ of $\mathcal{A}$. In
the present case, we can calculate all regular triangulations (e.g.
by using TOPCOM \cite{TOPCOM}). It turns out that there are 16 regular
triangulations $T_{k}\,(k=1,...,16).$ We refrain listing all of them,
but we list those which correspond to orbifold points here; 
\[
\begin{aligned} &  & T_{9}= & \left\{ \{1,2,3,4\},\{1,2,4,5\}\right\} ,\\
 &  & T_{12}= & \left\{ \{1,3,4,5\},\{1,2,3,5\}\right\} , &  & T_{16}= & \left\{ \{1,3,5,6\},\{1,3,4,6\},\{1,2,3,5\}\right\} ,
\end{aligned}
\]
where numbers $i$ represent the point $\bar{v}_{i}$ of $\mathcal{A}$,
and $\left\{ 1,2,3,4\right\} $, for example, represents a simplex
of the triangulation of $\text{Conv.}(\mathcal{A})$. Since it is
straightforward to obtain the following results, although we need
some calculations, we only list local solutions of (\ref{eq:exGKZ})
which we obtain for the above three triangulations. 

\noindent$\bullet\,T_{9}$. Following \cite{GKZ2}, we can determine
the affine chart $A_{T_{9}}$ corresponding to the triangulation $T_{9}$
as 
\[
A_{T_{9}}=\mathrm{Spec}\mathbb{C}[\sigma_{T_{9}}^{\vee}\cap\mathbb{Z}^{3}]
\]
where we define a cone $\sigma_{T_{9}}^{\vee}=\mathbb{R}_{\geq0}\,^{t}(1,2,0)+\mathbb{R}_{\geq0}\,^{t}(-1,-2,-6)+\mathbb{R}_{\geq0}\,^{t}(-1,0,0)$
and make the identification $x^{m_{1}}y^{m_{2}}z^{m_{3}}$ with $^{t}(m_{1},m_{2},m_{3})\in\mathbb{Z}^{3}$.
This chart is singular at the origin, and we have to resolve the singularity.
However, as explained in the text, the orbifold regularity condition
does not depend on the resolution, and we can impose the regularity
condition at any affine chart of the resolution. Making a resolution,
we consider the following (smooth) affine coordinates which are related
to the affine coordinates $x,y,z$ for the LCSL: 
\[
(\tilde{x},\tilde{y},\tilde{z})=(\frac{1}{xy^{2}z^{2}},yz,xy).
\]
Using these coordinates, the local solutions $\tilde{w}_{k}(\tilde{x})$
are given by
\[
\begin{aligned} &  & \tilde{w}_{0}= & \tilde{x}^{\frac{5}{12}}\tilde{y}^{\frac{1}{3}}\tilde{z}^{\frac{1}{3}}(1+\cdots), & \tilde{w}_{1}= & \tilde{x}^{\frac{1}{2}}\tilde{y}^{\frac{1}{2}}\tilde{z}^{\frac{1}{2}}(1+\cdots), & \tilde{w}_{2}= & \tilde{x}^{\frac{7}{12}}\tilde{y}^{\frac{2}{3}}\tilde{z}^{\frac{2}{3}}(1+\cdots)\\
 &  & \tilde{w}_{3}= & \tilde{x}^{\frac{11}{12}}\tilde{y}^{\frac{4}{3}}\tilde{z}^{\frac{1}{3}}(1+\cdots), & \tilde{w}_{4}= & \tilde{x}^{\frac{13}{12}}\tilde{y}^{\frac{5}{3}}\tilde{z}^{\frac{2}{3}}(1+\cdots),
\end{aligned}
\]
where $\tilde{w}_{0}$ is the unique solution with the property that
$\tilde{w}_{k}/\tilde{w}_{0}$ are regular. In the above equations,
$\cdots$ represents higher order terms $O(\tilde{x})$ with respect
to $\tilde{x},\tilde{y},\tilde{z}$. 

\noindent$\bullet\,T_{12}$. In a similar way as above, we have 
\[
A_{T_{12}}=\mathrm{Spec}\mathbb{C}[\sigma_{T_{12}}^{\vee}\cap\mathbb{Z}^{3}]
\]
with $\sigma_{T_{12}}^{\vee}=\mathbb{R}_{\geq0}\,^{t}(-1,-2,0)+\mathbb{R}_{\geq0}\,^{t}(-1,-2,-5)+\mathbb{R}_{\geq0}\,^{t}(-1,0,0)$.
Making a resolution of the singularity of $A_{T_{12}}$ at the origin,
we take the following (smooth) affine coordinates: 
\[
(\tilde{x},\tilde{y},\tilde{z})=(\frac{1}{yz^{2}},\frac{1}{xz},z).
\]
Using these coordinates, the local solutions $\tilde{w}_{k}(\tilde{x})$
are given by
\[
\begin{aligned} &  & \tilde{w}_{0}= & \tilde{x}^{\frac{1}{5}}\tilde{y}^{\frac{1}{10}}(1+\cdots), & \tilde{w}_{1}= & \tilde{x}^{\frac{2}{5}}\tilde{y}^{\frac{7}{10}}(x+c_{1}z+\cdots),\\
 &  & \tilde{w}_{2}= & \tilde{x}^{\frac{3}{5}}\tilde{y}^{\frac{3}{10}}(1+\cdots), & \tilde{w}_{3}= & \tilde{x}^{\frac{4}{5}}\tilde{y}^{\frac{9}{10}}(x+c_{2}z+\cdots),\,\,\,\,\,\,\,\,\tilde{w}_{4}=\tilde{x}\tilde{y}(1+\cdots),
\end{aligned}
\]
where $c_{1}=-\frac{38400}{11},c_{2}=-\frac{153600}{119}$ and $\tilde{w}_{0}$
is the unique solution with the property that $\tilde{w}_{k}/\tilde{w}_{0}$
are regular. As above, $\cdots$ represents higher order terms.~

\noindent$\bullet\,T_{16}$. For this, we can calculate 
\[
A_{T_{16}}=\mathrm{Spec}\mathbb{C}[\sigma_{T_{16}}^{\vee}\cap\mathbb{Z}^{3}]
\]
with $\sigma_{T_{16}}^{\vee}=\mathbb{R}_{\geq0}\,^{t}(1,0,0)+\mathbb{R}_{\geq0}\,^{t}(-1,-2,-5)+\mathbb{R}_{\geq0}\,^{t}(-3,-2,0)$.
As affine coordinates of a resolution of $A_{T_{16}}$, we take 
\[
(\tilde{x},\tilde{y},\tilde{z})=(\frac{1}{x^{3}y},x,\frac{x}{z}).
\]
Using these coordinates, the local solutions $\tilde{w}_{k}(\tilde{x})$
are given by

\[
\begin{aligned} &  & \tilde{w}_{0}= & \tilde{x}^{\frac{1}{5}}\tilde{z}^{\frac{1}{2}}(1+\cdots), & \tilde{w}_{1}= & \tilde{x}^{\frac{2}{5}}\tilde{z}^{\frac{1}{2}}(1+\cdots), & \tilde{w}_{2}= & \tilde{x}^{\frac{3}{5}}\tilde{z}^{\frac{1}{2}}(1+\cdots)\\
 &  & \tilde{w}_{3}= & \tilde{x}^{\frac{4}{5}}\tilde{z}^{\frac{1}{2}}(1+\cdots), & \tilde{w}_{4}= & \tilde{x}\tilde{z}^{\frac{1}{2}}(1+\cdots),
\end{aligned}
\]
where $\tilde{w}_{0}$ is the unique solution with the property that
$\tilde{w}_{k}/\tilde{w}_{0}$ are regular. As above, $\cdots$ represents
higher order terms.
\begin{rem}
Although we do not describe details, all the regular solutions above
are given by hypergeometric series due to GKZ, which have the form
\[
\Phi_{\gamma}(\tilde{x},\tilde{y},\tilde{z})=\sum_{\ell\in L}\frac{a_{0}}{\prod_{i}\Gamma(1+\ell_{i}+\gamma_{i})}\prod_{i=0}^{6}a_{i}^{\ell_{i}+\gamma_{i}},
\]
where $L=\text{Ker}\left\{ A:\mathbb{Z}^{4}\rightarrow\mathbb{Z}^{4}\right\} $
with $A$ is the integral $4\times7$ matrix naturally defined by
$\mathcal{A}$. The exponent $\gamma$ is a solution of $A\gamma=\beta$,
and the affine coordinates $\tilde{x},\tilde{y},\tilde{z}$ are defined
from a $\mathbb{Z}$-basis of $L$ which is compatible with the corresponding
triangulation $T$. See {[}GKZ{]} for mode detail. 
\end{rem}

\section{\textbf{\label{sec:AppendixD-Elliptic}$\tau_{_{\mathrm{BCOV}}}$
for elliptic curves }}

\subsection{Period integrals and Picard-Fuchs equation}

Clingher-Doran's family (\ref{eq:def-eq-CD}) is a generalization
of the Weierstrass normal form of a family of elliptic curves; $\ty^{2}=4\tx^{3}-g_{2}\tx-g_{3}$
where $g_{2}=12E_{4}(\tau),g_{3}=216E_{6}(\tau)$ with the standard
Eisenstein series $E_{4}$ and $E_{6}$. It is a classical fact that
the weighted projective space $\mathbb{P}^{1}(2,3)$ of points $[g_{2},g_{3}]$
gives a coarse moduli space of elliptic curves. Period integrals of
this family of elliptic curve can be described in a quite parallel
way to Subsect. \ref{subsec:GKZ-for-CD-family}, and based on it,
we can find a family $\mathfrak{X}\rightarrow\mathbb{P}_{\Sigma}\simeq\mathbb{P}^{1}$,
i.e., a parameter space $\mathbb{P}_{\Sigma}$ over which we have
a universal family. Here we sketch calculations for this. 

We embed the Weierstrass normal form into more general cubics 
\[
F=a_{0}\tx\ty\tz+a_{1}\ty^{2}\tz+a_{2}\tx^{3}+a_{3}\tz^{3}+a_{4}\tx\tz^{2}+a_{5}\tx^{2}\tz+a_{6}\ty\tz^{2}
\]
in $\mathbb{P}^{2}$ with $(a_{0},a_{1},...,a_{6})=(0,1,-4,g_{3},g_{2},0,0)$.
We associate to the Newton polytope of $F$ a reflexive polytope $\Delta_{F}=\text{Conv.}\left\{ \left(\begin{smallmatrix}1\\
-1
\end{smallmatrix}\right),\left(\begin{smallmatrix}-1\\
2
\end{smallmatrix}\right),\left(\begin{smallmatrix}-1\\
-1
\end{smallmatrix}\right)\right\} $ in $\mathbb{R}^{2}$ where $\tx\ty\tz$ corresponds to the origin
and $\tx\tz^{2},\tx^{2}\tz,\ty\tz^{2}$ are on lines of codimension
one. These three monomials indicates that a (reductive) group $G$
is acting on $F$ as its automorphism. It is straightforward to find
that $G$ is generated by 
\[
(\tx,\ty,\tz)\mapsto(\tx,\ty+c_{1}\tx,\tz);\,\,(\tx,\ty,\tz)\mapsto(\tx,\ty+c_{2}\tz,\tz);\,\,(\tx,\ty,\tz)\mapsto(\tx+c_{3}\tz,\ty,\tz).
\]
As in the text, a preferred form $F_{toric}(a)$ for the period integral
(cf. (\ref{eq:w0-CD})) is given by $a_{4}=a_{5}=a_{6}=0$, since
in this form it is easy to find the Picard-Fuchs equation 
\begin{equation}
\left\{ \theta_{x}^{2}-12x(6\theta_{x}+5)(6\theta_{x}+1)\right\} w(x)=0\label{eq:Appendix-PF-E}
\end{equation}
with $x=\frac{a_{1}^{3}a_{2}^{2}a_{3}}{a_{0}^{6}}$. Namely the family
$\left\{ F_{toric}(a)\right\} _{a\in(\mathbb{C}^{*})^{4}}$ defines
the desired family $\mathfrak{X}\rightarrow\mathbb{P}^{1}$ by the
coordinate $x$. We write the singularities of this equation in Riemann's
$\mathcal{P}$-scheme,\def\RPe{\left\{ 
\small\begin{array}{cccccc} 
0 & \frac{1}{432} & \infty \\ 
\hline  
0 & 0 & \frac{1}{6}\\ 
0 & 0 & \frac{5}{6}\\
\end{array}  \right\}
}
\begin{equation}
\RPe\label{eq:Appendix-RP-elliptic}
\end{equation}
Using the actions of $G$, we can find that the Weierstrass normal
form and $F_{toric}(a)$ are related by
\[
(a_{0},a_{1},a_{2},a_{3})=(-2i\sqrt{2}(3\,g_{2})^{\frac{1}{4}},1,-4,-\frac{g_{2}^{\frac{3}{2}}}{3\sqrt{3}}+g_{3}).
\]
From this, we can find relations 
\begin{equation}
x=\frac{\sqrt{3}g_{2}^{\frac{3}{2}}-9g_{3}}{864\sqrt{3}g_{2}^{\frac{3}{2}}},\,\,\,\frac{1}{x(1-432x)}=\frac{1728g_{2}^{3}}{g_{2}^{3}-27g_{3}^{2}}=\frac{1}{j(\tau)}.\label{eq:AppE-x-mirror-map}
\end{equation}
As a power series solution, we have $w_{0}(x)=\sum_{n\geq0}\frac{(6n)!}{(3n)!(2n)!n!}x^{n}$,
and find a relation $w_{0}(x(\tau))=E_{4}(\tau)^{\frac{1}{4}}.$ 

\subsection{$\tau_{_{\mathrm{BCOV}}}$ for the family $\mathfrak{X}\rightarrow\mathbb{P}^{1}$~}

As shown in (\ref{eq:Appendix-RP-elliptic}), the family degenerates
isomorphically at $x=0$ and $432x-1=0$. In fact, it is easy to verify
that the Picard-Fuchs equation (\ref{eq:Appendix-PF-E}) is invariant
under changing variable by $x=\frac{1}{432}-u$. Respecting this symmetry
and also imposing the regularity at the orbifold point $x=\infty$,
we determine $\tau_{_{\mathrm{BCOV}}}$ of this family by 
\begin{equation}
\tau_{_{\mathrm{BCOV}}}=\left(\frac{1}{w_{0}(x)}\right)^{3}\frac{dx}{d\tau}(1-432x)^{-1-\frac{1}{12}}x^{-1-\frac{1}{12}}=\frac{1}{\eta(\tau)^{2}}.\label{eq:Appendix-tau-E}
\end{equation}
Here the power of $\left(\frac{1}{w_{0}(x)}\right)^{3}$ is determined
from the BCOV formula (\ref{eq:BCOV-CY3}) by setting $h_{X}^{1,1}=0$
and $\chi=0$. 

\subsection{Monodromy group $\mathcal{G}_{E}$}

Starting from $\left\{ F_{toric}(a)\right\} _{a\in(\mathbb{C}^{*})^{4}}$,
we obtain a universal family $\mathfrak{X}$ over $\mathbb{P}_{\Delta}=\mathbb{P}^{1}$,
while from the Weierstrass normal form, we only come to a coarse moduli
space $\mathbb{P}^{1}(2,3)\simeq\mathbb{P}^{1}$ which is isomorphic
to $\overline{\mathbb{H}_{+}/PSL(2,\mathbb{Z})}$. We can see the
difference in the monodromy group $\mathcal{G}_{E}$ of the family
$\mathfrak{X}\rightarrow\mathbb{P}^{1}$. Here we briefly summarize
calculations of $\mathcal{G}_{E}$.

As an integral basis of the Picard-Fuchs equation, we take 
\[
I_{1}(x)=\frac{1}{2\pi i}w_{0}(x)\,\log x+w_{1}(x),\,\,I_{0}(x)=w_{0}(x),
\]
where $w_{0}(x)$ is the power series solution around $x=0$ above,
and $w_{1}(x)$ is a power series of the form $w_{1}(x)=cx+O(x^{2})$.
Arranging these solutions to a column vector $\,^{t}(I_{1}(x)\,I_{0}(x))$,
it is straightforward to obtain the monodromy matrices around $x=0,\frac{1}{432}$
and $\infty$, respectively, as 
\[
M_{0}=\left(\begin{matrix}1 & 1\\
0 & 1
\end{matrix}\right),\,M_{1}=\left(\begin{matrix}1 & 0\\
-1 & 1
\end{matrix}\right),\,M_{\infty}=\left(\begin{matrix}1 & -1\\
1 & 0
\end{matrix}\right).
\]
From these, we obtain $P\mathcal{G}_{E}=\langle TS,T\rangle$ where
$PSL(2,\mathbb{Z})=\langle S,T\rangle$ with $S=\left(\begin{smallmatrix}0 & -1\\
1 & 0
\end{smallmatrix}\right),T=\left(\begin{smallmatrix}1 & 1\\
0 & 1
\end{smallmatrix}\right)$, and we see that $P\mathcal{G}_{E}=\langle S,PSL(2,\mathbb{Z})\rangle$.
From the last relation, we see that the fundamental domain of $P\mathcal{G}_{E}$
is given by 
\[
\mathcal{F}_{E}=\mathcal{F}\cup(S\cdot\mathcal{F}),
\]
where $\mathcal{F}$ represents the standard fundamental domain of
$PSL(2,\mathbb{Z})$. 

\subsection{Connection matrix $C_{[0,1]}$ }

It is clear that $S$ acts on $\mathcal{F}_{E}$ involutively, exchanging
the two cusp points. Period map $\mathcal{P}:\mathbb{P}^{1}\rightarrow\mathbb{H}_{+}$
defined by $\tau(x)=\frac{I_{1}(x)}{I_{0}(x)}$ is multi-valued. We
can fix the branch of this map to have an isomorphism $\mathbb{P}^{1}\simeq\mathcal{F}_{E}$,
whose inverse map gives the mirror map $x=x(\tau)$. The $S$-transform
$S\s\tau$ of $\tau\in\mathcal{F}_{E}$ is contained in $\mathcal{F}_{E}$.
Correspondingly it holds that $x(S\s\tau)=\frac{1}{432}-x(\tau)$
under the isomorphism. We can derive the last relation as follows:
Since the local solutions around $x=0$ and $x=\frac{1}{432}$ have
overlaps in their domains of convergence, it is easy to obtain the
following relations: 
\[
I_{0}(x)=-i\,I_{1}(\frac{1}{432}-x),\,\,\,\,\,I_{1}(x)=i\,I_{0}(\frac{1}{432}-x).
\]
If we set the integral basis around $u:=\frac{1}{432}-x=0$ by $^{t}(I_{1}(u),I_{0}(u))$,
we have $C_{[0,1]}=\left(\begin{smallmatrix}0 & -i\\
i & 0
\end{smallmatrix}\right)=i\,S$ for the connection matrix. Note that this result gives us 
\[
w_{0}(x(S\s\tau))=-i\,\frac{I_{1}(x(\tau))}{I_{0}(x(\tau))}w_{0}(x(\tau))=-i\tau\,w_{0}(x(\tau)),
\]
which tells us the branch of the relation $w_{0}(\tau)=(E_{4}(\tau))^{\frac{1}{4}}$,
i.e., $(E_{4}(S\s\tau))^{\frac{1}{4}}=(\tau^{4}E_{4}(\tau))^{\frac{1}{4}}=-i\tau(E_{4}(\tau))^{\frac{1}{4}}$.
This also determines $(g_{2}(S\s\tau))^{\frac{3}{2}}=-\tau^{6}(g_{2}(\tau))^{\frac{3}{2}}$,
which gives rise to the relation $x(S\s\tau)=\frac{1}{432}-x(\tau)$
by (\ref{eq:AppE-x-mirror-map}).

~

~

{\small{}Shinobu Hosono}{\small\par}

{\small{}Department of Mathematics, Gakushuin University }{\small\par}

{\small{}Mejiro, Toshima-ku, Tokyo 171-8588, Japan }{\small\par}

{\small{}e-mail: hosono@math.gakushuin.ac.jp}{\small\par}

\selectlanguage{english}%
{\small{}~}{\small\par}

\selectlanguage{american}%
{\small{}Atsushi Kanazawa}{\small\par}

{\small{}Faculty of Policy Management, Keio University }{\small\par}

{\small{}Endo 5322, Fujisawa, Kanagawa, 252-0882, Japan}{\small\par}

{\small{}e-mail: atsushik@sfc.keio.ac.jp}{\small\par}

\newpage

\setcounter{page}{1}
%\overfullrule=0pt
%\input epsf
%\documentclass[10pt]{amsart}
%\usepackage{amscd, amsmath, amsthm}
%\begin{document}

\centerline{
{\bf K3 differential operators for genus zero groups $\Gamma_{0}(n)_{+}$}}
\vskip0.3cm
\centerline{
{\bf --- addendum to ``BCOV cusp forms of lattice polarized K3 surfaces''
---}}
\vskip0.3cm
\centerline{ S. Hosono and A. Kanazawa}
\vskip1cm

In the paper \cite{[1]}, 
we have introduced {\it K3 differential operators} (Proposition 4.9) for 
all genus zero groups of type $\Gamma_0(n)_+$. About the genus zero groups, 
all details including historical backgrounds can be found 
in the reference \cite{[2]}. In particular, the $q$-series expansions of 
Thompson series (up to additive constants) are given in Table 4 of \cite{[2]}. 
To have an easy access to the Table 4 in [2], we specify a genus 
zero group of $\Gamma_0(n)_+$ by the corresponding {\it type},  
which represents a conjugary class of the Monster group 
following ATLAS table \cite{[3]}. For a quick view of the correspondence, 
we refer to a table in Appendix B of \cite{[1]}. 

For each type $m$, we have the Thompson series which is defined by  
\[
T_m(q) = \frac{1}{q}+0+H_1(m) q+H_2(m) q^2+\cdots
\]
with the so-called head representations $H_r$ of the Monstor group. 
For $m$, we have the corresponding genus zero group $\Gamma_0(n)_+$, 
and K3 differential operator $\mathcal{D}_m$ whose mirror map is given by 
the Thompson series as  
\[
\frac{1}{x(q)}=T_m(q)+c_0
\eqno{(1)}
\]
with some constant $c_0$. For several types $m$, this relation was 
first discovered in [4]; and it was conjectured in general for all 
genus zero groups. Proposition 4.9 in \cite{[1]} gives us a complete list 
of K3 differential operators for genus zero groups $\Gamma_0(n)_+$, 
which verifies partly the conjecture by Lian and Yau \cite{[4]}. 
It is not clear, however, whether every K3 differential operator (listed below) 
comes from a family of K3 surfaces or not. 

In the following table, the first column shows the type $m$, and the 
corresponding K3 differential operator $\mathcal{D}_m$ is in 
the second column. In the third column, the number $N_L$ of LCSL  
points of $\mathcal{D}_m$ and the constant term $c_0$ in 
(1) are shown in the form; 
$
\begin{matrix} N_{L} \\ c_0 \end{matrix}
$. 
We can verify that 
$N_{L}$ is equal to the number of the cusps $N_c$ for the corresponding 
group $\Gamma_0(n)_+$. The constant $c_0$ is not uniquely determined. If 
$N_L\geq 2$, $c_0$ is determined so that $T_n+c_n$ has an eta product except 
for $n=25,27,49,50,54$. 
If $N_L=1$, we tried to determine $c_0$ so that the unique power series 
$w_0(x)=1+\sum_n a_n x^n$ has integer coefficients $a_n$ 
at least for lower $n\leq 20$ (cf. Remark 4.9 of \cite{[1]}). Almost all cases, 
we were able to find such $c_0$, but we couldn't find such values of 
$c_0$ for some cases, which are indicated by $\dagger$ in $c_0^\dagger$.

\hfill
$\;$
\newpage

%\centerline{{\bf Table of K3 differential operators }}
%\vskip0.3cm
\begin{tabular}{|c|c|c|}
\hline 
$m$ & K3 differential operator $\mathcal{D}_m$ & 
$\begin{matrix} N_L \\ c_0 \end{matrix}$ 
\tabularnewline
\hline 
\hline 
${}{1A}$ 
&  
$\theta ^3-24 x (2 \theta +1) (6 \theta +1) (6 \theta +5)$
& 
$\begin{matrix} 1 \\ 744 \end{matrix}$ 
\tabularnewline
\hline 
${}{2A}$ 
&  
$
\theta ^3-8 x (2 \theta +1) (4 \theta +1) (4 \theta +3)
$
& $\begin{matrix} 1 \\ 104 \end{matrix}$ 
\tabularnewline
\hline 
${}{3A}$ 
&  
$
\theta ^3-6 x (2 \theta +1) (3 \theta +1) (3 \theta +2)
$
& $\begin{matrix} 1 \\ 42 \end{matrix}$ 
\tabularnewline
\hline 
${}{4A}$ 
&  
$
\theta ^3-8 x (2 \theta +1)^3
$
& $\begin{matrix} 2 \\ 24 \end{matrix}$ 
\tabularnewline
\hline 
${}{5A}$ 
&  
$\begin{matrix}
\theta^3
-2 x \left(23 \theta ^3+33 \theta ^2+17 \theta +3\right)
+4 x^2 \theta  \left(19 \theta ^2+21 \theta +6\right) \\ 
-24 x^3 \theta  (2 \theta +1) (3 \theta +2) 
-16 x^4 (\theta +1) (2 \theta +1) (2 \theta +3)
\end{matrix}
$
& $\begin{matrix} 1 \\ 16 \end{matrix}$ 
\tabularnewline
\hline 
${}{6A}$ 
&  
$\begin{matrix}
\theta ^3
-2 x \left(11 \theta ^3+30 \theta ^2+16 \theta +3\right)
-36 x^2 \theta  (2 \theta +1) (8 \theta +5) \\
+648 x^3 (\theta +1) (2 \theta +1) (2 \theta +3)
\end{matrix}$
& $\begin{matrix} 1 \\ 14 \end{matrix}$ 
\tabularnewline
\hline 
${}{7A}$ 
&  
$\begin{matrix}
\theta ^3
+x \left(-32 \theta ^3-39 \theta ^2-21 \theta -4\right)
+3 x^2 \theta  \left(43 \theta ^2+51 \theta +16\right) \\
+18 x^3 (\theta +1) (3 \theta +2) (3 \theta +4)
\end{matrix}$
& $\begin{matrix} 1 \\ 9 \end{matrix}$ 
\tabularnewline
\hline 
${}{8A}$ 
&  
$\begin{matrix}
%pf[8 A]=
\theta ^3
-4 x (\theta +1) \left(5 \theta ^2+4 \theta +1\right)
-16 x^2 \theta  \left(5 \theta ^2+6 \theta +2\right) \\
+64 x^3 (\theta +1)^3
\end{matrix}$
& $\begin{matrix} 2 \\ 8 \end{matrix}$ 
\tabularnewline
\hline 
${}{9A}$ 
&  
$\begin{matrix}
%pf[9 A]=
\theta ^3
-3 x (3 \theta +1) \left(3 \theta ^2+2 \theta +1\right)
+27 x^2 \theta  \left(5 \theta ^2+6 \theta +2\right) \\
+243 x^3 (\theta +1)^3
\end{matrix}$
& $\begin{matrix} 2 \\ 6 \end{matrix}$ 
\tabularnewline
\hline 
${}{10A}$ 
&  
$\begin{matrix}
%pf[10 A]=
\theta ^3
-2 x \left(21 \theta ^3+9 \theta ^2+5 \theta +1\right)
+4 x^2 \theta  \left(299 \theta ^2+87 \theta +28\right) \\
-120 x^3 \theta  \left(299 \theta ^2+87 \theta +28\right)
+3600 x^4 \theta  \left(74 \theta ^2+87 \theta +28\right) \\
+108000 x^5 (\theta +1) (4 \theta +3) (4 \theta +5)
\end{matrix}$
& $\begin{matrix} 1 \\ 4 \end{matrix}$ 
\tabularnewline
\hline 
${}{11A}$ 
&  
$\begin{matrix}
%pf[11 A]=
\theta ^3
-2 x \left(2 \theta ^3+15 \theta ^2+9 \theta +2\right)
-x^2 \theta  \left(41 \theta ^2+312 \theta +112\right) \\
-4 x^3 \theta  \left(902 \theta ^2+1050 \theta +371\right)  \\
+8 x^4 (\theta +1) \left(1473 \theta ^2+2814 \theta +1520\right) \\
-4906 x^5 (\theta +1) (\theta +2) (2 \theta +3)
\end{matrix}$
& $\begin{matrix} 1 \\ 6 \end{matrix}$ 
\tabularnewline
\hline 
${}{12A}$ 
&  
$\begin{matrix}
%pf[12 A]=
\theta ^3
-2 x \left(2 \theta ^3+15 \theta ^2+9 \theta +2\right) \\
-32 x^2 \theta  \left(8 \theta ^2+9 \theta +3\right) 
+1024 x^3 (\theta +1)^3
\end{matrix}$
& $\begin{matrix} 2 \\ 6 \end{matrix}$ 
\tabularnewline
\hline 
${}{13A}$ 
&  
$\begin{matrix}
%pf[13 A]=
18 \theta ^3
+12 x \left(55 \theta ^3-108 \theta ^2-72 \theta -18\right) \\
+x^2 \theta  \left(40337 \theta ^2-65556 \theta -25416\right)\\
-6 x^3 \theta  \left(712163 \theta ^2+645294 \theta +240298\right) \\
+x^4 \left(102725191 \theta ^3+250140666 \theta ^2+235007148 \theta +79888024\right) \\
-6 x^5 (6 \theta +7) \left(30044871 \theta ^2+79818815 \theta +58759441\right)
\\
+36 x^6 (6 \theta +7) (6 \theta +13) (4122537 \theta +7564391) \\
-46966608 x^7 (6 \theta +7) (6 \theta +13) (6 \theta +19)
\end{matrix}
$
& $\begin{matrix} 1 \\ 12 \end{matrix}$ 
\tabularnewline
\hline 
${}{14A}$ 
&  
$\begin{matrix}
%pf[14 A]=
5 \theta ^3
+x \left(31 \theta ^3-165 \theta ^2-105 \theta -25\right) \\
-x^2 \theta  \left(2497 \theta ^2+2838 \theta +1046\right) \\
+x^3 (\theta +1) \left(16081 \theta ^2+30692 \theta +16941\right) \\
-13818 x^4 (\theta +1) (\theta +2) (2 \theta +3)
\end{matrix}$
& $\begin{matrix} 1 \\ 6 \end{matrix}$ 
\tabularnewline
\hline 
${}{15A}$ 
&  
$\begin{matrix}
%pf[15 A]=
\theta ^3
+x \left(14 \theta ^3-39 \theta ^2-25 \theta -6\right) \\
-3 x^2 \theta  \left(277 \theta ^2+311 \theta +114\right) \\
+40 x^3 (\theta +1) \left(191 \theta ^2+355 \theta +186\right) \\
-900 x^4 (\theta +1) (\theta +2) (31 \theta +45)
+36000 x^5 (\theta +1) (\theta +2) (\theta +3)
\end{matrix}$
& $\begin{matrix} 1 \\ 7 \end{matrix}$ 
\tabularnewline
\hline 
\end{tabular}

\begin{tabular}{|c|c|c|}
\hline
${}{16C}$ 
&  
$\begin{matrix}
%pf[16 C]=
\theta ^3
-2 x (3 \theta +1) \left(3 \theta ^2+2 \theta +1\right) \\
+8 x^2 \theta  \left(11 \theta ^2+15 \theta +7\right)
-16 x^3 (\theta +1) \left(11 \theta ^2+7 \theta +3\right) \\
+16 x^4 (3 \theta +2) \left(3 \theta ^2+4 \theta +2\right)
-32 x^5 (\theta +1)^3
\end{matrix}
$
& $\begin{matrix} 3 \\ 4 \end{matrix}$ 
\tabularnewline
\hline 
${}{17A}$ 
&  
$\begin{matrix}
%pf[17 A]=
\theta ^3
+x \left(-41 \theta ^3-9 \theta ^2-5 \theta -1\right) 
+x^2 \theta  \left(183 \theta ^2+234 \theta +86\right) \\
+x^3 \left(917 \theta ^3+2709 \theta ^2+2921 \theta +1123\right) \\
+2 x^4 (2 \theta +3) \left(241 \theta ^2+717 \theta +575\right) \\
+140 x^5 (\theta +2) (2 \theta +3) (2 \theta +5)
\end{matrix}
$
& $\begin{matrix} 1 \\ 2 \end{matrix}$ 
\tabularnewline
\hline 
${}{18B}$ 
&  
$\begin{matrix}
%pf[18 B]=
16 \theta ^3
-8 x \left(47 \theta ^3+21 \theta ^2+13 \theta +3\right) 
+36 x^2 \theta  \left(42 \theta ^2+49 \theta +18\right) \\
+18 x^3 \left(332 \theta ^3+1032 \theta ^2+1141 \theta +444\right)
-891 x^4 (2 \theta +3)^3
\end{matrix}$
& $\begin{matrix} 2\; \\ 2^\dagger \end{matrix}$ 
\tabularnewline
\hline 
${}{19A}$ 
&  
$\begin{matrix}
%pf[19 A]=
882 \theta ^3
-294 x \left(5 \theta ^3+54 \theta ^2+36 \theta +9\right)  \\
+98 x^2 \theta  \left(151 \theta ^2-1080 \theta -441\right) 
-294 x^3 \theta  \left(3570 \theta ^2+4260 \theta +1771\right) \\
+98 x^4 \left(24019 \theta ^3+72444 \theta ^2+82866 \theta +34597\right) \\
+98 x^5 \left(6041 \theta ^3+23418 \theta ^2+34287 \theta +18954\right) \\
-314286 x^6 (\theta +2) (3 \theta +5) (3 \theta +7)
\end{matrix}
$
& $\begin{matrix} 1 \\ 3 \end{matrix}$ 
\tabularnewline
\hline 
${}{20A}$ 
&  
$\begin{matrix}
%pf[20 A]=
\theta ^3
+4 x \left(\theta ^3-6 \theta ^2-4 \theta -1\right)
-16 x^2 \theta  \left(16 \theta ^2+18 \theta +7\right) \\
+8 x^3 \left(152 \theta ^3+444 \theta ^2+466 \theta +173\right)
-160 x^4 (2 \theta +3)^3
\end{matrix}
$
& $\begin{matrix} 2 \\ 4 \end{matrix}$ 
\tabularnewline
\hline 
${}{21A}$ 
&  
$\begin{matrix}
%pf[21 A]=
11 \theta ^3
+x \left(347 \theta ^3-660 \theta ^2-462 \theta -121\right) \\
-x^2 \theta  \left(27527 \theta ^2+27684 \theta +11373\right) \\
+x^3 \left(516169 \theta ^3+1352262 \theta ^2+1356837 \theta +493276\right) \\
-6 x^4 (3 \theta +4) \left(237216 \theta ^2+656335 \theta +500699\right) \\
+360 x^5 (3 \theta +4) (3 \theta +7) (5067 \theta +9551) \\
-885600 x^6 (3 \theta +4) (3 \theta +7) (3 \theta +10)
\end{matrix}
$
& $\begin{matrix} 1 \\ 9 \end{matrix}$ 
\tabularnewline
\hline 
${}{22A}$ 
& 
$\begin{matrix}
%pf[22 A]=
\theta ^3
-2 x \left(5 \theta ^3+6 \theta ^2+4 \theta +1\right)
+2 x^2 \theta  \left(41 \theta ^2+6 \theta +2\right) \\
-2 x^3 \theta  \left(258 \theta ^2+327 \theta +143\right)
-4 x^4 \left(96 \theta ^3+339 \theta ^2+429 \theta +196\right) \\
+4 x^5 (2 \theta +3) \left(323 \theta ^2+981 \theta +810\right)
-560 x^6 (\theta +2) (2 \theta +3) (2 \theta +5)
\end{matrix}
$
& $\begin{matrix} 1 \\ 2 \end{matrix}$ 
\tabularnewline
\hline 
${}{23AB}$ 
&  
$\begin{matrix}
%pf[23 AB]=
\theta ^3
+x \left(8 \theta ^3-21 \theta ^2-15 \theta -4\right) 
-x^2 \theta  \left(251 \theta ^2+291 \theta +128\right) \\ 
+x^3 \left(1148 \theta ^3+3285 \theta ^2+3671 \theta +1492\right) \\
-2 x^4 \left(1121 \theta ^3+4977 \theta ^2+7951 \theta +4500\right) \\
+4 x^5 (\theta +2) \left(491 \theta ^2+1958 \theta +2082\right) \\
-x^6 (\theta +2) (\theta +3) (371 \theta +956)
-418 x^7 (\theta +2) (\theta +3) (\theta +4)
\end{matrix}
$
& $\begin{matrix} 1 \\ 3 \end{matrix}$ 
\tabularnewline
\hline 
${}{24B}$ 
&  
$ \begin{matrix}
%pf[24 B]=
\theta ^3
-2 x \left(\theta ^3+6 \theta ^2+4 \theta +1\right)
-4 x^2 \theta  \left(10 \theta ^2+12 \theta +5\right) \\
+8 x^3 \left(2 \theta ^3-7 \theta -6\right)
-16 x^4 \left(11 \theta ^3+48 \theta ^2+72 \theta +37\right) \\
+96 x^5 (\theta +2)^3
\end{matrix}
$
& $\begin{matrix} 2 \\ 2 \end{matrix}$ 
\tabularnewline
\hline 
${}{25A}$ 
&  
$ \begin{matrix}
%pf[25 A]=24 
\theta ^3
-8 x \left(5 \theta ^3+45 \theta ^2+33 \theta +9\right)
-8 x^2 \theta  \left(202 \theta ^2+231 \theta +104\right) \\
+4 x^3 \left(986 \theta ^3+3237 \theta ^2+4449 \theta +2393\right) \\
+4 x^4 \left(1142 \theta ^3+4527 \theta ^2+7967 \theta +5337\right) \\
-4 x^5 \left(3568 \theta ^3+21660 \theta ^2+48460 \theta +39175\right) \\
-2 x^6 \left(1900 \theta ^3+12900 \theta ^2+38249 \theta +44097\right) \\
+2 x^7 \left(7764 \theta ^3+70272 \theta ^2+220703 \theta +239610\right) \\
+2 x^8 (2 \theta +7) \left(932 \theta ^2+6272 \theta +12851\right) \\
-x^9 (2 \theta +7) (2 \theta +9) (1466 \theta +5963) \\
-275 x^{10} (2 \theta +7) (2 \theta +9) (2 \theta +11)
\end{matrix}
$
& $\begin{matrix} 3\; \\ 2^\dagger \end{matrix}$ 
\tabularnewline
\hline 
${}{26A}$ 
&  
$\begin{matrix}
%pf[26 A]=
13 \theta ^3
+x \left(503 \theta ^3-819 \theta ^2-611 \theta -169\right) \\
-x^2 \theta  \left(39505 \theta ^2+38904 \theta +17544\right) \\
+x^3 \left(835657 \theta ^3+2137227 \theta ^2+2234459 \theta +859173\right) \\
-4 x^4 \left(2141159 \theta ^3+8582394 \theta ^2+12553492 \theta +6539391\right)\\
+24 x^5 (4 \theta +7) \left(489557 \theta ^2+1799756 \theta +1778906\right) \\
-8 x^6 (4 \theta +7) (4 \theta +11) (1036418 \theta +2421687) \\
+2360544 x^7 (4 \theta +7) (4 \theta +11) (4 \theta +15)
\end{matrix}
$
& $\begin{matrix} 1 \\ 8 \end{matrix}$ 
\tabularnewline
\hline 
\end{tabular}

\begin{tabular}{|c|c|c|}
\hline
${}{27A}$ 
&  
$ \begin{matrix}
%pf[27 A]=
8 \theta ^3
+8 x \left(5 \theta ^3-15 \theta ^2-11 \theta -3\right)
-40 x^2 \theta  \left(25 \theta ^2+30 \theta +14\right) \\
+8 x^3 \left(341 \theta ^3+972 \theta ^2+1152 \theta +495\right) \\
-8 x^4 \left(494 \theta ^3+2199 \theta ^2+3877 \theta +2520\right)\\
+8 x^5 \left(254 \theta ^3+1545 \theta ^2+3541 \theta +2940\right)\\
+8 x^6 \left(175 \theta ^3+1260 \theta ^2+3099 \theta +2592\right)\\
-8 x^7 (\theta +3) \left(499 \theta ^2+2955 \theta +4676\right)\\
+8 x^8 (\theta +3) (\theta +4) (379 \theta +1310)
-1320 x^9 (\theta +3) (\theta +4) (\theta +5)
\end{matrix}
$
& $\begin{matrix} 3 \\ 2 \end{matrix}$ 
\tabularnewline
\hline 
${}{28B}$ 
&  
$ \begin{matrix}
%pf[28 B]=
\theta ^3
+x \left(8 \theta ^3-21 \theta ^2-15 \theta -4\right) 
-x^2 \theta  \left(247 \theta ^2+279 \theta +120\right)\\
+2 x^3 \left(615 \theta ^3+1761 \theta ^2+1914 \theta +752\right)\\
-16 x^4 \left(150 \theta ^3+669 \theta ^2+1041 \theta +562\right)
+1408 x^5 (\theta +2)^3
\end{matrix}
$
& $\begin{matrix} 2 \\ 3 \end{matrix}$ 
\tabularnewline
\hline 
${}{29A}$ 
&  
$ \begin{matrix}
%pf[29 A]=
9 \theta ^3
+3 x \left(7 \theta ^3-45 \theta ^2-33 \theta -9\right)
-3 x^2 \theta  \left(301 \theta ^2+348 \theta +158\right) \\
+3 x^3 \left(821 \theta ^3+2418 \theta ^2+2858 \theta +1252\right) \\
-3 x^4 \left(415 \theta ^3+1935 \theta ^2+3191 \theta +1851\right) \\
-3 x^5 \left(495 \theta ^3+2880 \theta ^2+5810 \theta +4042\right) \\
+6 x^6 (2 \theta +5) \left(173 \theta ^2+847 \theta +1095\right) \\
-444 x^7 (\theta +3) (2 \theta +5) (2 \theta +7)
\end{matrix}
$
& $\begin{matrix} 1 \\ 2 \end{matrix}$ 
\tabularnewline
\hline 
${}{30A}$ 
&  
$ \begin{matrix}
%pf[30 B]=
\theta ^3
+x \left(4 \theta ^3-21 \theta ^2-15 \theta -4\right)
-x^2 \theta  \left(199 \theta ^2+219 \theta +92\right) \\
+2 x^3 (\theta +1) \left(463 \theta ^2+905 \theta +606\right) \\
-12 x^4 (\theta +2) (2 \theta +3) (25 \theta +31)
-432 x^5 (\theta +2) (2 \theta +3) (2 \theta +5)
\end{matrix}
$
& $\begin{matrix} 1 \\ 3 \end{matrix}$ 
\tabularnewline
\hline 
${}{31AB}$ 
&  
$ \begin{matrix}
%pf[31 AB]=
11 \theta ^3
+x \left(742 \theta ^3-1023 \theta ^2-825 \theta -242\right) \\
-x^2 \theta  \left(79298 \theta ^2+78345 \theta +39414\right) \\
+2 x^3 \left(1308021 \theta ^3+3326079 \theta ^2+3698316 \theta +1516375\right) \\
+x^4 \left(-45916439 \theta ^3-183209175 \theta ^2-284070204 \theta -161450532\right) \\
+x^5 \left(488163212 \theta ^3+2638139415 \theta ^2+5180998833 \theta +3620750162\right) \\
+x^6 \left(-3258667723 \theta ^3-22178264979 \theta ^2-52548830340 \theta -43060503364\right) \\
+6 x^7 (3 \theta +8) \left(744793176 \theta ^2+4123737497 \theta +5940087425\right) \\ 
-9 x^8 (3 \theta +8) (3 \theta +11) (384663357 \theta +1256839126) \\
+1161721710 x^9 (3 \theta +8) (3 \theta +11) (3 \theta +14)
\end{matrix}
$
& $\begin{matrix} 1 \\ 9 \end{matrix}$ 
\tabularnewline
\hline 
${}{32A}$ 
&  
$\begin{matrix}
%pf[32 A]=
48 \theta ^3
-8 x \left(53 \theta ^3+63 \theta ^2+39 \theta +9\right)
+4 x^2 \theta  \left(68 \theta ^2+393 \theta +218\right) \\
+2 x^3 \left(3336 \theta ^3+1428 \theta ^2+21 \theta -368\right) \\
+x^4 \left(-36288 \theta ^3-42576 \theta ^2-21860 \theta +655\right) \\
+2 x^5 \left(51200 \theta ^3+85440 \theta ^2+63818 \theta +19745\right)  \\
-32 x^6 \left(5988 \theta ^3+12753 \theta ^2+11366 \theta +3837\right) \\
+8 x^7 \left(31360 \theta ^3+80544 \theta ^2+80689 \theta +28394\right) \\
-4 x^8 \left(57984 \theta ^3+172560 \theta ^2+187412 \theta +75565\right) \\ 
+8 x^9 \left(18000 \theta ^3+60432 \theta ^2+72006 \theta +31009\right) \\ 
-16 x^{10} \left(3376 \theta ^3+12588 \theta ^2+16588 \theta +7581\right) \\
+992 x^{11} (2 \theta +3)^3
\end{matrix}
$
& $\begin{matrix} 4\; \\ 2^\dagger \end{matrix}$ 
\tabularnewline
\hline 
${}{33AB}$ 
&  
$ \begin{matrix}
%pf[33 B]=
\theta ^3
+x \left(-13 \theta ^3-6 \theta ^2-4 \theta -1\right)
+3 x^2 \theta  \left(10 \theta ^2+12 \theta +5\right) \\
+x^3 \left(42 \theta ^3+108 \theta ^2+103 \theta +33\right)
+x^4 \left(109 \theta ^3+492 \theta ^2+789 \theta +446\right) \\
+x^5 (\theta +2) \left(19 \theta ^2+118 \theta +147\right) 
-2 x^6 (\theta +2) (\theta +3) (130 \theta +331) \\
+72 x^7 (\theta +2) (\theta +3) (\theta +4)
\end{matrix}
$
& $\begin{matrix} 1 \\ 1 \end{matrix}$ 
\tabularnewline
\hline 
\end{tabular}

\begin{tabular}{|c|c|c|}
\hline
${}{34A}$ 
&  
$ \begin{matrix}
%pf[34 A]=
112 \theta^3
+8 x \left(43 \theta ^3-231 \theta ^2-175 \theta -49\right)
-4 x^2 \theta  \left(3466 \theta ^2+3897 \theta +1810\right) \\
+2 x^3 \left(24372 \theta ^3+73032 \theta ^2+87623 \theta  +39124\right) \\
+x^4 \left(-10984 \theta ^3-63540 \theta ^2-117998 \theta -74239\right) \\
-56 x^5 \left(2236 \theta ^3+13224 \theta ^2+26961 \theta +18902\right) \\
+3152 x^6 (2 \theta +5) (4 \theta +9) (4 \theta +11)
\end{matrix}
$
& $\begin{matrix} 1\; \\ 2^\dagger \end{matrix}$ 
\tabularnewline
\hline 
${}{35A}$ 
&  
$ \begin{matrix}
%pf[35 A]=
\theta ^3
+x \left(3 \theta ^3-6 \theta ^2-4 \theta -1\right)
-x^2 \theta  \left(26 \theta ^2+36 \theta +17\right) \\
+x^3 \left(-18 \theta ^3-102 \theta ^2-167 \theta -101\right) \\
+x^4 \left(-251 \theta ^3-1170 \theta ^2-2053 \theta -1308\right) \\
-x^5 (\theta +2) \left(253 \theta ^2+1108 \theta +1383\right) \\
-4 x^6 (\theta +2) (\theta +3) (131 \theta +356)
-532 x^7 (\theta +2) (\theta +3) (\theta +4)
\end{matrix}
$
& $\begin{matrix} 1 \\ 1 \end{matrix}$ 
\tabularnewline
\hline 
${}{36A}$ 
&  
$ \begin{matrix}
%pf[36 A]=
\theta ^3
-x (3 \theta +1) \left(3 \theta ^2+2 \theta +1\right)
-6 x^2 \theta  \left(12 \theta ^2-3 \theta -1\right) \\
+2 x^3 \theta  \left(284 \theta ^2+405 \theta +199\right)
+6 x^4 \theta  \left(1156 \theta ^2+75 \theta +89\right) \\
-6 x^5 \theta  \left(11927 \theta ^2+10401 \theta +4939\right) \\
+18 x^6 \left(8968 \theta ^3+11586 \theta ^2+5960 \theta +2553\right) \\ 
+18 x^7 \left(11788 \theta ^3+14184 \theta ^2-5086 \theta -19947\right) \\
-27 x^8 \left(30109 \theta ^3+44628 \theta ^2+7040 \theta -6990\right) \\
-27 x^9 \left(19871 \theta ^3+39147 \theta ^2+9715 \theta +29949\right) \\
+486 x^{10} \left(2664 \theta ^3+4503 \theta ^2+2623 \theta +561\right) \\
+486 x^{11} \left(2892 \theta ^3+6453 \theta ^2+5465 \theta +1657\right) 
+360126 x^{12} (\theta +1)^3
\end{matrix}
$
& $\begin{matrix} 4 \\ 2 \end{matrix}$ 
\tabularnewline
\hline 
${}{38A}$ 
&  
$ \begin{matrix}
%pf[38 A]=
81 \theta ^3
+9 x \left(239 \theta ^3-351 \theta ^2-279 \theta -81\right) \\ 
+4 x^2 \theta  \left(17725 \theta ^2-24678 \theta -12195\right) \\
-x^3 \theta  \left(3086561 \theta ^2+3244998 \theta +1524970\right) \\
+x^4 \left(37779277 \theta ^3+104515257 \theta ^2+118146733 \theta +48647753\right) \\
-4 x^5 \left(57515231 \theta ^3+251097714 \theta ^2+406894475 \theta +237268425\right) \\
+2 x^6 \left(394121816 \theta ^3+2334406791 \theta ^2+4853928845 \theta +3517728050\right) \\
-2 x^7 (2 \theta +5) \left(382331563 \theta ^2+1896349895 \theta +2463921642\right) \\
+152 x^8 (2 \theta +5) (2 \theta +7) (2520778 \theta +7524021) \\
-72570120 x^9 (2 \theta +5) (2 \theta +7) (2 \theta +9)
\end{matrix}
$
& $\begin{matrix} 1 \\ 4 \end{matrix}$ 
\tabularnewline
\hline 
${}{39A}$ 
&  
$ \begin{matrix}
%pf[39 A]=
27 \theta ^3
+9 x \left(38 \theta ^3-81 \theta ^2-63 \theta -18\right) \\
+3 x^2 \theta  \left(2641 \theta ^2-4347 \theta -2070\right) \\ 
+x^3 \theta  \left(171998 \theta ^2-284013 \theta -131229\right) \\
-69 x^4 \theta  \left(85495 \theta ^2+89508 \theta +41008\right) \\
+8 x^5 \left(5597701 \theta ^3+15741750 \theta ^2+17627454 \theta +7229767\right) \\
-138 x^6 \left(938927 \theta ^3+4168992 \theta ^2+6623258 \theta +3718384\right) \\
+4 x^7 \left(28364771 \theta ^3+170987367 \theta ^2+343956795 \theta +230941456\right) \\
+2640216 x^8 (2 \theta +5) (3 \theta +7) (3 \theta +8)
\end{matrix}
$
& $\begin{matrix} 1 \\ 3 \end{matrix}$ 
\tabularnewline
\hline 
${}{41A}$ 
&  
$\begin{matrix}
%pf[41 A]=
\theta ^3
+4 x \left(9 \theta ^3-12 \theta ^2-10 \theta -3\right)
-8 x^2 \theta  \left(217 \theta ^2+243 \theta +128\right) \\
+x^3 \left(26494 \theta ^3+74343 \theta ^2+88577 \theta +38487\right) \\
-4 x^4 \left(54068 \theta ^3+237045 \theta ^2+400047 \theta +246361\right)  \\
+4 x^5 \left(270602 \theta ^3+1603767 \theta ^2+3459565 \theta +2660558\right) \\
+x^6 \left(-3489323 \theta ^3-26006067 \theta ^2-68068153 \theta -62098537\right) \\
+4 x^7 \left(1821615 \theta ^3+16341414 \theta ^2+50218666 \theta +52756576\right) \\
-16 x^8 (2 \theta +7) \left(297749 \theta ^2+2079338 \theta +3728485\right) \\
+8 x^9 (2 \theta +7) (2 \theta +9) (221314 \theta +883687) \\
-284512 x^{10} (2 \theta +7) (2 \theta +9) (2 \theta +11)
\end{matrix}
$
& $\begin{matrix} 1 \\ 4 \end{matrix}$ 
\tabularnewline
\hline 
${}{42A}$ 
&  
$\begin{matrix}
%pf[42 A]=
48 \theta ^3
-8 x \left(53 \theta ^3+45 \theta ^2+33 \theta +9\right) 
+4 x^2 \theta  \left(194 \theta ^2+237 \theta +112\right) \\
+2 x^3 \left(588 \theta ^3+2232 \theta ^2+3379 \theta +1870\right) \\
+x^4 \left(-2104 \theta ^3-9036 \theta ^2-13898 \theta -7449\right) \\
+4 x^5 (\theta +2) \left(60 \theta ^2+744 \theta +1163\right) \\
-16 x^6 (\theta +2) (\theta +3) (442 \theta +1285) 
+7360 x^7 (\theta +2) (\theta +3) (\theta +4)
\end{matrix}
$
& $\begin{matrix} 1\; \\ 1^\dagger \end{matrix}$ 
\tabularnewline
\hline 
\end{tabular}

\begin{tabular}{|c|c|c|}
\hline
${}{44AB}$ 
&  
$\begin{matrix}
%pf[44 AB]=
\theta ^3
+2 x \left(6 \theta ^3-9 \theta ^2-7 \theta -2\right)
-8 x^2 \theta  \left(29 \theta ^2+33 \theta +16\right) \\
+2 x^3 \left(598 \theta ^3+1683 \theta ^2+1935 \theta +810\right) \\
-16 x^4 \left(208 \theta ^3+915 \theta ^2+1505 \theta +894\right) \\
+32 x^5 \left(162 \theta ^3+963 \theta ^2+2021 \theta +1484\right) \\
-64 x^6 \left(71 \theta ^3+531 \theta ^2+1353 \theta +1173\right) 
+1536 x^7 (\theta +3)^3
\end{matrix}
$
& $\begin{matrix} 2 \\ 2 \end{matrix}$ 
\tabularnewline
\hline 
${}{45A}$ 
&  
$ \begin{matrix}
%pf[45 A]=
16 \theta ^3
-8 x \left(23 \theta ^3+15 \theta ^2+11 \theta +3\right)
+12 x^2 \theta  \left(38 \theta ^2+49 \theta +24\right) \\
+2 x^3 \left(424 \theta ^3+1596 \theta ^2+2633 \theta +1614\right) \\
-3 x^4 \left(1000 \theta ^3+4596 \theta ^2+8270 \theta +5493\right) \\
+6 x^5 \left(88 \theta ^3+348 \theta ^2+503 \theta +282\right) \\
+36 x^6 \left(142 \theta ^3+1083 \theta ^2+2807 \theta +2469\right)
-2808 x^7 (\theta +3)^3
\end{matrix}
$
& $\begin{matrix} 2\; \\ 1^\dagger \end{matrix}$ 
\tabularnewline
\hline 
${}{46CD}$ 
&  
$\begin{matrix}
%pf[46 CD]=
\theta ^3
-x (2 \theta +1) \left(3 \theta ^2+3 \theta +2\right)
+x^2 \theta  (\theta +1) (\theta +2) \\
+x^3 (2 \theta +3) \left(15 \theta ^2+45 \theta +52\right)
-2 x^4 (\theta +2) \left(17 \theta ^2+68 \theta +84\right) \\
-10 x^5 (2 \theta +5) \left(\theta ^2+5 \theta +8\right)
+5 x^6 (\theta +3) \left(9 \theta ^2+54 \theta +88\right) \\
-2 x^7 (\theta +3) (\theta +4) (2 \theta +7)
-20 x^8 (\theta +3) (\theta +4) (\theta +5)
\end{matrix}
$
& $\begin{matrix} 1 \\ 1 \end{matrix}$ 
\tabularnewline
\hline 
${}{47AB}$ 
&  
$\begin{matrix}
%pf[47 AB]=
\theta ^3
+x \left(6 \theta ^3-9 \theta ^2-7 \theta -2\right)
-x^2 \theta  \left(61 \theta ^2+75 \theta +38\right) \\
+12 x^3 \left(9 \theta ^3+24 \theta ^2+29 \theta +12\right)  \\
+x^4 \left(-269 \theta ^3-1182 \theta ^2-2184 \theta -1504\right) \\
+2 x^5 \left(106 \theta ^3+624 \theta ^2+1418 \theta +1179\right) \\
+x^6 \left(-205 \theta ^3-1557 \theta ^2-4358 \theta -4380\right) \\ 
+x^7 \left(-126 \theta ^3-1089 \theta ^2-3341 \theta -3598\right) \\ 
+2 x^8 \left(161 \theta ^3+1662 \theta ^2+5958 \theta +7388\right) \\
-2 x^9 (\theta +4) \left(214 \theta ^2+1691 \theta +3477\right) \\
+x^{10} (\theta +4) (\theta +5) (325 \theta +1446) \\
-132 x^{11} (\theta +4) (\theta +5) (\theta +6)
\end{matrix}
$
& $\begin{matrix} 1 \\ 1 \end{matrix}$ 
\tabularnewline
\hline 
${}{49Z}$ 
&  
$ \begin{matrix}
%pf[49 Z]=
432 \theta ^3
-72 x \left(39 \theta ^3+75 \theta ^2+63 \theta +19\right)
-36 x^2 \theta  \left(30 \theta ^2-81 \theta -76\right) \\
+2 x^3 \left(14796 \theta ^3+72036 \theta ^2+148941 \theta +116234\right) \\
+x^4 \left(-13824 \theta ^3-150984 \theta ^2-451062 \theta -462041\right) \\
-2 x^5 \left(64800 \theta ^3+504792 \theta ^2+1516734 \theta +1680341\right) \\
+3 x^6 \left(27432 \theta ^3+346140 \theta ^2+1418634 \theta +1990651\right) \\ 
+18 x^7 \left(16992 \theta ^3+181752 \theta ^2+704575 \theta +974633\right)  \\
+x^8 \left(-184032 \theta ^3-2862648 \theta ^2-14756022 \theta -25753097\right)\\
-6 x^9 \left(68868 \theta ^3+936468 \theta ^2+4429179 \theta +7248311\right) \\
+6 x^{10} \left(33300 \theta ^3+624726 \theta ^2+3879420 \theta +8061365\right) \\
+2 x^{11} \left(152496 \theta ^3+2510892 \theta ^2+14003784 \theta +26380243\right) \\
-9 x^{12} \left(12168 \theta ^3+265764 \theta ^2+1918830 \theta +4601999\right)\\
-6 x^{13} \left(17172 \theta ^3+327960 \theta ^2+2091075 \theta +4440746\right)\\
+x^{14} (6 \theta +43) \left(5652 \theta ^2+93792 \theta +384547\right)\\
+4 x^{15} (6 \theta +43) (6 \theta +49) (84 \theta +473) \\
-19 x^{16} (6 \theta +43) (6 \theta +49) (6 \theta +55)
\end{matrix}
$
& $\begin{matrix} 4\; \\ 1^\dagger \end{matrix}$ 
\tabularnewline
\hline 
${}{50A}$ 
&  
$\begin{matrix}
%pf[50 A]=
1216 \theta ^3
+64 x \left(1071 \theta ^3-1311 \theta ^2-1159 \theta -361\right) \\
-640 x^2 \theta  \left(6938 \theta ^2+7977 \theta +4540\right) \\
+80 x^3 \left(1235916 \theta ^3+3540549 \theta ^2+4439647 \theta +2016171\right)\\
-80 x^4 \left(15527084 \theta ^3+69650223 \theta ^2+123709031 \theta +80706657\right)\\
+32 x^5 \left(316607084 \theta ^3+1924672605 \theta ^2+4377358125 \theta +3597085400\right) \\
-4 x^6 \left(14265057072 \theta ^3+109312727208 \theta ^2+302353965627 \theta +296792076433\right) \\
+20 x^7 \left(11416566608 \theta ^3+105538882728 \theta ^2+343613961789 \theta +390817201500\right) \\
-40 x^8 \left(16371371232 \theta ^3+177239719848 \theta ^2+663865518584 \theta +857151934777\right) \\
+15 x^9 \left(89056576448 \theta ^3+1105164713776 \theta ^2+4681291119652 \theta +6759099722789\right) \\
-x^{10} \big(1888277885504 \theta ^3+26426519949360 \theta ^2 +124906791611500 \theta +199331847523425 \big) \\
+22 x^{11} (4 \theta +21) \left(19909272688 \theta ^2+205694511960 \theta +538657069103\right) \\
-4840 x^{12} (4 \theta +21) (4 \theta +25) (12327956 \theta +69871153) \\
+3588908400 x^{13} (4 \theta +21) (4 \theta +25) (4 \theta +29)
\end{matrix}
$
& $\begin{matrix} 3\; \\ 4^\dagger \end{matrix}$ 
\tabularnewline
\hline 
\end{tabular}

\begin{tabular}{|c|c|c|}
\hline
${}{51A}$ 
&  
$\begin{matrix}
%pf[51 A]=
\theta ^3
+x \left(2 \theta ^3-9 \theta ^2-7 \theta -2\right) 
-x^2 \theta  \left(39 \theta ^2+45 \theta +22\right) \\ 
+x^3 \left(58 \theta ^3+153 \theta ^2+163 \theta +59\right) \\
-2 x^4 \left(46 \theta ^3+192 \theta ^2+305 \theta +174\right) \\ 
+x^5 (\theta +2) \left(158 \theta ^2+629 \theta +801\right) \\
-x^6 (\theta +3) \left(35 \theta ^2+186 \theta +259\right) \\
-x^7 (\theta +3) \left(134 \theta ^2+777 \theta +1180\right) \\
+3 x^8 (\theta +3) (\theta +4) (43 \theta +143) 
-120 x^9 (\theta +3) (\theta +4) (\theta +5)
\end{matrix}
$
& $\begin{matrix} 1 \\ 1 \end{matrix}$ 
\tabularnewline
\hline 
${}{54A}$ 
&  
$\begin{matrix}
%pf[54 A]=
24 \theta ^3
+8 x \left(85 \theta ^3-99 \theta ^2-87 \theta -27\right)
-8 x^2 \theta  \left(2665 \theta ^2+3012 \theta +1706\right) \\
+12 x^3 \left(19412 \theta ^3+54285 \theta ^2+68019 \theta +30869\right) \\
-12 x^4 \left(121898 \theta ^3+532521 \theta ^2+943777 \theta +618045\right) \\
+12 x^5 \left(506870 \theta ^3+2995185 \theta ^2+6773107 \theta +5586387\right)\\
-18 x^6 \left(987400 \theta ^3+7341870 \theta ^2+20105135 \theta +19755909\right) \\
+18 x^7 \left(2072236 \theta ^3+18555768 \theta ^2+59548649 \theta +67524810\right) \\
-18 x^8 \left(3102556 \theta ^3+32481420 \theta ^2+119405669 \theta +152980491\right) \\
+81 x^9 \left(703184 \theta ^3+8425716 \theta ^2+34889540 \theta +49743479\right) \\
-27 x^{10} \left(1255512 \theta ^3+16946388 \theta ^2+78041746 \theta +122430879\right) \\
+81 x^{11} \left(36920 \theta ^3+558444 \theta ^2+2849282 \theta +4901357\right) \\
+81 x^{12} (2 \theta +11) \left(76764 \theta ^2+842676 \theta +2345251\right) \\
-1944 x^{13} (2 \theta +11) (2 \theta +13) (1202 \theta +7203) \\
+293544 x^{14} (2 \theta +11) (2 \theta +13) (2 \theta +15)
\end{matrix}
$
& $\begin{matrix} 3\; \\ 2^\dagger \end{matrix}$ 
\tabularnewline
\hline 
${}{55A}$ 
&  
$\begin{matrix}
%pf[55 A]=
\theta ^3
-x (5 \theta +2) \left(2 \theta ^2+\theta +1\right)
+3 x^2 \theta  \left(47 \theta ^2+11 \theta +6\right) \\
-x^3 \theta  \left(666 \theta ^2+879 \theta +479\right) \\
+x^4 \left(-269 \theta ^3-1032 \theta ^2-2372 \theta -1704\right) \\
+4 x^5 \left(1099 \theta ^3+4881 \theta ^2+8610 \theta +5592\right) \\
+6 x^6 \left(71 \theta ^3+496 \theta ^2+1676 \theta +1960\right) \\
-2 x^7 \left(4878 \theta ^3+36435 \theta ^2+94555 \theta +84860\right) \\ 
-8 x^8 (\theta +3) \left(435 \theta ^2+2643 \theta +4084\right) \\
+2596 x^9 (\theta +3) (\theta +4) (2 \theta +7)
\end{matrix}
$
& $\begin{matrix} 1 \\ 1 \end{matrix}$ 
\tabularnewline
\hline 
${}{56A}$ 
&  
$\begin{matrix}
%pf[56 A]=
\theta ^3
+x \left(2 \theta ^3-9 \theta ^2-7 \theta -2\right) 
-x^2 \theta  \left(39 \theta ^2+45 \theta +22\right) \\
+8 x^3 \left(7 \theta ^3+18 \theta ^2+17 \theta +4\right) \\ 
-4 x^4 \left(22 \theta ^3+84 \theta ^2+133 \theta +74\right) \\
+16 x^5 \left(18 \theta ^3+105 \theta ^2+237 \theta +196\right) \\
-4 x^6 \left(55 \theta ^3+399 \theta ^2+1062 \theta +1016\right) \\
+8 x^7 \left(30 \theta ^3+261 \theta ^2+777 \theta +790\right) \\
-16 x^8 \left(23 \theta ^3+240 \theta ^2+844 \theta +1000\right) 
+128 x^9 (\theta +4)^3
\end{matrix}
$
& $\begin{matrix} 2 \\ 1 \end{matrix}$ 
\tabularnewline
\hline 
${}{59AB}$ 
&  
$\begin{matrix}
%pf[59 AB]=
3 \theta ^3
+x \left(19 \theta ^3-36 \theta ^2-30 \theta -9\right)
-x^2 \theta  \left(278 \theta ^2+318 \theta +169\right) \\
+x^3 \left(862 \theta ^3+2460 \theta ^2+3081 \theta +1426\right) \\
+x^4 \left(-1195 \theta ^3-5364 \theta ^2-9433 \theta -6153\right) \\
+3 x^5 \left(83 \theta ^3+558 \theta ^2+995 \theta +485\right) \\
+x^6 \left(1534 \theta ^3+11226 \theta ^2+31787 \theta +33036\right) \\
+x^7 \left(-2546 \theta ^3-22734 \theta ^2-73553 \theta -84696\right) \\
+x^8 \left(1711 \theta ^3+17952 \theta ^2+65258 \theta +81885\right) \\
-3 x^9 \left(67 \theta ^3+840 \theta ^2+3240 \theta +3842\right) \\
-2 x^{10} \left(486 \theta ^3+6516 \theta ^2+30400 \theta +49119\right) \\
+2 x^{11} (\theta +5) \left(424 \theta ^2+4231 \theta +10982\right) \\
-2 x^{12} (\theta +5) (\theta +6) (98 \theta +557)
-344 x^{13} (\theta +5) (\theta +6) (\theta +7)
\end{matrix}
$
& $\begin{matrix} 1 \\ 1 \end{matrix}$ 
\tabularnewline
\hline 
${}{60B}$ 
&  
$\begin{matrix}
%pf[60 B]=
\theta ^3
+x \left(10 \theta ^3-9 \theta ^2-7 \theta -2\right)
-x^2 \theta  \left(79 \theta ^2+93 \theta +46\right) \\
+8 x^3 \left(28 \theta ^3+75 \theta ^2+85 \theta +34\right) \\
-4 x^4 \left(175 \theta ^3+762 \theta ^2+1285 \theta +794\right) \\
+16 x^5 \left(62 \theta ^3+363 \theta ^2+761 \theta +564\right) \\
-64 x^6 \left(23 \theta ^3+171 \theta ^2+433 \theta +373\right) 
+1024 x^7 (\theta +3)^3
\end{matrix}
$
& $\begin{matrix} 2 \\ 1 \end{matrix}$ 
\tabularnewline
\hline 
\end{tabular}

\begin{tabular}{|c|c|c|}
\hline
${}{62AB}$ 
&  
$\begin{matrix}
%pf[62 AB]=
3 \theta ^3
+x \left(19 \theta ^3-36 \theta ^2-30 \theta -9\right)
-x^2 \theta  \left(278 \theta ^2+318 \theta +169\right) \\
+x^3 \left(856 \theta ^3+2433 \theta ^2+3000 \theta +1345\right) \\
+x^4 \left(-1245 \theta ^3-5535 \theta ^2-9958 \theta -6618\right) \\
+x^5 \left(747 \theta ^3+4635 \theta ^2+11534 \theta +10576\right) \\
+x^6 \left(1195 \theta ^3+8562 \theta ^2+21563 \theta +18864\right) \\
+x^7 \left(-3473 \theta ^3-30726 \theta ^2-94787 \theta -101358\right) \\
+x^8 \left(4753 \theta ^3+49425 \theta ^2+176490 \theta +215912\right) \\
-x^9 (\theta +4) \left(4433 \theta ^2+35212 \theta +72459\right) \\
+4 x^{10} (\theta +4) (\theta +5) (584 \theta +2601)
-1032 x^{11} (\theta +4) (\theta +5) (\theta +6)
\end{matrix}
$
& $\begin{matrix} 1 \\ 1 \end{matrix}$ 
\tabularnewline
\hline 
${}{66A}$ 
&  
$\begin{matrix}
%pf[66 A]=
1024 \theta ^3
+256 x \left(13 \theta ^3-66 \theta ^2-54 \theta -16\right)
-128 x^2 \theta  \left(970 \theta ^2+1029 \theta +515\right) \\
+128 x^3 \left(4033 \theta ^3+12078 \theta ^2+14875 \theta +7017\right) \\
-32 x^4 \left(5918 \theta ^3+34311 \theta ^2+52631 \theta +26021\right) \\
-48 x^5 \left(45932 \theta ^3+271530 \theta ^2+601769 \theta +484188\right) \\
+8 x^6 \left(262162 \theta ^3+2012157 \theta ^2+5256965 \theta +4681203\right)\\
+24 x^7 (\theta +3) \left(139021 \theta ^2+824433 \theta +1358168\right) \\
-18 x^8 (\theta +3) (\theta +4) (131788 \theta +474773)
-2311065 x^9 (\theta +3) (\theta +4) (\theta +5)
\end{matrix}
$
& $\begin{matrix} 1\; \\ {\frac{3}{2}}\dagger \end{matrix}$ 
\tabularnewline
\hline 
${}{69AB}$ 
&  
$\begin{matrix}
%pf[69 AB]=
3 \theta ^3
+x \left(11 \theta ^3-36 \theta ^2-30 \theta -9\right)
-5 x^2 \theta  \left(44 \theta ^2+48 \theta +25\right) \\
+x^3 \left(652 \theta ^3+1884 \theta ^2+2295 \theta +1054\right) \\
+x^4 \left(-566 \theta ^3-2556 \theta ^2-3997 \theta -2193\right) \\
+x^5 \left(-52 \theta ^3-285 \theta ^2-1346 \theta -1950\right) \\
+3 x^6 \left(38 \theta ^3+237 \theta ^2+515 \theta +368\right) \\
-4 x^7 \left(238 \theta ^3+2079 \theta ^2+6463 \theta +7077\right) \\
+x^8 \left(1987 \theta ^3+20904 \theta ^2+75449 \theta +93240\right) \\
+x^9 (\theta +4) \left(201 \theta ^2+1437 \theta +3560\right) \\
-x^{10} (\theta +4) (\theta +5) (1363 \theta +6183)
-385 x^{11} (\theta +4) (\theta +5) (\theta +6)
\end{matrix}
$
& $\begin{matrix} 1 \\ 1 \end{matrix}$ 
\tabularnewline
\hline 
${}{70A}$ 
&  
$\begin{matrix}
%pf[70 A]=
400 \theta ^3
+40 x \left(27 \theta ^3-105 \theta ^2-85 \theta -25\right)
+4 x^2 \theta  \left(5084 \theta ^2-5985 \theta -3020\right) \\
-6 x^3 \theta  \left(39692 \theta ^2+45558 \theta +23781\right) \\
+x^4 \left(572304 \theta ^3+1662552 \theta ^2+1994356 \theta +874283\right) \\
+x^5 \left(-105208 \theta ^3-720336 \theta ^2-1296832 \theta -828769\right) \\
-5 x^6 \left(412240 \theta ^3+2398596 \theta ^2+5290012 \theta +4255295\right)\\
+x^7 \left(2585736 \theta ^3+19390140 \theta ^2+50625486 \theta +45792265\right)\\
+64 x^8 (\theta +3) \left(5624 \theta ^2+26469 \theta +36840\right) \\
-1675840 x^9 (\theta +3) (\theta +4) (2 \theta +7)
\end{matrix}
$
& $\begin{matrix} 1\; \\ 1^\dagger \end{matrix}$ 
\tabularnewline
\hline 
${}{71AB}$ 
&  
$\begin{matrix}
%pf[71 AB]=
\theta ^3
+x \left(8 \theta ^3-15 \theta ^2-13 \theta -4\right)
-x^2 \theta  \left(143 \theta ^2+159 \theta +88\right) \\
+3 x^3 \left(200 \theta ^3+567 \theta ^2+701 \theta +318\right) \\ 
-2 x^4 \left(561 \theta ^3+2475 \theta ^2+4271 \theta +2698\right) \\
+12 x^5 \left(95 \theta ^3+564 \theta ^2+1274 \theta +1045\right) \\
-3 x^6 \left(283 \theta ^3+2115 \theta ^2+6174 \theta +6624\right) \\
+2 x^7 \left(155 \theta ^3+1425 \theta ^2+5154 \theta +6847\right) \\
+2 x^8 \left(327 \theta ^3+3384 \theta ^2+12380 \theta +15854\right) \\
-3 x^9 \left(374 \theta ^3+4455 \theta ^2+18525 \theta +26724\right) \\
+2 x^{10} \left(565 \theta ^3+7584 \theta ^2+35425 \theta +57286\right) \\
-4 x^{11} \left(251 \theta ^3+3750 \theta ^2+19231 \theta +33768\right) \\
+3 x^{12} \left(231 \theta ^3+3798 \theta ^2+21316 \theta +40776\right) \\
-3 x^{13} (\theta +6) \left(160 \theta ^2+1917 \theta +5885\right) \\
+x^{14} (\theta +6) (\theta +7) (101 \theta +646)
-126 x^{15} (\theta +6) (\theta +7) (\theta +8)
\end{matrix}
$
& $\begin{matrix} 1 \\ 1 \end{matrix}$ 
\tabularnewline
\hline 
\end{tabular}

\begin{tabular}{|c|c|c|}
\hline
${}{78A}$ 
&  
$\begin{matrix}
%pf[78 A]=
9 \theta ^3
+3 x \left(11 \theta ^3-36 \theta ^2-30 \theta -9\right)
+x^2 \theta  \left(313 \theta ^2-720 \theta -375\right) \\
-x^3 \theta  \left(5846 \theta ^2+6105 \theta +3088\right) \\
+x^4 \left(17642 \theta ^3+50307 \theta ^2+59625 \theta +26153\right) \\
+x^5 \left(-16566 \theta ^3-72963 \theta ^2-118994 \theta -68817\right) \\
+x^6 \left(8513 \theta ^3+51375 \theta ^2+123130 \theta +109230\right) \\
+x^7 \left(5985 \theta ^3+37107 \theta ^2+87299 \theta +77042\right) \\
+x^8 \left(-51959 \theta ^3-460107 \theta ^2-1402813 \theta -1468803\right) \\
+14 x^9 (2 \theta +7) \left(1766 \theta ^2+12092 \theta +21363\right) \\
-11676 x^{10} (\theta +4) (2 \theta +7) (2 \theta +9)
\end{matrix}
$
& $\begin{matrix} 1 \\ 1 \end{matrix}$ 
\tabularnewline
\hline 
${}{87AB}$ 
&  
$\begin{matrix}
%pf[87 AB]=
4 \theta ^3
+4 x \left(16 \theta ^3-15 \theta ^2-13 \theta -4\right)
+x^2 \theta  \left(915 \theta ^2-972 \theta -544\right) \\
-x^3 \theta  \left(13150 \theta ^2+14049 \theta +7331\right) \\
+3 x^4 \left(23581 \theta ^3+63823 \theta ^2+78458 \theta +34528\right) \\
+x^5 \left(-249716 \theta ^3-1070154 \theta ^2-1889116 \theta -1233207\right) \\
+2 x^6 \left(316163 \theta ^3+1844226 \theta ^2+4150643 \theta +3415212\right) \\
+x^7 \left(-1229502 \theta ^3-9062481 \theta ^2-24684291 \theta -24165398\right) \\
+x^8 \left(1833562 \theta ^3+16319166 \theta ^2+51995469 \theta +58496871\right) \\
+x^9 \left(-2095894 \theta ^3-21855723 \theta ^2-79694435 \theta -100976538\right) \\
+x^{10} \left(1721047 \theta ^3+20572332 \theta ^2+84371249 \theta +118444956\right) \\
-12 x^{11} \left(75554 \theta ^3+1018933 \theta ^2+4645483 \theta +7156124\right) \\
+2 x^{12} (\theta +5) \left(69349 \theta ^2+699736 \theta +1749345\right) \\
+18 x^{13} (\theta +5) (\theta +6) (7914 \theta +43085) \\
-89709 x^{14} (\theta +5) (\theta +6) (\theta +7)
\end{matrix}
$
& $\begin{matrix} 1 \\ 1 \end{matrix}$ 
\tabularnewline
\hline 
${}{92AB}$ 
&  
$\begin{matrix}
%pf[92 AB]=
\theta ^3
+x \left(16 \theta ^3-15 \theta ^2-13 \theta -4\right)
-x^2 \theta  \left(211 \theta ^2+243 \theta +136\right) \\
+x^3 \left(1108 \theta ^3+3075 \theta ^2+3857 \theta +1732\right) \\
-2 x^4 \left(1949 \theta ^3+8457 \theta ^2+15199 \theta +10108\right) \\
+4 x^5 \left(2511 \theta ^3+14757 \theta ^2+33833 \theta +28464\right) \\
+x^6 \left(-20123 \theta ^3-148971 \theta ^2-412614 \theta -413176\right) \\
+2 x^7 \left(16089 \theta ^3+143565 \theta ^2+464766 \theta +535568\right) \\
-8 x^8 \left(5147 \theta ^3+53730 \theta ^2+198658 \theta +257548\right) \\
+8 x^9 \left(5268 \theta ^3+62967 \theta ^2+261601 \theta +375926\right) \\
-16 x^{10} \left(2109 \theta ^3+28398 \theta ^2+130937 \theta +206306\right) \\
+32 x^{11} \left(627 \theta ^3+9390 \theta ^2+47577 \theta +81506\right) \\
-64 x^{12} \left(129 \theta ^3+2127 \theta ^2+11761 \theta +21806\right) \\
+1664 x^{13} (\theta +6)^3
\end{matrix}
$
& $\begin{matrix} 2 \\ 1 \end{matrix}$ 
\tabularnewline
\hline 
${}{94AB}$ 
&  
$\begin{matrix}
%pf[94 AB]=
5 \theta ^3
+x \left(51 \theta ^3-90 \theta ^2-80 \theta -25\right)
-x^2 \theta  \left(1052 \theta ^2+1158 \theta +661\right) \\
+x^3 \left(5636 \theta ^3+16038 \theta ^2+20283 \theta +9531\right) \\
+x^4 \left(-12726 \theta ^3-57042 \theta ^2-98449 \theta -62768\right) \\
+x^5 \left(5020 \theta ^3+32655 \theta ^2+54742 \theta +22049\right) \\
+x^6 \left(34568 \theta ^3+254835 \theta ^2+742308 \theta +799590\right) \\
+x^7 \left(-65090 \theta ^3-585210 \theta ^2-1900367 \theta -2193581\right) \\
+x^8 \left(12505 \theta ^3+136740 \theta ^2+424739 \theta +354504\right) \\
+x^9 \left(77745 \theta ^3+928845 \theta ^2+3976036 \theta +6015924\right) \\
+x^{10} \left(-62251 \theta ^3-842856 \theta ^2-3902055 \theta -6167414\right) \\
+x^{11} \left(-34429 \theta ^3-513300 \theta ^2-2687519 \theta -4907160\right) \\
+3 x^{12} \left(15621 \theta ^3+257979 \theta ^2+1439906 \theta +2714784\right) \\
+x^{13} (\theta +6) \left(9683 \theta ^2+115236 \theta +369481\right) \\
-8 x^{14} (\theta +6) (\theta +7) (1801 \theta +11744) \\
-4440 x^{15} (\theta +6) (\theta +7) (\theta +8)
\end{matrix}
$
& $\begin{matrix} 1 \\ 1 \end{matrix}$ 
\tabularnewline
\hline 
\end{tabular}

\begin{tabular}{|c|c|c|}
\hline
${}{95AB}$ 
&  
$ \begin{matrix}
%pf[95 AB]=
\theta ^3
+x \left(6 \theta ^3-15 \theta ^2-13 \theta -4\right)
-x^2 \theta  \left(125 \theta ^2+135 \theta +74\right) \\
+2 x^3 \left(262 \theta ^3+759 \theta ^2+953 \theta +455\right) \\
+x^4 \left(-655 \theta ^3-3066 \theta ^2-5056 \theta -3056\right) \\
-2 x^5 \left(514 \theta ^3+2907 \theta ^2+7457 \theta +7226\right) \\
+x^6 \left(3663 \theta ^3+27303 \theta ^2+77642 \theta +80596\right) \\
+x^7 \left(-2082 \theta ^3-19149 \theta ^2-62397 \theta -71486\right) \\
+x^8 \left(-4031 \theta ^3-41796 \theta ^2-155776 \theta -205600\right) \\
+2 x^9 \left(2816 \theta ^3+33780 \theta ^2+141378 \theta +205211\right) \\
-4 x^{10} \left(107 \theta ^3+1509 \theta ^2+7819 \theta +14385\right) \\
-6 x^{11} (\theta +5) \left(442 \theta ^2+4395 \theta +10938\right) \\
+120 x^{12} (\theta +5) (\theta +6) (13 \theta +71)
-640 x^{13} (\theta +5) (\theta +6) (\theta +7)
\end{matrix}
$
& $\begin{matrix} 1 \\ 1 \end{matrix}$ 
\tabularnewline
\hline 
${}{105A}$ 
&  
$\begin{matrix}
%pf[105 A]=
112 \theta ^3
+8 x \left(51 \theta ^3-189 \theta ^2-161 \theta -49\right)
-4 x^2 \theta  \left(2458 \theta ^2+2595 \theta +1376\right) \\
+2 x^3 \left(17008 \theta ^3+48924 \theta ^2+58107 \theta +25722\right) \\
+x^4 \left(-33384 \theta ^3-147204 \theta ^2-224558 \theta -115513\right) \\
+6 x^5 \left(3184 \theta ^3+17508 \theta ^2+39283 \theta +32662\right) \\
-4 x^6 \left(19814 \theta ^3+148059 \theta ^2+415393 \theta +423138\right) \\
+8 x^7 \left(2861 \theta ^3+26589 \theta ^2+77179 \theta +69798\right) \\
+4 x^8 \left(15728 \theta ^3+167496 \theta ^2+608482 \theta +752619\right) \\
+96 x^9 (\theta +4) \left(868 \theta ^2+7007 \theta +14397\right) \\
+48 x^{10} (\theta +4) (\theta +5) (1006 \theta +4443) \\
-33984 x^{11} (\theta +4) (\theta +5) (\theta +6)
\end{matrix}
$
& $\begin{matrix} 1\; \\ 1^\dagger \end{matrix}$ 
\tabularnewline
\hline 
${}{110A}$ 
&  
$\begin{matrix}
%pf[110 A]=
13 \theta ^3
+x \left(523 \theta ^3-624 \theta ^2-546 \theta -169\right) \\
-3 x^2 \theta  \left(7988 \theta ^2+8940 \theta +5021\right) \\
+x^3 \left(385606 \theta ^3+1076049 \theta ^2+1325658 \theta +592731\right) \\
+x^4 \left(-3568186 \theta ^3-15575109 \theta ^2-27075910 \theta -17328220\right) \\
+x^5 \left(21839700 \theta ^3+129002166 \theta ^2+286042517 \theta +229554485\right) \\
-3 x^6 \left(31296937 \theta ^3+232650999 \theta ^2+624852043 \theta +595871452\right) \\
+x^7 \left(290710695 \theta ^3+2602602891 \theta ^2+8194990757 \theta +9003786049\right) \\
-4 x^8 \left(162439196 \theta ^3+1700194029 \theta ^2+6135245796 \theta +7608504938\right) \\
+2 x^9 \left(514091078 \theta ^3+6157604979 \theta ^2+25040250991 \theta +34538728740\right) \\
-2 x^{10} (2 \theta +9) \left(274135057 \theta ^2+2463526113 \theta +5639381396\right) \\
+1320 x^{11} (2 \theta +9) (2 \theta +11) (134044 \theta +669401) \\
-26029080 x^{12} (2 \theta +9) (2 \theta +11) (2 \theta +13)
\end{matrix}
$
& $\begin{matrix} 1 \\ 3 \end{matrix}$ 
\tabularnewline
\hline 
${}{119AB}$ 
&  
$ \begin{matrix}
%pf[119 AB]=
11 \theta ^3
+x \left(383 \theta ^3-429 \theta ^2-385 \theta -121\right)
-x^2 \theta  \left(13907 \theta ^2+15630 \theta +9082\right) \\
+x^3 \left(185695 \theta ^3+517518 \theta ^2+653456 \theta +298236\right) \\
+x^4 \left(-1466322 \theta ^3-6391539 \theta ^2-11366471 \theta -7474399\right) \\
+3 x^5 \left(2621124 \theta ^3+15463800 \theta ^2+35023680 \theta +28962833\right) \\
+x^6 \left(-30481337 \theta ^3-226356705 \theta ^2-620372371 \theta -610718664\right) \\
+x^7 \left(88108207 \theta ^3+788093196 \theta ^2+2531016954 \theta +2875397306\right) \\
+x^8 \left(-192542015 \theta ^3-2013670887 \theta ^2-7410595865 \theta -9519098329\right) \\
+3 x^9 \left(106194863 \theta ^3+1271037894 \theta ^2+5271993246 \theta +7545734338\right) \\
+x^{10} \left(-395463473 \theta ^3-5329767801 \theta ^2-24598308230 \theta -38796473612\right) \\
+3 x^{11} \left(120203819 \theta ^3+1801129755 \theta ^2+9155073262 \theta +15771861160\right) \\
+x^{12} \left(-232808411 \theta ^3-3838864788 \theta ^2-21309215317 \theta -39808006260\right) \\
+3 x^{13} (\theta +6) \left(33297643 \theta ^2+399363585 \theta +1210021238\right) \\
-119 x^{14} (\theta +6) (\theta +7) (212357 \theta +1379743)
+2786385 x^{15} (\theta +6) (\theta +7) (\theta +8)
\end{matrix}
$
& $\begin{matrix} 1 \\ 2 \end{matrix}$ 
\tabularnewline
\hline 
\end{tabular}

\end{document}